\documentclass{amsart}

\usepackage{amsmath}
\usepackage{amsxtra}
\usepackage{amscd}
\usepackage{amsthm}
\usepackage{amsfonts}
\usepackage{amssymb}
\usepackage{eucal}

\newtheorem{cor}{Corollary}
\newtheorem{lem}{Lemma}
\newtheorem{prop}{Proposition}

\newtheorem{thm}{Theorem}

\theoremstyle{remark}
\newtheorem{rem}{Remark}

\numberwithin{thm}{section}
\numberwithin{cor}{section}
\numberwithin{propconstr}{section}
\numberwithin{prop}{section}
\numberwithin{lem}{section}
\numberwithin{rem}{section}

\newcommand{\thmref}[1]{Theorem~\ref{#1}}
\newcommand{\secref}[1]{Section~\ref{#1}}
\newcommand{\lemref}[1]{Lemma~\ref{#1}}
\newcommand{\propref}[1]{Proposition~\ref{#1}}
\newcommand{\corref}[1]{Corollary~\ref{#1}}

\newcommand{\remref}[1]{Remark~\ref{#1}}

\newcommand{\nc}{\newcommand}

\nc{\ssec}{\section}
\nc{\sssec}{\subsection}
\nc{\bi}{\bibitem}

\nc{\on}{\operatorname}

\nc{\ZZ}{{\mathbb Z}}
\nc{\NN}{{\mathbb N}}
\nc{\CC}{{\mathbb C}}
\nc{\DD}{{\mathbb D}}

\nc{\A}{{\mathcal A}}
\nc{\B}{{\mathcal B}}
\renewcommand{\O}{{\mathcal O}}
\nc{\K}{{\mathcal K}}
\nc{\M}{{\mathcal M}}
\nc{\N}{{\mathcal N}}
\nc{\V}{{\mathcal V}}
\nc{\R}{{\mathcal R}}
\nc{\T}{{\mathcal T}}
\nc{\I}{{\mathcal I}}
\nc{\J}{{\mathcal J}}
\renewcommand{\P}{{\mathcal P}}
\nc{\U}{{\mathcal U}}
\renewcommand{\L}{{\mathcal L}}

\nc{\Tf}{{\mathfrak T}}
\nc{\Gr}{{\mathfrak G}{\mathfrak r }}
\nc{\tGr}{\widetilde\Gr}
\nc{\Hf}{{\mathfrak H}}
\nc{\Mf}{{\mathfrak M}}
\nc{\Xf}{{\mathfrak X}}
\nc{\gf}{{\mathfrak v}}
\nc{\af}{{\mathfrak a}}
\nc{\gff}{{\mathfrak g}}
\nc{\hf}{{\mathfrak h}}
\nc{\z}{{\mathfrak z}}

\nc{\conf}{H_{\nabla}}
\nc\Vir{\on{Vir}}
\renewcommand\d{{\partial}}
\nc{\one}{{\mathbf 1}}
\nc{\unit}{{\mathbf u}{\mathbf n}{\mathbf i}{\mathbf t}}
\nc{\surj}{{\twoheadrightarrow}}
\nc\Cliff{\on{Cliff}}
\nc\Aut{\on{Aut}}
\nc\Act{\on{Act}}
\nc\tAct{\widetilde{\on{Act}}}
\nc\Maps{\on{Maps}}
\nc\tMaps{\widetilde{\Maps}}
\nc\Bun{\on{Bun}}
\nc\aut{\on{aut}}
\nc\can{\on{can}}

\nc\Spec{\on{Spec}}
\nc\Em{\on{E-m}}

\nc\tilgf{{\widetilde\gf}}
\nc\tilL{{\widetilde\L}}
\nc\tily{{\widetilde y}}
\nc\tilz{{\widetilde z}}
\nc\tilf{{\widetilde f}}
\nc\tilrho{{\widetilde\rho}}
\nc\tilphi{{\widetilde\phi}}

\title{Notes on 2d conformal field theory and string theory} 

\author{Dennis Gaitsgory}

\address{Department of Mathematics, Harvard University, Cambridge, MA
02138, USA\footnote{Current address for 1998-1999: School of Math., IAS,
Olden Lane, Princeton, 08540 NJ, USA; email: gaitsgde@math.harvard.edu}}

\date{May 1999.   This paper will appear in
{\em Quantum Fields and Strings: A course for Mathematicians} (P.~Deligne, P.~Etingof, D.S.~Freed, L.~Jeffrey, D.~Kazhdan,
J.~Morgan, D.R.~Morrison and E.~Witten eds.), Amer. Math. Soc., Providence (1999).}

\begin{document}

\maketitle

\centerline{\bf \Large Introduction}
\ssec{Contents of these notes}
\sssec{}
 These notes appeared out of an attempt to write down explanations and solutions 
of home assignments for the course
on String Theory which was given by E.~D'Hoker in Spring 1997 at IAS.

Our original plan was to deal only with questions that admitted a rigorous 
mathematical formulation or axiomatization.
A variety of such questions arises in the most basic ingredient of 
String Theory, namely in $2$-dimensional 
conformal field theory (CFT). 

We shall restrict our attention to a tiny part of 2D CFT: our goal is to explain 
the mathematical 
formalism that describes the operator product expansion (OPE) of quantum fields. 
To make things as simple as possible 
(and since we wish to stay within the realm of algebraic geometry) we will deal 
only with fields from the so-called 
holomorphic sector. 

Fortunately, appropriate mathematical objects that imitate the OPE operation on 
quantum fields have been found
several years ago by A.~Beilinson and V.~Drinfeld (they called them ``chiral 
algebras''). We shall adopt their point
of view completely and the present paper (and especially Chapters 
1, 3 and 4) can be regarded as an example-oriented digest of the
long-awaited text of Beilinson and Drinfeld.

The language of chiral algebras is in many ways equivalent to the language of 
vertex operator algebras (VOA's). More 
precisely,
a vertex operator algebra is the same as a ``universal chiral algebra'', 
i.e. one 
that is defined naturally over
an arbitrary Riemann surface. Since we will be interested in chiral algebras 
that 
are universal in the above sense, our
treatment of such matters as the Sugawara construction 
(Sect. 9), fermion algebras and boson-fermion correspondence 
(Sect. 11), the BRST complex (Sect. 12) and of many others, is nothing 
more than 
a reinterpretation 
in terms of algebraic geometry of constructions that have been known to the 
founders of the theory of 
vertex operator algebras for many years (cf. \cite{FF},\cite{FF1},\cite{FLM},
\cite{FHL},\cite{FZh}).

Therefore, we wish to emphasize that hardly any 
part of the contents of these notes is the invention of their author. 
Besides the classical sources mentioned just above we have borrowed the material 
from the partially written 
manuscripts \cite{BD} and \cite{BFM}.

\sssec{}
 Let us discuss in some detail the mathematical contents of the paper.

\medskip

Chapter 1 deals with the basics of the theory of chiral algebras.

\smallskip

In Section 3 we introduce the main character, namely the notion of a chiral 
algebra. A chiral algebra on a curve $X$
is a D-module with a binary operation that satisfies the axioms of a Lie-* 
bracket. 
One of the features of chiral algebras
is that these objects are not so easy to construct. The more elementary 
constructions go through the procedure of taking
the chiral universal enveloping algebra of a more accessible object-- a Lie-* 
algebra (the latter notion is introduced
in Section 4). Finally, in Section 5 we introduce the notion of the space of 
conformal blocks of a chiral algebra (the
dual to the space of conformal blocks should be thought of as a set of possible 
correlation functions of a quantum field
theory (QFT) with a given OPE of the fields).

\medskip

In Chapter 2 we make an attempt to define axiomatically what a 2D CFT is. 
However, 
it should be clear form the very
beginning that our list of axioms is far from being
complete: we are dealing with the short-distance singularities of quantum
fields only and we do not mention at all what happens when a Riemann surface 
approaches the boundary of the moduli space.
We hope to be able to address this more subtle and interesting issue in the 
future.

\smallskip

In Section 6 we introduce the notion of a local $\O$-module over the universal 
curve $\Xf$ over the moduli space of smooth 
curves.
A local $\O$-module is a quasi-coherent sheaf which is ``universal'' 
(in the same 
sense as a VOA is a ``universal'' chiral
algebra). We need this notion since fields of a CFT are clearly local in 
the above 
sense, i.e. the fields and
the OPE structure on them over a given region of a Riemann surface do not 
``feel'' 
a change of a metric at a remote point.
In Section 7 we give the axioms of a CFT data in the simplest case of central 
charge $0$ and in Section 8 we generalize
the discussion to the case of an arbitrary central charge. 
In particular, we show 
(what should be a standard thing by now)
that the space of conformal blocks of a CFT chiral algebra acquires a 
(projective) 
connection along the moduli space 
of smooth curves.

\medskip

Chapter 3 is devoted to the examples. (Unfortunately, all our examples are 
trivial 
or almost trivial from the point
of view of physics. They correspond to either non-existing physical 
theories or to 
theories that are free.)

We discuss three types of examples here: Kac-Moody chiral algebras, i.e. those 
associated to Lie algebras with an invariant
quadratic form, in particular, the Heisenberg chiral algebra which 
corresponds to 
a Lie algebra with a zero bracket 
(all this in Section 9); the dilaton twist of the Heisenberg chiral algebra 
(Section 10) and the fermionic chiral
algebra, which is referred to as the $bc$-system in the physics literature 
(Section 11). In all these examples we write
down explicitly the energy-momentum tensor (or equivalently, the Sugawara 
construction) and perform explicitly
the corresponding calculations. We should mention that there are many more ways
to make the same calculations and ours
are not always the optimal ones: we just wanted to exhibit as many techniques 
as possible of manipulating with
chiral algebras. Also, the reader should notice that the way we present these 
calculations is very close to how a physicist
would have done it.

\medskip

In Chapter 4 we discuss the BRST complex and related matters.

In Section 12, given a Lie-* algebra $\B$, we construct (following Beilinson) 
its central extension $\B'$ and define a
BRST differential on a chiral algebra obtained by tensoring the twisted chiral 
universal enveloping algebra of $\B'$
and the corresponding fermionic algebra. This enables one to write a BRST 
reduction of any chiral algebra endowed
with an action of $\B$ (as usual in such situations, this action must be of 
``Harish-Chandra'' type). In Section 13 we apply
the discussion of Sect. 12 in the case of the BRST reduction of a chiral 
algebra corresponding to a CFT of central charge $26$
with respect to the Virasoro algebra. This procedure is used by physicists for 
what can be mathematically reformulated
as writing down ``explicitly'' the connection along $\Mf$ on the space of 
conformal blocks. However, we do not pretend to
having really understood why the physicists really need the BRST procedure.

\medskip

Starting from Chapter 5, this paper becomes much more involved mathematically 
and sometimes we are forced to omit proofs
of certain important statements. Our goals in Chapters 5 and 6 are the 
construction of the bosonic chiral algebra 
(in the VOA terminology it is sometimes called the lattice VOA) and the proof 
of \thmref{confbose} that computes its space of 
conformal blocks on a given curve in terms of the canonical line bundle over 
the Jacobian 
(or, rather a direct product of several copies of the Jacobian) of this curve. 
Usually,
in the theory of VOA's, the sought-for bosonic chiral algebra is 
constructed by an explicit formula (cf. \cite{FLM})
(which is the most classical example of a ``vertex operator'') and the 
fact that this formula works remains for us a miracle.  

\smallskip

The construction of the bosonic chiral algebra described in this paper is quite 
different: it is very geometric and uses 
a sophisticated algebro-geometric object called the Beilinson-Drinfeld 
Grassmannian (cf. \cite{MiVi}, \cite{FM}). In fact,
our exposition covers a part of a baby (=abelian) version of a new 
construction of the geometric Langlands 
correspondence due again to Beilinson and Drinfeld.

\smallskip

In Section 14 we describe a new way of thinking about chiral algebras (due to 
Beilinson): now a chiral algebra is interpreted as a 
system of $\O$-modules on various powers of our curve endowed with what we 
shall call ``factorization
isomorphisms''. In Section 16 (after some geometric preliminaries in 
Section 15) we demonstrate in a relatively simple situation
an interplay between the two pictures.

\smallskip

In Section 17 we describe the canonical line bundle over the B.-D. Grassmannian 
corresponding to a torus, which arises from
an even quadratic form on the lattice of $1$-parameter subgroups of this torus. 
Finally, in Section 18 we 
introduce the bosonic chiral algebra and prove \thmref{confbose}.

\ssec{Some background on D-modules} \label{Dmod}
\sssec{}
 First, some notation:

For a scheme (or a stack) $Z$, $\O_Z$ (resp., $\Theta_Z$, $\Omega_Z$, $D_Z$) 
will 
denote the structure sheaf
(resp., the tangent sheaf, the sheaf of $1$-forms, the sheaf of 
differential operators) on $Z$. 
The terms ``an $\O$-module on $Z$'' and ``a quasi-coherent sheaf on $Z$'' will 
mean for us the same thing 
(with a few exceptions in Sections 6-7, where
we will be dealing with projective limits of coherent sheaves). 
If $f:Z_1\to Z_2$ 
is a morphism of algebraic
varieties, we shall denote by $f_*$ and $f^*$ the direct and the inverse image 
functors respectively,
in the category of quasi-coherent sheaves. 
The functor of the inverse image in the 
category of sheaves
(and not $\O$-modules) will be used only once 
(in \secref{centralalgebroid}) and 
will be denoted by 
$f^{\bullet}$.

\smallskip

For an algebraic curve $X$, $\Delta$ will denote the closed embedding of the 
diagonal divisor in $X^n$ and 
$j:X^n\setminus\Delta\hookrightarrow X^n$ will denote the open embedding of its 
complement. A symbol 
like $\Delta_{x_i=x_j}$
will designate an embedding of the closed subset of the diagonal divisor 
corresponding to points 
$(x_1,\dotsc, x_n)\in X^n$
with $x_i=x_j$. By $p_i:X^n\to X$ we shall denote the 
projection on the $i$-th factor.

\smallskip

$\Mf$ will denote the moduli stack of smooth complete curves and 
$\pi:\Xf\to\Mf$ 
will denote
the universal curve over $\Mf$. $\Xf^n$ will denote the $n$-fold 
Cartesian product 
of $\Xf$
with itself over $\Mf$; $\Delta$ will as above denote the embedding of the 
divisor 
of diagonals into $\Xf^n$ 
and $p_i:\Xf^n\to\Xf$ will stand for the corresponding projections. 
If $L_1$ and 
$L_2$ are quasi-coherent 
sheaves on
$\Xf^{k_1}$ and $\Xf^{k_2}$ respectively we let $L_1\boxtimes L_2$ denote (by 
slightly 
abusing the notation) the sheaf
on $\Xf^{k_1+k_2}$ obtained by restricting the corresponding sheaf on 
$\Xf^{k_1}\times \Xf^{k_2}$.

\sssec{}
 In this paper we will be extensively using the language of
D-modules. Unfortunately, there are very few sources with an appropriate
exposition of the theory, the standard reference being the book
\cite{Bo}.\footnote{One of the best readings on D-modules is an unpublished
series of lectures by J.~Bernstein (\cite{Ber}) (obtainable by e-mail from
the author of these notes).}

\smallskip

For a smooth variety $Z$ a (right) D-module is a (right) module over the sheaf 
$D_Z$ of differential 
operators on $Z$. By definition, 
there exists a natural forgetful functor {D-modules} $\to$ {$\O$-modules}. 

When $Z$ is a singular scheme, one proceeds as follows: we embed $Z$
(locally) as a closed subscheme into a non-singular variety $Z'$ and consider
D-modules on $Z'$ that are set-theoretically supported on $Z$. (The
correctness of this definition follows from the well-known Kashiwara
theorem). In this case, there still is a forgetful functor {D-modules} $\to$
{$\O$-modules} which sends a D-module $\M$ over $Z'$ as above to the
$\O_Z'$-module of its sections that are scheme-theoretically supported on
$Z$.  This enables one to make sense of right D-modules over a strict
ind-scheme that are set-theoretically supported on a closed
finite-dimensional subscheme (cf. Section~ 6).

\medskip

For a smooth variety $Z$ there exists a functor 
$$h:\,\{\text{ D-modules on }Z\}\, \to\, \{\text{ Sheaves in
the \'etale topology on }Z\}$$ defined as follows: 
$$\M\to \M/\M\cdot \Theta_Z$$
(we shall often abuse the notation and denote by $h$ also the projection $\M\to 
h(\M)$.)

The functor $h$ is right-exact and we let $Lh$ denote the corresponding left 
derived functor.
(the functor $Lh$ applied to holonomic D-modules with regular singularities 
realizes the celebrated 
Riemann-Hilbert correspondence).
The functor 
$$R\Gamma\circ Lh:\,\{\text{ D-modules }\} \to\, \{\text{ Vector spaces }\}$$ 
is called the De Rham cohomology
functor, denoted $\M\to DR^{\bullet}(Z,\M)$.

When $f:Z_1\to Z_2$ is a map between smooth algebraic varieties, one 
defines the 
direct image functor 
$$f_*:\{\text{ D-modules on $Z_1$}\}\to \{\text{ D-modules on $Z_2$}\}$$ 
(to distinguish it from the corresponding functor
on the category of $\O$-modules; in cases where confusion can occur, we shall 
specify explicitly 
which one we are using.) When $f$ is a projection $Z_1=Z_2\times Y\to Z_2$ the 
definition of $f_*$
is a straightforward generalization of that of $DR^{\bullet}$ (the latter 
corresponds to the case
$Z_2=\on{pt}$). When $f$ is a closed embedding, $f_*$ is defined using 
Kashiwara's 
theorem (see above);
in this case we shall replace the notation $f_*$ by $f_{!}$ in order to avoid 
confusion with the functor
$f_*$ for $\O$-modules.

For instance, when $f$ is affine, the functor $f_*$ is right-exact and when 
$f$ is 
an open embedding,
$f_*$ is left exact.

For ($f,Z_1,Z_2$) as above one defines also the inverse image functor 
$$f^{!}:\,\{\text{ D-modules on }Z_2\}\, \to \{\text{ D-modules on }Z_1\}.$$ 
We refer the reader to \cite{Bo} for the definition.

\medskip

When $Z$ is a smooth variety, in addition to right D-modules one can 
consider left 
D-modules (it may seem surprising, but one
cannot define left D-modules for non-smooth schemes). The categories of 
left and 
right D-modules are equivalent and the functor
that realizes this equivalence sends a left D-module $^l\M$ to 
$\M:={}^l\M\otimes\Omega^{top}_Z$;in particular, the basic left D-module
$\O_Z$ goes over to the basic right D-module $\Omega^{top}_Z$. 

If $^l\M$ is a left D-module on $Z$ and $\N$ is another left
(resp., right) D-module, the tensor product $^l\M\otimes\N$
is naturally a left (resp., right) D-module.
If now $\M_1$ and
$\M_2$ are two right D-modules, we shall denote by $\M_1\overset{!}\otimes\M_2$ 
the right D-module
$$({}^l\M_1\otimes {}^l\M_2)\otimes\Omega_Z^{top}.$$
The operation $\overset{!}\otimes$ makes the category of right D-modules over a 
smooth variety into a symmetric tensor category.
For instance, one can talk about Lie algebras in the category of right 
D-modules 
and of their universal enveloping algebras
(this will be important in \secref{unenv}).

\medskip

In addition, in this paper we will be dealing with D-modules on algebraic 
stacks in the smooth topology
($\Mf$ is an example of such). It follows from the definition of a stack that 
these objects are perfectly
suitable for having D-modules on them. For example, if a stack has 
the form $Z/G$ , where $Z$ is a scheme and $G$ is a group
acting on $Z$, a D-module on $Z/G$ is the same as a strongly $G$-equivariant 
D-module on $Z$ (cf. \cite{BB}, Section 1 for
the discussion of equivariant D-modules.)

\medskip

The following observation will be used several times in this paper:

\smallskip

For a smooth variety $Z$ consider the co-induction functor
$$\{\text{ $\O$-modules on $Z$}\}\,\to\{\text{ D-modules on $Z$}\}$$
which sends an $\O$-module $\L$ to $\L\underset{\O_Z}\otimes D_Z$.

\begin{lem}
We have:

\smallskip

\noindent{\em (a)}
For a closed embedding $f:Z_1\to Z_2$ and an $\O$-module $\L$ over $Z_1$ 
there is 
a canonical isomorphism:
$$f_{!}(\L\underset{\O_{Z_1}}\otimes D_{Z_1})\simeq 
f_*(\L)\underset{\O_{Z_2}}\otimes D_{Z_2}.$$

\smallskip

\noindent{\em (b)}
For two $\O$-modules $\L_1$ and $\L_2$ over $Z$ the space 
$$\on{Hom}_{\O-\on{mod}}(\L_1,\L_2\underset{\O_Z}\otimes D_Z)\simeq 
\on{Hom}_{\on{D}-\on{mod}}
(\L_1\underset{\O_Z}\otimes D_Z,\L_2\underset{\O_Z}\otimes D_Z)$$ is 
canonically 
isomorphic to the space
of differential operators from $\L_1$ to $\L_2$.
\end{lem}

\sssec{}
 Another notion that will be important for us is that of a Lie
algebroid.

\smallskip

A Lie algebroid $\af$ over a smooth variety $Z$ is a quasi-coherent sheaf $\af$ 
over $Z$ with the following 
additional structure:

\begin{itemize}

\item
An $\O$-module map $\phi:\af\to\Theta_Z$.

\item
A bi-differential operator $[\,\,,\,\,]:
\af\underset{\CC}\otimes\af\to\af$ satisfying the axioms of a Lie-* bracket.

\end{itemize}

The above structures should be compatible in the following way: if $a_1$ and 
$a_2$ are sections of $\af$ and $f$ is a
a function over $Z$, we must have:
$$[f\cdot a_1,a_2]=f\cdot [a_1,a_2]+a_1\cdot \on{Lie}_{\phi(a_2)}(f).$$

The most basic example of a Lie algebroid is $\af=\Theta_Z$.

\smallskip

In a natural way one defines a notion of a (left) module 
over a Lie algebroid. In particular, if $\af$ is as above,
any left D-module on $Z$ is automatically an $\af$-module.

\smallskip

Let now $f:Z_1\to Z_2$ be a map between smooth varieties and let $\af$ be a Lie 
algebroid on $Z_2$. Consider
the $\O$-module $\af'$ over $Z_1$ given by
$$\on{ker}[f^*(\af)\oplus\Theta_{Z_1}\to f^*(\Theta_{Z_2})],$$
where the maps $f^*(\af)\to f^*(\Theta_{Z_2})$ and $\Theta_{Z_1}\to 
f^*(\Theta_{Z_2})$ are $f^*(\phi)$ and $df$,
respectively.

\begin{lem}
The $\O$-module $\af'$ on $Z_1$ has a canonical
structure of Lie algebroid. Moreover, the functor $f^*$ maps naturally left 
$\af$-modules to left $\af'$-modules.
\end{lem}

We shall call $\af'$ the pull-back of $\af$ with respect to $f$.

\smallskip

There are several other notions and constructions connected to Lie algebroids 
that will be used in this paper
(such as Picard Lie algebroids and their connection with rings of twisted 
differential operators;
Lie algebroids in the equivariant setting, in particular, Harish-Chandra Lie 
algebroids). They all are discussed
in a clear and self-contained way in the first two sections of \cite{BB}.

\bigskip

\noindent{\bf Acknowledgments}
While working on these notes, the author benefited from numerous 
discussions with 
many mathematicians and physicists
(the complete list would be too long to include here). 
Special thanks are to the 
instructors and participants of
the QFT program at IAS and most of all, to E.~D'Hoker. Quite separately, 
I wish to express my gratitude to A.~Beilinson, 
A.~Braverman and D.~Kazhdan, without whom this paper would have never been 
written. All I know about chiral algebras I was 
taught by A.~Beilinson and this paper is largely a write-up of what I 
heard from him on various occasions. Initially,
the project of writing these notes was meant to be a joint one 
with A.~Braverman 
and the fact of his collaboration was very important in the early stages of 
the work. Most of what I managed to understand at the QFT program was via 
discussions with D.~Kazhdan; in addition, his constant attention and 
support were crucial factors in the course of preparation of
these notes.

\newpage

\centerline{\bf Chapter I. Chiral Algebras}

\medskip

 In this chapter we will be working over a fixed smooth curve $X$ (in 
\secref{confblocks} we will assume, 
in addition, that $X$ is complete). The symbols $\O$ 
(resp., $\Omega$, $\Theta$) with an omitted subscript will 
mean the corresponding objects for $X$. However, all the results of this 
section are valid in the case when $X$ is 
in fact a family of curves of the above type over an arbitrary base; in 
particular, $X$ can be taken to be the universal
curve $\Xf$ over the moduli stack $\Mf$ of smooth complete curves.

\ssec{Definition of chiral algebras}
\sssec{} \label{defnchiralalg}
 Let $^l\A$ be a (left) D-module on $X$, which we shall think of as 
consisting of ``fields''
of our CFT. Consider the corresponding right D-module $\A={}^l\A\otimes\Omega$; 
its sections should 
be viewed as fields with values in $1$-forms, i.e. as currents. 

\smallskip

A chiral algebra structure on $\A$ is a D-module map 
(called ``chiral bracket''):
$$j_*j^*(\A\boxtimes\A){\overset{\{,\}}\longrightarrow}\Delta_{!}(\A),$$
that satisfies the following two conditions:

\begin{itemize}

\item Antisymmetry:

\smallskip

\noindent If $f(x,y)\cdot a\boxtimes b$ is a section of $j_*j^*(\A\boxtimes\A)$, 
then
$$\{f(x,y)\cdot a\boxtimes b\}=-\sigma_{1,2}(\{f(y,x)\cdot b\boxtimes a\}),$$ 
where 
$\sigma_{1,2}$ is the transposition acting on $\Delta_{!}(\A)$.

\item Jacobi identity:

\smallskip

\noindent If $a\boxtimes b\boxtimes c\cdot f(x,y,z)$ is a section of the
restriction of $\A\boxtimes \A\boxtimes\A$ to the complement of the divisor of 
diagonals in 
$X\times X\times X$, then the element  
$$\{\{f(x,y,z)\cdot a\boxtimes b\}\boxtimes c\}+
\sigma_{1,2,3}(\{\{f(z,x,y)\cdot b\boxtimes c\}\boxtimes a\})+ 
\sigma_{1,2,3}(\{\{f(y,z,x)\cdot c\boxtimes a\}\boxtimes b\})$$
of $\Delta_{x=y=z}{}_{!}(\A)$ vanishes.
(Here $\sigma_{1,2,3}$ denotes the lift of the cyclic automorphism of $X\times 
X\times X$:
$(x,y,z)\to (y,z,x)$ to $\Delta_{x=y=z}{}_{!}(\A)$.)

\end{itemize} 

\medskip

A basic example of a chiral algebra is $\A=\Omega$. The bracket operation is 
given by the
canonical map $\can_\Omega:j_*j^*(\Omega\boxtimes\Omega)\to\Delta_{!}(\Omega)$ 
(see below).

\medskip

In most examples we shall deal with chiral algebras that possess a unit: a 
unit in a chiral
algebra $\A$ is a map $\unit:\Omega\to \A$ such that the square below becomes 
commutative.
$$
\CD
j_*j^*(\Omega\boxtimes \A) @>{\unit\boxtimes\on{id}}>> j_*j^*(\A\boxtimes \A) \\
@V{\can_\A}VV             @V{\{,\}}VV   \\
\Delta_{!}(\A)  @>{\on{id}}>>  \Delta_{!}(\A),   
\endCD
$$
where the left vertical arrow is the canonical map that comes from the short 
exact 
sequence
$$0\to \Omega\boxtimes\A\to j_*j^*(\Omega\boxtimes \A)\to\Delta_{!}(\A)\to 0.$$

\medskip

Along with ordinary chiral algebras one can consider the corresponding
super-objects. In this case, the above antisymmetry and Jacobi identity 
conditions
get transformed according to the sign rules of the super world. 

\sssec{}  \label{defnchiralmodules}
 Now let $\A$ be a chiral algebra and let $\M$ be a (right) D-module on $X$. 

\smallskip

A chiral $\A$-module structure on $\M$ is a D-module map 
$$\rho:j_*j^*(\A\boxtimes \M)\to \Delta_{!}(\M)$$ 
such that if $f(x,y,z)\cdot a\boxtimes b\boxtimes m$ is a section of the 
restriction of
$\A\boxtimes\A\boxtimes\M$ to the complement of the divisor of diagonals in 
$X\times X\times X$,
we have:
$$\rho(\{f(x,y,z)\cdot a\boxtimes b\},m)=
\rho(a,\rho(f(x,y,z)\cdot b,m))-\sigma_{1,2}(\rho(b,\rho(f(y,x,z)\cdot a,m)),$$
as sections of $\Delta_{x=y=z}{}_{!}(\M)$.

When $\A$ has a unit we require moreover that
the induced map $j_*j^*(\Omega\boxtimes \M)\to \Delta_{!}(\M)$ 
coincide with the map $\can_\M$.

\smallskip

It is easy to check that the chiral bracket makes $\A$ into a chiral module over 
itself.

\sssec{}  \label{modes}
 Let us now explain the concept of the ``mode expansion'' of a field. 

\medskip

Let $x$ be a point of $X$ and let $\M$ be a chiral $\A$-module supported at $x$; 
we shall
denote by $M$ the underlying vector space (i.e. $M\simeq h(\M)\simeq DR^0(\M)$).

\smallskip

\noindent A typical example of such an $\A$-module is 
$\M:=i_x{}_{!}i_x^{!}(\A)[1]$, called the vacuum module at $x$, in this case 
$A_x\simeq i_x^{!}(\A)[1]$.

\medskip

Let $\z$ be a coordinate near $x$, i.e. $\z$
is a regular function on some $U$ with $\z(x)=0,\,\,d\z(x)\neq 0$. Then any 
section
$a\in\Gamma(U,\A)$ gives rise to a collection $a_n,\,n\in \ZZ$, of 
elements in $\on{End}(M)$:

\smallskip

\noindent 
By passing to the De Rham cohomology, the structure map $\rho$ gives rise 
to a map
$DR^0(U\setminus x,\A)\otimes M\to M$ and we take $a_n$ to be the endomorphism 
of $M$ corresponding
to the image of $a\cdot \z^n$ under the projection
$$\Gamma(U\setminus x,\A)\to DR^0(U\setminus x,\A).$$

\begin{rem}

We shall see shortly that for any $(x\in U)$ as above,
$DR^0(U\setminus x,\A)$ has a natural structure of Lie algebra. 
The vector space $A_x$ together with the action of
 $DR^0(U\setminus x,\A)$ on it is a prototype of the Hilbert 
space of our CFT.

\end{rem}

\begin{lem} \label{finite}
For any $m'\in M$ and $a$ as above we have 
$$a_n\cdot m'=0 \text{ for } n>>0.$$
\end{lem}

\bigskip

Now let $a$ and $b$ be two sections of $\A$ over $U$; we can produce a third 
section $c\in\Gamma(U,\A)$
by setting
$$c=(\on{id}\boxtimes h)(\{\frac{1}{\z_1-\z_2}\cdot a\boxtimes b\})\in 
(h\boxtimes\on{id})(\Delta_{!}(\A))\simeq\A,$$
where $\z_1=\z\boxtimes 1$ and $\z_2=1\boxtimes \z$ are the corresponding 
functions on $X\times X$.

\begin{prop}  \label{normalordering}
Modes of $c$ can be expressed through the modes of $a$ and $b$ by the ``normal 
ordering'' formula: 
$$c_k=:a\cdot b:_k={\underset{n\geq 0}\Sigma} a_{k-1-n}\cdot b_{n}
+{\underset{n\geq 0}\Sigma}b_{-1-n}\cdot a_{k+n}.$$
\end{prop}

Note that both sides of the above formula depend on the choice of the coordinate 
~$\z$.

\begin{proof}

Choose a section $m\in\M$ which goes over to $m'\in M$ under the identification
$DR^0(X,\M)\simeq M$.

\smallskip

The element $c_k\cdot m'$ is the image of
$$\rho(\{\frac{z_1^k\boxtimes 1}{\z_1-\z_2}\cdot a\boxtimes b\},m)\in 
\Delta_{!}(\M)$$
under the identification $DR^0(X\times X\times X,\Delta_{1,2,3}{}_{!}(\M))\simeq 
DR^0(X,\M)\simeq M$.

By definition, the above expression can be rewritten as
$$\rho(\z_1^{k-1}\cdot a,\rho(\frac{1}{1-\z_2/\z_1}\cdot b,m))+
\rho(\z_2^k\cdot \z_1^{-1}\cdot b,\rho(\frac{1}{1-\z_2/\z_1}\cdot a,m)).$$
(In this formula we have omitted $\sigma_{1,2}$ since it does not affect the image 
of our element
in De Rham cohomology.)

Now, \lemref{finite} implies that we can expand $\frac{1}{1-\z_2/\z_1}$ into a 
power series
$1+\z_2/\z_1+\z_2^2/\z_1^2+\ldots $ and our assertion follows.

\end{proof}

\ssec{Lie-* algebras and construction of chiral algebras}  
We shall now make a brief detour and discuss the auxiliary notion of a Lie-* 
algebra.
These are objects that do not appear naturally in QFT; however, they are useful 
(and even necessary)
for a construction of any non-trivial chiral algebra.

\sssec{}  \label{Liestaralg}
 A Lie-* algebra on $X$ is a (right) D-module $\B$ with a map (called ``Lie-* 
bracket''):
$$\B\boxtimes\B{\overset{\{,\}}\longrightarrow}\Delta_{!}(\B),$$
which is antisymmetric and satisfies the Jacobi identity in the sense similar to 
what we had
in the definition of chiral algebras.

\medskip

If $\B$ is a Lie-* algebra, it follows from the definition that $h(\B)$ is a sheaf
(in the \'etale topology) of ordinary Lie algebras; moreover it acts on $\B$ by 
endomorphisms
of the D-module structure that are derivations of the Lie-* structure.

\smallskip

In particular, for an affine subset $U\subset X$, $DR^0(U,\B)$ is a Lie algebra.

\sssec{}  \label{Liestarmodules}
 As for modules over a Lie-* algebra, there are two types of these:

\medskip

\noindent A Lie-* module over a Lie-* algebra $\B$ is a D-module $\M$ with a map
$$\rho:\B\boxtimes\M\to\Delta_{!}(\M)$$ such that for a section $a\boxtimes 
b\boxtimes m$
of $\B\boxtimes\B\boxtimes\M$ the two sections 
$$\rho(\{a\boxtimes b\},m) \text{ and 
}\rho(a,\rho(b,m))-\sigma_{1,2}(\rho(b,\rho(a,m)))$$
of $\Delta_{1,2,3}{}_{!}(\B)$ coincide.
 
\medskip

\noindent A chiral module over a Lie-* algebra is again a D-module $\M$, but with 
an operation
$$\rho:j_*j^*(\B\boxtimes\M)\to\Delta_{!}(\M)$$ such that every section
$$f(x,y,z)\cdot a\boxtimes b\boxtimes m\in\Gamma(X\times X\times X\setminus 
(\Delta_{x=z}\cup \Delta_{y=z}),
\B\boxtimes\B\boxtimes\M)$$ satisfies an identity similar to the one in the 
definition of chiral modules over a chiral algebra.

\medskip

For example, if $\M$ is a Lie-* (resp., chiral) module over a Lie-* algebra 
$\B$ which is supported at 
$x\in U\subset X$, the Lie algebra $DR^0(U,\B)$ 
(resp., $DR^0(U\setminus x,\B)$) 
acts on the underlying 
vector space $M=DR^0(U,\M)$. 

\medskip

We have an obvious forgetful functor from the category of chiral modules over a 
Lie-*
algebra $\B$ to the category of Lie-* modules.

\smallskip

In addition, there is a forgetful functor from the category
of chiral algebras to the category of Lie-* algebras. If $\M$ is a chiral 
module over
a chiral algebra $\A$ it is automatically a chiral module over the 
corresponding Lie-* algebra.

\sssec{}   \label{firstexamples}
 Let us now consider a few examples of Lie-* algebras.

\medskip

{\bf Example 1} 

Let $\gff$ be a Lie algebra with an ad-invariant quadratic form $Q$ on it. 
Consider the D-module
$\B(\gff,Q):=\gff\otimes D_X \oplus\Omega$. We define a Lie-* algebra structure
on $\B(\gff,Q)$ in the following way:

\smallskip

The Lie-* bracket $\{,\}$ factors through $\B(\gff,Q)\boxtimes\B(\gff,Q)\to 
\gff\otimes D_X\boxtimes \gff\otimes D_X$
and 
$$\{\xi_1\otimes 1\boxtimes \xi_2\otimes 1\}=[\xi_1,\xi_2]\otimes \one\oplus 
Q(\xi_1,\xi_2)\otimes \one',$$
where $\one$ is the unit section of $D_X\subset \Delta_{!}(D_X)$ and $\one'$ is 
the canonical antisymmetric section
of $\Omega\boxtimes\Omega(2\cdot\Delta)/\Omega\boxtimes\Omega$.

\smallskip

If $\Spec(\hat\O_x)$ (resp., $\Spec(\hat\K_x)$) is the formal disc 
(resp., formal punctured disc) around $x\in X$,
we can naturally identify the Lie algebras 
$DR^0(\Spec(\hat\O_x),\B(\gff,Q))$ and
$DR^0(\Spec(\hat\K_x)$,
$\B(\gff,Q))$ with
$\gff\otimes\hat\O_x$ and $\gff\otimes\hat\K_x\oplus\CC$. For instance, when 
$\gff$ is semisimple, 
$DR^0(\Spec(\hat\K_x),\B(\gff,Q))$
is the corresponding Kac-Moody algebra.

\medskip

{\bf Example 2} 

Let $\T$ denote the D-module $\Theta\otimes D_X$. The action of $\Theta$ on itself 
by Lie derivations
is a bi-differential operator $\Theta{\underset{\CC}\otimes}\Theta\to\Theta$; 
hence, according to \secref{Dmod},
it gives rise to a D-module map
$$\T\boxtimes\T\simeq \Theta\otimes D_X\boxtimes \Theta\otimes D_X\to 
\Theta{\underset{\O_{X\times X}}\otimes} 
D_{X\times X}\simeq \Delta_{!}(\T),$$
which we take to be minus\footnote{The minus sign is chosen in
order to be consistent with the physics literature.} the structure map for the 
Lie-* algebra $\T$.

\smallskip

It is easy to check that the resulting map satisfies the axioms of a Lie-* 
bracket. The Lie algebras
$DR^0(\Spec(\hat\O_x),\T)$ and $DR^0(\Spec(\hat\K_x),\T)$ identify with 
$\Vir^+:=\Theta\underset{\O}\otimes\hat\O_x$
and $\Vir:=\Theta\underset{\O}\otimes\hat\K_x$, respectively (the bracket being 
minus the usual one).

\medskip

{\bf Example 3}

Let $\Theta'$ be an extension $0\to \Omega\to\Theta'\to\Theta\to 0$. Assume
that the sheaf $\Theta$ of vector fields on $X$ acts 
on $\Theta'$ ``by Lie
derivations''\footnote{We say that $\Theta$ acts on an $\O$-module $\L$ by
Lie derivations if the action map $\Theta{\underset{\CC}\otimes}\L\to\L$ is a
bi-differential operator that satisfies $$\xi(f\cdot
l)=f\cdot\xi(l)+\on{Lie}_\xi(f)\cdot l,$$ for $\xi\in\Theta$, $l\in\Theta'$
and $f\in\O$.}  in a way compatible with the above filtration and with its
action on $\Omega$ and on $\Theta$.  We shall endow the D-module
$$\T':=\Theta'\otimes D_X/\on{ker}(\Omega\otimes D_X\to\Omega)$$
with a Lie-* algebra structure. As a D-module $\T'$ is an extension
$$0\to\Omega\to \T'\to\T\to 0$$
and the Lie-* algebra structure on it will be such that $\Omega$ is central (i.e. 
is annihilated by the Lie-* bracket)
and $\T'$ will be a central extension of the Lie-* algebra $\T$ considered in 
Ex.2.
 
\smallskip

According to \secref{Dmod}, the action of $\Theta$ on $\Theta'$ can be viewed as a 
map
$$\Theta\boxtimes\Theta'\to \Delta_*(\Theta'){\underset{\O_{X\times X}}\otimes} 
D_{X\times X}\simeq \Delta_{!}(\Theta'\otimes D_X).$$
Therefore, we obtain a D-module map
$$(\Theta\otimes D_X)\boxtimes (\Theta'\otimes D_X)\to \Delta_{!}(\Theta'\otimes 
D_X)\to
\Delta_{!}(\T').$$

\smallskip

It is easy to see that $\Theta\otimes D_X\boxtimes \Omega\otimes D_X$ lies
in the kernel of the above map, which means that it factors as $\T\boxtimes 
\T\to\Delta_{!}(\T')$.
We define now the Lie-* bracket on $\T'$ as minus the composition:
$$\T'\boxtimes\T'\to \T\boxtimes\T\to\Delta_{!}(\T').$$

\smallskip

We shall see in \secref{linalg} that the set of isomorphism classes of extensions 
$\Theta'$ as above
is in bijection with $\CC$. For any $\Theta'$ of this form, the 
Lie algebra 
$DR^0(\Spec(\hat\K_x),\T')$ is the corresponding Virasoro extension 
$$0\to \CC\to \Vir'\to \Vir\to 0.$$

\sssec{} \label{unenv}
 A fact of crucial importance for us is that the forgetful functor 
from  chiral algebras to Lie-* algebras admits 
a left adjoint, i.e. that for every Lie-* algebra
$\B$ there exists a (unital) chiral algebra $\U(\B)$ (called the universal 
enveloping
algebra of $\B$) such that 
$$\on{Hom}_{\on{Lie}-*}(\B,\A)\simeq\on{Hom}_{\on{ch}}(\U(\B),\A),$$
functorially in $\A$. 

\smallskip

Let us sketch the construction of $\U(\B)$ (we shall assume for simplicity that 
$\B$ is torsion-free as
an $\O$-module).

\smallskip

The construction will be local and we can assume that $X$ is affine. Consider the 
D-module direct images
$B_1:=p_1{}_{*}(\Omega\boxtimes\B)$ and 
$B_2:=p_1{}_{*}(j_*j^*(\Omega\boxtimes\B))$ on $X$ (here $p_1$ is the projection 
$(x,y)\in X\times X\longrightarrow x\in X$). There is an obvious map $B_1\to B_2$ 
and from the fact that $\B$ is a Lie-* 
algebra we infer that both $B_1$ and $B_2$ are Lie algebras in the category of 
D-modules.\footnote{The category
of right D-modules is a tensor category under the operation $\overset{!}\otimes$, 
cf. \secref{Dmod}.}

\smallskip

We define $\U(\B)$ to be the ``vacuum representation'', i.e. 
$$\U(\B):=U(B_2){\underset{U(B_1)}\otimes}\CC,$$
where ``$U$'' is the universal enveloping algebra in the ordinary sense. 
By definition, for $x\in X$, the fiber of $\U(\B)$ at $x$ 
is identified with 
$U(DR^0(X\setminus x,\B)){\underset{U(DR^0(X,\B))}\otimes}\CC$.

\smallskip

Since $X$ is affine, the cokernel $$\on{coker}(DR^0(X,\B)\to DR^0(X\setminus 
x,\B))$$ is identified with
$$\on{coker}(DR^0(\Spec(\hat\O_x),\B)\to DR^0(\Spec(\hat\K_x),\B))\simeq 
B_x:=i^{!}_x(\B)[1].$$
This description implies that, first of all, $\U(\B)$ constructed in the above way 
will not change if we shrink $X$ to
a smaller open affine subset, and secondly, that $\U(\B)$ is naturally a chiral 
module over the Lie-* algebra $\B$. 

\medskip

Note, that by the construction, $\U(\B)$ carries a filtration 
$\U(\B)\simeq{\underset{n}\cup}\U(\B)_n,\,\, n\geq 0$
induced from the filtration on $U(B_2)$. We have $\U(\B)_0\simeq\Omega$, 
$\U(\B)_1\simeq \U(\B)_0\oplus \B$.

\begin{lem} \label{filtration}
Consider the chiral action of $j_*j^*(\B\boxtimes\U(\B))\to\Delta_{!}(\U(\B))$. We 
have:

\smallskip

\noindent {\em (a)}
$\on{Im}(\B\boxtimes \U(\B)_n)=\Delta_{!}(\U(\B)_n)$.

\smallskip

\noindent {\em (b)}
$\on{Im}(j_*j^*(\B\boxtimes \U(\B)_n))=\Delta_{!}(\U(\B)_{n+1})$.

\smallskip

\noindent {\em (c)}
$\U(\B)_n/\U(\B)_{n-1}\simeq \on{Sym}^n(\B)$.

\end{lem}

\sssec{}
 We shall now endow $\U(\B)$ with a chiral algebra structure. If fact, we will 
show that whenever $\M$ is a chiral $\B$-module, there exists a canonical map
$$j_*j^*(\U(\B)\boxtimes\M)\to\Delta_{!}(\M).$$

\smallskip

Consider the D-modules\footnote{Here $j_{y\neq z}$ and $j_{x\neq z\neq y}$ denote 
the open embeddings of the complements to the divisors 
$\Delta_{y=z}$ and $\Delta_{x=z}\cup \Delta_{y=z}$ in $X\times X\times X$ 
respectively.} $$j_{y\neq z}{}_*j_{y\neq z}^*(\Omega\boxtimes\Omega\boxtimes\B) 
\text{ and }
j_{x\neq z\neq y}{}_*j_{x\neq z\neq y}^*(\Omega\boxtimes\Omega\boxtimes\B)$$ on 
$X\times X\times X$
and let $B'_1$ and $B'_2$ denote 
their direct (D-module) images with respect to the projection $(x,y,z)\to (x,y)$.

\smallskip

We have an embedding $B'_1\to B'_2$ and the Lie-* algebra structure on $\B$ makes 
$B'_1$ and $B'_2$ 
into Lie algebras in the category of D-modules on $X\times X$. For any $(x,y)\in 
X\times X$, the fiber
of $B'_1$ (resp., of $B'_2$) at this point is identified with $DR^0(X\setminus 
y,\B)$ (resp., with 
$DR^0(X\setminus \{x,y\},\B)$). In particular, for a chiral $\B$-module $\M$ and 
for $(x,y)$ as above, the action of
$DR^0(X\setminus y,\B)$ on $M_y:=i_y^{!}(\M)[1]$ defines on $\Omega\boxtimes\M$ a 
structure of module over $B'_1$.

\smallskip

Consider now the induced module 
$\M':=U(B'_2){\underset{U(B'_1)}\otimes}(\Omega\boxtimes\M)$. It is clear that
$j^*(\M')\simeq j^*(\U(\B)\boxtimes \M)$ and that $\Delta^{!}(\M')[1]\simeq \M$. 
The canonical
map $$j_*j^*(\M')\to \Delta_{!}\Delta^{!}(\M')[1]$$ yields, therefore, a map
$j_*j^*(\U(\B)\boxtimes\M)\to\Delta_{!}(\M)$.

\smallskip

It is easy to show now that if we set $\M$ to be $\U(\B)$, then the resulting 
operation on $\U(\B)$ will satisfy
the axioms of a chiral bracket. Moreover, the above construction shows that any 
$\M$ which is chiral
module over $\B$ acquires a structure of chiral module over $\U(\B)$.

\ssec{Conformal blocks, correlation functions} \label{confblocks}
 From now on we shall assume that $\A$ is a chiral algebra possessing a unit. 
\sssec{}
 We have a surjection of D-modules $j_*j^*(\A\boxtimes\A)\to\Delta_{!}(\A)$; let 
$\A^{(2)}$ 
denote its kernel. 

\medskip

For a complete curve $X$ the space of conformal blocks $\conf(X,\A)$ of $\A$ is by 
definition the vector space 
$DR^2(X\times X,\A^{(2)})$. From the long exact cohomology sequence we obtain that
$$\conf(X,\A)\simeq\on{coker}(DR^1(X\times X\setminus 
\Delta(X),j^*(\A\boxtimes\A))\to DR^1(X,\A)).$$

\medskip

Let $x$ be a point of $X$ and recall that $A_x$ is a 
$DR^0(X\setminus x,\A)$-module.

\begin{prop}  \label{confbasic}
The space $\conf(X,\A)$ is identified naturally with the space of coinvariants
$$\on{coker}(DR^0(X\setminus x,\A)\otimes A_x\to A_x).$$
\end{prop}

\begin{proof}

We have a commutative square:
$$
\CD
DR^0(X\setminus x,\A)\otimes A_x  @>>> A_x  \\
@VVV   @VVV \\
DR^1(X\times X,j_*j^*(\A\boxtimes\A)) @>>> DR^1(X,\A) 
\endCD
$$
such that the vertical arrows are surjective. To prove the statement it remains to 
show
that the kernel of the map $A_x\to DR^1(X,\A)$ contains the image 
of\break 
$DR^0(X\setminus x,\A)\otimes A_x$. 

\smallskip

The unit in $\A$ provides a canonical element $\unit_x\in \A_x$
and the map\break 
$DR^0(X\setminus x,\A)\simeq DR^0(X\setminus x,\A)\otimes \unit_x\to A_x$ 
coincides
with the canonical map\break
 $DR^0(X\setminus x,\A)\to A_x$.  

Since the image of the latter is precisely the kernel of the map $A_x\to 
DR^1(X,\A)$ our
assertion follows.

\end{proof}

The following observation is useful for actual computations of the space of 
conformal blocks:

\begin{prop}  \label{actualcompuatations}
Let $\A$ be a chiral universal enveloping algebra of a Lie-* algebra $\B$. Then 
the natural map
$$A_x/(DR^0(X\setminus x,\B)\otimes A_x)\to A_x/(DR^0(X\setminus x,\A)\otimes 
A_x)\simeq\conf(X,\A)$$
is an isomorphism.
\end{prop}

\begin{proof}

We will prove a slightly more general assertion: we can replace $A_x$ in the 
formulation of the proposition
by the underlying space $M$ of any chiral $\A$-module $\M$ supported at $x$.

\smallskip

Let $\A_n$ be the filtration on $\A$ as in \secref{unenv}. We must show that the 
surjection
$$M/(DR^0(X\setminus x,\B)\otimes M)\twoheadrightarrow M/(DR^0(X\setminus 
x,\A_n)\otimes M)$$
is an isomorphism for any $n\geq 1$. We shall argue by induction.

\smallskip

For $a\in\Gamma(X\setminus x,\A_{n+1})$ and $m\in M$ consider the element 
$\rho(h(a),m)\in M$. By the construction,
it belongs to $\on{Im}(DR^0(X\setminus x,\A_{n+1})\otimes M)$ and it would suffice 
to show that
it also belongs to $\on{Im}(DR^0(X\setminus x,\A_n)\otimes M)$. 

\smallskip

Using the \lemref{filtration} and the fact that $X\setminus x$ is affine, for 
$a\in \Gamma(X\setminus x,\A_{n+1})$ as above,
we can find a section 
$b\boxtimes a'\cdot f(x,y)\in\Gamma(X\setminus x\times X\setminus 
x,j_*j^*(\B\boxtimes\A_n))$ 
such that $$(h\boxtimes\on{id})(\{b\boxtimes a'\cdot f(x,y)\})=a.$$

\smallskip

Using the Jacobi identity, the element $\rho(h(a),m)\in M$ is the image under the 
projection 
$$h\boxtimes h\boxtimes h:\Delta_{1,2,3}{}_{!}(\M)\to M$$ of the element
$$\rho(a',\rho(b,m\cdot f(x,y)))-
\rho(b,\rho(a',m\cdot f(y,x))).$$

Now, the first term obviously belongs to $\on{Im}(DR^0(X\setminus x,\A_n)\otimes 
M)$ and the second term
even belongs to 
$$\on{Im}(DR^0(X\setminus x,\B)\otimes M)\subset\on{Im}(DR^0(X\setminus 
x,\A_n)\otimes M).$$

\end{proof}

\sssec{}  
 Now let $x_1,\dotsc, x_k$ be a non-empty collection of distinct points of $X$. 
Each $A_{x_i}$ is a module over the 
Lie algebra $DR^0(X\setminus \{x_1,\dotsc, x_k\},\A)$; 
hence, so is the tensor product $A_{x_1}\otimes\ldots\otimes A_{x_k}$.

\begin{prop}  \label{confmany}
The space $\conf(X,\A)$ is canonically isomorphic to the space of coinvariants
$$\on{coker}(DR^0(X\setminus \{x_1,\dotsc,x_k\},\A)\otimes 
(A_{x_1}\otimes\ldots\otimes A_{x_k})\to
A_{x_1}\otimes\ldots\otimes A_{x_k}).$$
\end{prop}

\begin{proof}

When $k=1$ our assertion coincides with that of \propref{confbasic} above.
To treat the case $k\geq 2$ we shall proceed by induction, so we can assume that 
the assertion is true for $k-1$.

\medskip

We shall construct two mutually inverse maps $\phi$ and $\psi$ between 
$$(A_{x_1}\otimes\ldots\otimes A_{x_{k-1}})_{DR^0(X\setminus 
\{x_1,\dotsc,x_{k-1}\},\A)} \text{ and }
(A_{x_1}\otimes\ldots\otimes A_{x_k})_{DR^0(X\setminus \{x_1,\dotsc,x_k\},\A)}:$$

\noindent For $a_1\otimes\ldots\otimes a_{k-1}$ in $A_{x_1}\otimes\ldots\otimes 
A_{x_{k-1}}$ we set 
$\phi(a_1\otimes\ldots\otimes a_{k-1})\in A_{x_1}\otimes\ldots\otimes A_{x_k}$ to 
be 
$a_1\otimes\ldots\otimes a_{k-1}\otimes\unit_{x_{k}}$. It is clear that $\phi$ 
induces a well-defined map between
the spaces of coinvariants.

\medskip

\noindent Now let $a_1\otimes\ldots\otimes a_{k-1}\otimes b'$ 
be an element of $A_{x_1}\otimes\ldots\otimes A_{x_{k-1}}\otimes 
A_{x_k}$. Since $X\setminus\{x_1,\dotsc,x_{k-1}\}$ is affine, we can find an 
element 
$b\in DR^0(X\setminus\{x_1,\dotsc,x_k\},\A)$ that maps to $b'$ under the canonical 
map
$DR^0(X\setminus\{x_1,\dotsc,x_k\},\A)\to A_{x_k}$. 

\smallskip

We set $\psi(a_1\otimes\ldots\otimes a_{k-1}\otimes b')\in 
A_{x_1}\otimes\ldots\otimes A_{x_{k-1}}$ to be
$-b(a_1\otimes\ldots\otimes a_{k-1})$. It is again easy to see that $\psi$ is 
well-defined as a map between the
coinvariants and that $\phi$ and $\psi$ are mutually inverse.

\smallskip

It remains to show that the isomorphism 
$$\on{coker}(DR^0(X\setminus \{x_1,\dotsc,x_k\},\A)\otimes 
(A_{x_1}\otimes\ldots\otimes A_{x_k})\to
A_{x_1}\otimes\ldots\otimes A_{x_k})\simeq \conf(X,\A)$$
does not depend on the ordering of the points $x_1,\dotsc,x_k$. 
Obviously, it is enough to do this when $k=2$.

\smallskip

However, for two distinct points $x_1,x_2\in X$, our map $A_{x_1}\otimes 
A_{x_2}\to \conf(X,\A)$ coincides with the canonical map
$$A_{x_1}\otimes A_{x_2}\simeq i_{x_1\times x_2}^*(\A\boxtimes\A)[2]\simeq 
i_{x_1\times x_2}^*(\A^{(2)})[2]\to DR^2(X\times X,\A^{(2)})\simeq \conf(X,\A)$$
and our assertion follows, since the fact that $\{,\}$ is antisymmetric implies 
that 
the isomorphism $DR^2(X\times X,\A^{(2)})\simeq \conf(X,\A)$ is invariant with 
respect to the transposition.

\end{proof}

\sssec{} \label{confsummary}
 By globalizing the argument of \propref{confmany} we obtain that there exists a 
canonical map
of D-modules on $X^n$: 
$$\langle\ldots\rangle:j_*j^*(\A\boxtimes\ldots\boxtimes\A)\to
j_*j^*(\Omega\boxtimes\ldots\boxtimes\Omega)\otimes\conf(X,\A)$$ with the 
following properties:

\begin{itemize}

\item

For $x_1,\dotsc,x_n$ distinct points on $X$ the induced map 
$\langle\ldots\rangle:A_{x_1}\otimes\ldots\otimes 
A_{x_n}\to \conf(X,\A)$
coincides with the composition 
$$A_{x_1}\otimes\ldots\otimes A_{x_1}\to (A_{x_1}\otimes\ldots\otimes 
A_{x_n})_{DR^0(X\setminus \{x_1,\dotsc,x_n\},\A)}\simeq
\conf(X,\A).$$
In particular, the map $\langle\ldots\rangle$ is equivariant with respect to the 
action of the symmetric group $S^n$
on $j_*j^*(\A^{\boxtimes n})$ and on $j_*j^*(\Omega^{\boxtimes n})$.

\item

The square 
$$
\CD
j_*j^*(\A^{\boxtimes n-1}\boxtimes\Omega) 
@>{\langle\ldots\rangle\boxtimes\on{id}}>>
j_*j^*(\Omega^{\boxtimes n})\otimes\conf(X,\A) \\
@V{\on{id}^{\boxtimes n-1}\boxtimes\unit}VV    @V{\on{id}^{\boxtimes n}}VV \\
j_*j^*(\A^{\boxtimes n}) @>{\langle\ldots\rangle}>> 
j_*j^*(\Omega^{\boxtimes n})\otimes\conf(X,\A)
\endCD
$$
is commutative.

\item

The square
$$
\CD
j_*j^*(\A^{\boxtimes n}) @>{\langle\ldots\rangle}>> 
j_*j^*(\Omega^{\boxtimes n})\otimes\conf(X,\A) \\
@V{\{,\}\boxtimes\on{id}^{\boxtimes n-1}}VV   
@V{\can_\Omega\boxtimes\on{id}^{\boxtimes n-1}}VV  \\
\Delta_{x_1=x_2}{}_{!}(j_*j^*(\A^{\boxtimes n-1})) @>{\langle\ldots\rangle}>>
\Delta_{x_1=x_2}{}_{!}(j_*j^*(\Omega^{\boxtimes n-1}))\otimes\conf(X,\A).
\endCD
$$
commutes as well.

\end{itemize}

(The first two of the above properties are evident from the definitions. The third 
one can be easily deduced 
from the second one using the fact that the upper horizontal arrow in the first 
commutative diagram is surjective).

\medskip

By passing to the corresponding left D-modules, we obtain a canonical 
$S^n$-invariant map
$$\langle\ldots\rangle:j_*j^*({}^l\A\boxtimes\ldots\boxtimes{}^l\A)\to
j_*j^*(\O\boxtimes\ldots\boxtimes\O)\otimes\conf(X,\A).$$

\smallskip

\begin{rem}

For a functional $\chi$ on $\conf(X,\A)$ and a section 
$a_1\boxtimes\ldots\boxtimes a_n\in {}^l\A^{\boxtimes n}$ we can 
produce, therefore, a function
$\langle a_1,\dotsc,a_n\rangle_\chi$ on $X^n$ with poles on the diagonal 
divisor. 
This function should be thought of
as an analog of the $n$-point correlation function of the fields 
$a_1,\dotsc,a_n$.

\end{rem}

\newpage

\centerline{\bf Chapter II. CFT Data (Algebraic Version)}

\medskip

Recall that $\Xf$ denotes the universal curve over the moduli stack $\Mf$ of 
smooth complete curves.
In this section the default meaning of $\O$ (resp., $\Omega$,
$\Theta$, $D$) is the structure sheaf (resp., the sheaf of relative $1$-forms, 
the sheaf of vertical
vector fields, vertical differential operators) on $\Xf$; D-modules (unless 
specified otherwise) are right modules over the sheaf of rings $D$.

\ssec{Local $\O$-modules on $\Xf$}
 On the intuitive level, we want to call a sheaf $\L$ on $\Xf$ local if the 
following holds:
whenever $X$ and $X'$ are two curves such that their (analytic) open subsets 
$U\subset X$
and $U'\subset X'$ are identified, the sheaves $(\L|X)|U$ and $(\L|X')|U'$ are 
identified as well.
For instance, if we are dealing with a $2$-dimensional CFT which is defined for 
every Riemann surface, its
fields form naturally a local sheaf over $\Xf$. 

\sssec{}  \label{local}
 Let $\L$ be an $\O$-module on $\Xf$, i.e. for a scheme $S$ with a map to $\Mf$ we 
have
an $\O$-module $\L^S$ on the corresponding family of curves $X^S$. We say that 
$\L$ is local if it possesses
the following additional structure:

\medskip

Let $I$ be a local Artinian scheme and let $X^{S\times I}$ be a family of curves 
over $X\times I$ extending
$X^S$. Let, in addition, $x_1^{S\times I},\dotsc,x_n^{S\times I}$ be sections 
$S\times I\to X^{S\times I}$
(we shall denote by $x_1^S,\dotsc,x_n^S$ their restrictions to $S\hookrightarrow 
S\times I$). Assume now that
$\phi^{S,I}$ is an isomorphism:
$$X^{S\times I}\setminus\{x_1^{S\times I},\dotsc,x_n^{S\times I}\}\simeq 
(X^S\setminus\{x_1^S,\dotsc,x_n^S\})\times I.$$
A local structure on $\L$ is lifting of $\phi^{S,I}$ to an isomorphism of sheaves
$$\beta^{S,I}:\phi^{S,I}{}^*
(\L^S|X^S\setminus\{x_1^S,\dotsc,x_n^S\}\otimes\O(I))\to 
\L^{S\times I}|X^{S\times I}\setminus\{x_1^{S\times I},\dotsc,
x_n^{S\times I}\}.$$
These $\beta^{S,I}$'s for various 
$X^{S\times I},x_1^{S\times I},\dotsc,
x_n^{S\times I},\phi^I$ must satisfy
two conditions. The first condition is a compatibility in the obvious sense in the 
situation when a map 
$S\times I\to\Mf$ factors as $S\times I\to S'\times I'\to\Mf$. To formulate the 
second condition we need
to introduce some notation.

\bigskip

Let $X^S$ be again a curve over the 
base $S$ and let $y^S_1,\dotsc, y^S_n$ 
and $z^S_1,\dotsc, z^S_k$ be two disjoint
collections of sections $S\to X^S$. Let $I$ and $J$ be two local Artinian schemes 
and let
$(X^{S\times I},y^{S\times I}_1,\dotsc, y^{S\times I}_n,\phi^I)$
and $(X^{S\times J},z^{S\times J}_1,\dotsc, z^{S\times J}_k,\phi^J)$ be two pieces 
of data as above.
(In what follows, to simplify the notation, we shall abbreviate 
$\{y_1,\dotsc,y_n\}$ and $\{z_1,\dotsc, z_m\}$
to $y$ and $z$ respectively.)

\medskip

We have sections $\tily^J:S\times J\to X^{S\times J}$ and $\tilz^I:S\times I\to 
X^{S\times I}$
and we can construct a curve $X^{S\times I\times J}$ over $S\times I\times J$ 
equipped with sections
$y^{S\times I\times J},z^{S\times I\times J}$
such that:

\begin{align}
&X^{S\times I\times J}\setminus 
y^{S\times I\times J}\,\,{\overset{\phi^{S\times J,I}}\simeq}
(X^J\setminus \tily^J)\times I \\
&X^{S\times I\times J}\setminus 
z^{S\times I\times J}\,\,{\overset{\phi^{S\times I,J}}\simeq}
(X^I\setminus \tilz^I)\times J \\
&X^{S\times I\times J}
\setminus\{y^{S\times I\times J},z^{S\times I\times J}\}\,
{\overset{\phi^{S,I\times J}}\simeq} 
(X^S\setminus\{y^S,z^S\})\times I\times J .
\end{align}

Note that the composition:
$$(\phi^{S,J}\times\on{id}_I)\circ \phi^{S\times J,I}:
X^{S\times I\times J}\setminus\{y^{S\times I\times J},z^{S\times I\times J}\}\to 
(X^S\setminus\{y^S,z^S\})\times I\times J$$
coincides with $\phi^{S,I\times J}$.

\medskip

The second condition on $\beta$ reads as follows:

\medskip

\noindent We need that the two isomorphisms between 
\begin{align*}
&\phi^{S,I\times J}{}^*((\L^S)|X
\setminus\{y^S,z^S\}\otimes \O(I)\otimes\O(J))\\
\intertext{and}
&(\L^{S\times I\times J})|X^{S\times I\times J}
\setminus\{y^{S\times I\times J},z^{S\times I\times J}\},
\end{align*}
namely,
$\beta^{S,I\times J}$ and 
$(\beta^{S,J}\times\on{id}_I)\circ((\phi^{S,J}\times\on{id}_I)^*(\beta^{S\times 
J,I}))$,
coincide.

\bigskip

The most basic example of a local $\O$-module on $\Xf$ is the sheaf $\O$ itself.
Local $\O$-modules on $\Xf$ form an abelian category (morphisms in this category 
are, by definition, maps of 
$\O$-modules compatible with the $\beta^{S,I}$'s).
It is easy to see that the sheaves $\Omega^i$, $D$, etc. are all local 
$\O$-modules in a natural way.

\bigskip

In a similar way one defines local $\O$-modules on $\Xf^n$.

\smallskip

An $\O$-module $\L$ on $\Xf^n$ assigns, by definition, to every
scheme $S$ mapping to $\Mf$ an $\O$-module $\L^S$ over 
$$X^S{\underset{S}\times}\ldots{\underset{S}\times}X^S,$$
where $X^S$ is as before. A local structure on $\L$ attaches to each data 
$X^{S,I},x_k^{S,I},\phi^{S,I}$ as above
an isomorphism
\begin{align*}
&\beta^{S,I}:\phi^{S,I}{}^*(\L^S|(X^S\setminus \{x^S_1,\dotsc,x^S_n\})
{\underset{S}\times}\ldots{\underset{S}\times}(X^S\setminus 
\{x^S_1,\dotsc,x^S_n\})\otimes\O(I))\simeq \\
&\L^{S\times I}|(X^{S\times I}\setminus \{x^{S\times I}_1,\dotsc,x^{S\times 
I}_n\})
{\underset{S\times I}\times}\ldots{\underset{S\times I}\times}(X^{S\times 
I}\setminus \{x^{S\times I}_1,\dotsc,x^{S\times I}_n\}).
\end{align*}

These isomorphisms should satisfy two conditions analogous to what we had in 
the case $n=1$.

\medskip

For example, if $\L_1$ and $\L_2$ are local $\O$-modules on $\Xf$, then 
$\L_1\boxtimes \L_2$ is a local $\O$-module
on $\Xf^2$. In addition, if $\L$ is a local $\O$-module on $\Xf^2$, so is 
$j_*j^*(\L)$, and $\Delta^*(\L)$
is a local $\O$-module back on $\Xf$.

\smallskip

In particular, we have the following assertion:

\begin{lem} \label{locchiral}
Let $\B$ be a local $\O$-module on $\Xf$ with the following additional structure: 
$\B$ is a D-module
and a Lie-* algebra such that:

\smallskip

\noindent {\em (a)} The map $\B\otimes D\to\B$ is a map of local $\O$-modules on 
$\Xf$.

\smallskip

\noindent {\em (b)} The Lie-* bracket $\B\boxtimes\B\to\Delta_{!}(\B)$ is a map of 
local $\O$-modules
on $\Xf^2$ (the RHS is a local $\O$-module on $\Xf^2$ in a natural way).

\smallskip

\noindent Then the chiral universal enveloping algebra $\U(\B)$ is a local 
$\O$-module on $\Xf$ 
and conditions {\em (a)} and {\em (b)} hold for $\U(\B)$ as well.
\end{lem}

\sssec{} \label{anotherlocal}
 We shall now describe a typical way of producing local $\O$-modules.

\medskip

Consider the standard $1$-dimensional local ring $\hat\O:=\CC[[t]]$ and the 
corresponding
local field $\hat\K:=\CC((t))$. 
Let $\Aut^+_0$ (resp., $\Aut^+$, $\Aut$) be the group scheme 
of automorphisms of $\hat\O$ that preserve
the maximal ideal (resp., the group ind-scheme of all automorphisms of 
$\hat\O$, the group ind-scheme of all automorphisms of $\hat\K$).
We shall
denote by $\Vir^+_0\subset\Vir^+\subset\Vir$ the corresponding Lie algebras.

\smallskip

Consider the scheme $\hat\Xf$ that classifies the following data: $$(X\in\Mf,\, 
x\in X,\alpha:
\hat\O_x\to\hat\O),$$ where $\alpha$ is an isomorphism that preserves the maximal 
ideal. The group
$\Aut^+_0$ acts on this scheme by changing $\alpha$ and the 
quotient
is identified with $\Xf$. Let $q$ denote the projection of $\hat\Xf$ to $\Xf$.

\smallskip

If now $V$ is a representation of the group $\Aut^+_0$, the $\O_{\Xf}$-module 
$V\otimes\O_{\hat\Xf}$ is 
$\Aut^+_0$-equivariant and hence descends to an $\O$-module $\on{Ass}(V)$ over 
$\Xf$. The next assertion is evident
from the definitions:

\begin{lem}
For $V$ as above, the $\O$-module $\on{Ass}(V)$ has a natural structure of local 
$\O$-module.
\end{lem}

\smallskip

For instance, for $i\in\ZZ$ let $L_i$ be the $1$-dimensional $\Aut^+_0$-module of 
\secref{linalg}. ($\Aut^+_0$ acts
on $L_i$ via the $i$-th power of the standard character.) We have:
$$\on{Ass}(L_i)\simeq \Omega^i.$$

\sssec{} \label{Lieder}
 Let $\L$ be a local $\O$-module on $\Xf$. We shall now construct a 
map\footnote{For two $\O$-modules $\L_1$ and $\L_2$ on $\Xf$, it will be sometimes 
more convenient to use the
notation $\L_1\boxtimes \L_2(\infty\cdot\Delta)$ for $j_*j^*(\L_1\boxtimes 
\L_2)$.}
$$\nabla_{loc}:\L\to p_1{}_*(\L\boxtimes\Omega^2(\infty\cdot\Delta)).$$
Let $x\in X\in\Mf$ be a point of $\Xf$, let $l$ be a local section of $\L$ on
$\Xf$ and let $\xi$ be an element of $\hat\K_x/\hat\O_x\otimes\Theta$.

\medskip

Constructing $\nabla_{loc}$ amounts to our being able to evaluate 
$(\on{id}\boxtimes\on{Res}_x)(\nabla_{loc}(l)\otimes(1\boxtimes\xi))$ as a local 
section
of $\L|X$ defined outside of $x$. (The triple $(X,x,\xi)$ should be thought of as 
living over 
a base $S$, where $S$ is a scheme mapping to $\Mf$; we shall denote by $\L|X$ the 
$\O$-module $\L^S$ over $X^S$.)

\medskip

However, $\xi$ as above is the same as a deformation of $X$ over the scheme 
$I=\on{Spec}(\CC[t]/t^2)$
which is trivialized outside of $x$. As $\L$ is local, we have a section
$\beta^I{}^{-1}(l)$ of $((\L|X)|X\setminus x)\otimes \CC[t]/t^2$ and we set
$$(\on{id}\times\on{Res}_x)(\nabla_{loc}(l)\otimes(1\boxtimes\xi)):=\d_t(\beta^I{}
^{-1}(l)).$$

\medskip

Since the isomorphism $\beta^I$ above respects the $\O$-module structure, the map 
$\nabla_{loc}$ satisfies an
analog of the Leibnitz rule: for a function $f$ on $\Xf$ and a local section $l$ 
of $\L$ we have:
$$\nabla_{loc}(f\cdot l)=(f\boxtimes 1)\cdot\nabla_{loc}(l)+ 
\nabla_{loc}(f)\otimes (l\boxtimes 1)\in 
p_1{}_*(\L\boxtimes\Omega^2(\infty\cdot\Delta)),$$
(here $\nabla_{loc}(f)$ is the corresponding operation for $\L=\O$).

\bigskip

Analogously, when $\L$ is a local $\O$-module on $\Xf^n$ we have a map
$$\nabla_{loc}:\L\to 
p_{1,\dotsc,n}{}_*(\L\boxtimes\Omega^2(\infty\cdot\Delta_{x_{n+1}\neq x_k})),$$
that satisfies the Leibnitz rule in the same sense as above.

\medskip

In particular, for a local $\O$-module $\L$ on $\Xf$, the $\O$-module 
$\L\boxtimes\Omega^2(\infty\cdot\Delta)$
on $\Xf^2$ is local as well and we can form a map
$$\nabla_{loc}^2:\L{\overset{\nabla_{loc}}\longrightarrow}p_1{}_*(\L\boxtimes\Omega^2(\infty\cdot\Delta))
{\overset{\nabla_{loc}}\longrightarrow}
p_1{}_*(\L\boxtimes\Omega^2\boxtimes\Omega^2(\infty\cdot\Delta)).$$

\begin{lem} 
Let $\sigma_{2,3}$ be the transposition of the last 2 coordinates acting on the 
sheaf
$p_1{}_*(\L\boxtimes\Omega^2\boxtimes\Omega^2(\infty\cdot\Delta))$. We have:
$$\sigma_{2,3}\circ\nabla_{loc}^2=\nabla_{loc}^2,$$
i.e. the image of $\nabla^2_{loc}$ is symmetric in the last two variables.
\end{lem}

The proof follows immediately from the second condition on the
$\beta$'s in the definition of local $\O$-modules.

\medskip

Let us consider a few examples:

\smallskip

Note that the cotangent bundle of $\hat\Xf$ is identified in an 
$\Aut^+_0$-equivariant way with 
$$q^*p_1{}_*(\O\boxtimes\Omega^2(\infty\cdot\Delta)).$$
Therefore, the De Rham differential on $\hat\Xf$ yields maps
$$d^0:\O\to \Omega_{\Xf}\simeq p_1{}_*(\O\boxtimes 
\Omega^2(\Delta))\hookrightarrow
p_1{}_*(\O\boxtimes \Omega^2(\infty\cdot\Delta)) \text{ and }$$
$$d^1:p_1{}_*(\O\boxtimes\Omega^2(\infty\cdot\Delta))\to p_1{}_*(\O\boxtimes
\Omega^2\boxtimes\Omega^2(\infty\cdot\Delta)).$$
We leave it to the reader to check that $d^0$ is the map $\nabla_{loc}$ for 
$\L=\O$ and that $d^1$ is the map $\nabla_{loc}$
for the local $\O$-module $\O\boxtimes\Omega^2(\infty\cdot\Delta)$ on $\Xf^2$.

\medskip

We claim, in fact, that for an $\O$-module $\L$ on $\Xf$, a structure of local 
$\O$-module on it is equivalent to the data
of $\nabla_{loc}$:

\begin{thm}  \label{descrlocal}
The following three categories are equivalent:

\begin{itemize}

\item
The category of local $\O$-modules on $\Xf$.

\item
The category of weakly $\Aut^+_0$-equivariant left D-modules on $\hat\Xf$.

\item
The category of $\O$-modules $\L$ on $\Xf$ endowed with a map 
$\nabla_{loc}:\L\to$\break
$p_1{}_*(\L\boxtimes\Omega^2(\infty\cdot\Delta))$ such that:

\smallskip

\noindent {\em (a)} For any function $f$ on $\Xf$ and a local section $l$ of $\L$
$$\nabla_{loc}(f\cdot l)=(f\boxtimes 1)\cdot\nabla_{loc}(l)+ 
\nabla_{loc}(f)\otimes (l\boxtimes 1)
\in p_1{}_*(\L\boxtimes\Omega^2(\infty\cdot\Delta)).$$

\smallskip

\noindent {\em (b)} If we define (assuming {\em (a)}) 
$$\nabla_{loc}:
\L\boxtimes\Omega^2(\infty\cdot\Delta)\to 
p_1{}_*(\L\boxtimes\Omega^2\boxtimes\Omega^2)(\infty\cdot\Delta)$$
according to the Leibnitz rule applied to $\nabla_{loc}$ acting on $\L$ and on 
$\Omega^2$,
then the composition 
$$\nabla_{loc}{}^2:=\nabla_{loc}\circ\nabla_{loc}:
\L\to p_1{}_*(\L\boxtimes\Omega^2\boxtimes\Omega^2)(\infty\cdot\Delta)$$ satisfies
$\sigma_{2,3}\circ\nabla_{loc}{}^2=\nabla_{loc}{}^2$.

\smallskip

Morphisms in this category are $\O$-module maps $\L_1\to\L_2$ compatible with the 
$\nabla_{loc}$'s.

\end{itemize}

\end{thm}

The proof of this theorem is quite straightforward and we shall omit it, 
especially since this statement will not be
used in the sequel.

\medskip

We shall conclude this subsection by showing that if $\L$ is a local $\O$-module 
on $\Xf$, there is a natural 
action of $\Theta$ on it ``by Lie derivations'' (cf. \secref{firstexamples}).

\smallskip

\noindent Consider the composition 
$$\nabla_{loc}':\L{\overset{\nabla_{loc}}\longrightarrow} p_1{}_*(\L\boxtimes 
\Omega^2(\infty\cdot\Delta))\to
p_1{}_*(\L\boxtimes\Omega^2(\infty\cdot\Delta)/\L\boxtimes\Omega^2).$$
This map is $\O_{\Mf}$-linear and it satisfies:
$$\nabla'_{loc}(f\cdot l)=(f\boxtimes 1)\cdot\nabla'_{loc}(l)+(l\boxtimes 
1)\otimes df,$$
where $df\in\Omega$ is viewed as an element of 
$\O\boxtimes\Omega^2(\infty\cdot\Delta)/\O\boxtimes\Omega^2$ via
$$\Omega\simeq \O\boxtimes\Omega^2(\Delta)/\O\boxtimes\Omega^2
\subset\O\boxtimes\Omega^2(\infty\cdot\Delta)/\O\boxtimes\Omega^2.$$

\medskip

If $l$ is a local section of $\L$ and $\xi$ is a vertical vector field on $\Xf$ 
(i.e. a local section of $\Theta$), we set
$$\on{Lie}_\xi(l):=(\on{id}\boxtimes\on{Res}_\Delta)(\nabla'_n(l)\otimes(1\boxtimes\xi))\in \L.$$
The properties of $\nabla_{loc}$ imply that the map $(l,\xi)\to \on{Lie}_\xi(l)$ 
is a bi-differential operator and is
$\O_{\Mf}$-linear in both variables. In addition, we have:
$$\on{Lie}_\xi(f\cdot l)=f\cdot \on{Lie}_\xi(l)+\on{Lie}_\xi(f)\cdot l.$$

\ssec{A formulation of CFT (central charge $0$)}
 \sssec{}  \label{defcft}
 A CFT data (of central charge $0$) for us will consist of 

\smallskip

\noindent (1) A local $\O$-module $\A$ on $\Xf$ endowed with a structure of a 
chiral algebra
(in particular, $\A$ is a right D-module).

\smallskip

\noindent (2) A map of local 
$\O$-modules $T:\Theta\to\A$ (called the energy-momentum tensor).

\medskip

These data must satisfy the following conditions:

\smallskip

\noindent (a) The chiral bracket on $\A$: 
$$j_*j^*(\A\boxtimes\A){\overset{\{,\}}\longrightarrow}\Delta_{!}(\A)$$
is a map of local $\O$-modules on $\Xf^2$. Moreover, the unit map 
$\unit:\Omega\to{}\A$ is also a map of local
$\O$-modules.

\smallskip

\noindent (b) Let $\xi$ be a vertical vector field on $\Xf$ and consider the map 
$\A\to\A$ given by
$$l\to (\on{id}\boxtimes h)(\{l,T(\xi)\})-\on{Lie}_\xi(l).$$
We need this map to coincide with the action of $\xi$ on $\A$ coming from the 
right D-module
structure on $\A$.

\bigskip

It follows from (a) and (b) that the map $\A\otimes D\to\A$ defining the right 
D-module structure on $\A$
is a map of local $\O$-modules. 

\medskip

Recall the Lie-* algebra $\T=\Theta\otimes D$ that we constructed in 
\secref{firstexamples}. The map
$T$ gives rise to a map $\T\to\A$ and conditions (a) and (b) above guarantee that 
this is a homomorphism
of Lie-* algebras. 

\begin{rem}
In the physics literature one often encounters the term ``primary field of weight 
$i$''. In our language, this is the same
as a map of local $\O$-modules $\Omega^{1-i}\to \A$.
\end{rem}

\sssec{}  \label{intralgebroid}
 Let $^l\A$ be the corresponding left D-module $^l\A:=\A\otimes\Theta$ and let 
$^lT$ denote
the corresponding map $\Omega^{-2}\to{}^l\A$. We will show now that $^l\A$ 
acquires a (left) 
module structure over some canonical Lie algebroid on $\Xf$ (cf. \secref{Dmod} for 
the notion
of a Lie algebroid). This will later on enable us 
to produce a connection along $\Mf$ on the space of conformal blocks of $\A$.

\smallskip

Recall the scheme $\hat\Xf$ that we introduced in the previous subsection. It is 
known that the natural
$\Aut^+_0$-action on it extends to an action of the whole of $\Aut$. Let us recall 
how this action is constructed.

\medskip

For a point $(x\in X\in\Mf,\alpha)$
of $\hat\Xf$ and a point $g\in\Aut$ with values in a local Artinian scheme 
$\Spec(I)$,
we must produce a $\Spec(I)$-valued point $(x_I\in X_I,\alpha_I)$ of $\hat\Xf$. 

\smallskip

We let $X_I$ be unchanged as a topological space; moreover, we set $X_I\setminus 
x_I$ to be
$(X\setminus x)\times \Spec(I)$. To fix the triple $(x_I\in X_I,\alpha_I)$ it 
remains to specify
which of the meromorphic functions on $(X\setminus x)\times \Spec(I)$ are regular 
at $x_I\in X_I$.

\smallskip

The set of meromorphic functions on $X\times \on{Spec}(I)$ 
embeds into $\hat\K\otimes I$ by means of $\alpha$. We declare a meromorphic 
function on
$(X\setminus x)\times \Spec(I)$ to be regular for $X_I$ at $x_I$ if its image in
$\hat\K\otimes I$ is transformed by $g:\hat\K\otimes I\to\hat\K\otimes I$ into an 
element of $\hat\O\otimes I$.

\medskip

In particular, the Lie algebra $\Vir$ acts on $\hat\Xf$ by vector fields and we
can form the Lie algebroid $\tilgf=\Vir\otimes\O_{\hat\Xf}$ (resp., 
$\tilgf^+=\Vir^+\otimes\O_{\hat\Xf}$,
$\tilgf^+_0=\Vir^+_0\otimes\O_{\hat\Xf}$) on $\hat\Xf$; 
all the three are Harish-Chandra Lie algebroids with respect to $\Aut^+_0$.
Let $\gf$, $\gf^+$ and $\gf^+_0$ denote the corresponding Lie algebroids on $\Xf$, 
i.e.
$\tilgf\simeq q^*(\gf)$ and similarly for $\gf^+$ and $\gf^+_0$;
the Harish-Chandra property implies that $\gf^+_0$ is an ideal in the other two.

\smallskip

We have:
$$\gf\simeq p_1{}_*(\O\hat\boxtimes\Theta(\infty\cdot\Delta)) \text{ and }
\gf/\gf^+_0\simeq 
p_1{}_*(\O\boxtimes\Theta(\infty\cdot\Delta)/\O\boxtimes\Theta(-\Delta)),$$
where the notation $\hat\boxtimes$ means completion along the diagonal divisor.

\smallskip

Moreover $\gf^+/\gf^+_0\hookrightarrow \gf/\gf^+_0$ is
identified with
$$\Theta\simeq \O\boxtimes\Theta/\O\boxtimes\Theta(-\Delta)\hookrightarrow 
p_1{}_*(\O\boxtimes\Theta(\infty\cdot\Delta)/\O\boxtimes\Theta(-\Delta)).$$

\bigskip

There are two maps $\gf\underset{\CC}\otimes {}^l\A\to{}^l\A$:

\medskip

\noindent The first map $\tau_{loc}:\gf\otimes {}^l\A\to{}^l\A$ is defined as 
follows:

\smallskip

\noindent A typical local section of $\gf\underset{\CC}\otimes {}^l\A$ can be 
written as
$(f\cdot 1\boxtimes\xi)\otimes a$, for $f\in$\break
$p_1{}_*(\O\hat\boxtimes\O(\infty\cdot\Delta))$,
$\xi\in\Theta$ and $a\in{}^l\A$. We set\footnote{Using \propref{descrlocal} one 
can rephrase the definition of $\tau_{loc}$ as follows: since
$q^*({}^l\A)$ is a weakly $\Aut^+_0$-equivariant D-module on $\hat\Xf$, it is also 
a weakly
equivariant module over the algebroid $\tilgf$; this defines on $^l\A$ a structure 
of $\gf$-module,
which is the same as a map $\tau_{loc}$.}
$$\tau_{loc}((f\cdot 1\boxtimes\xi),a)=
(\on{id}\boxtimes\on{Res}_\Delta)(f\cdot\nabla_{loc}(a)\otimes (1\boxtimes
\xi)).$$

\medskip

\noindent To define the second map $\tau_{em}$ (the subscript ``em'' here stands 
for ``energy-\break
momentum''), let
$(f\cdot 1\boxtimes\xi)\otimes a\in \gf\underset{\CC}\otimes {}^l\A$ be as above 
and we set\footnote{Here we have denoted by 
$(\on{id}\boxtimes h)(\{,\})$ the map
$^l\A\hat\boxtimes \A(\infty\cdot\Delta)\to{}^l\A$ that comes from the chiral 
bracket on $\A$.}
$$\tau_{em}((f\cdot 1\boxtimes\xi),a)=(\on{id}\boxtimes h)(\{f\cdot a\boxtimes 
T(\xi)\})\in{}^l\A.$$

\smallskip

Let $\tau$ be the difference: $$\tau=\tau_{loc}-\tau_{em}:\gf\otimes 
{}^l\A\to{}^l\A.$$

\begin{lem}
The map $\tau$ is a left action of the algebroid $\gf$ on $^l\A$. Moreover, $\tau$ 
annihilates the ideal
$\gf^+_0\subset\gf$ and the restriction of $\tau$ to $\Theta\simeq\gf^+/\gf^+_0$ 
coincides with 
the action of $\Theta$ on $^l\A$ coming from the structure of a left D-module on 
$^l\A$.
\end{lem}

The proof follows from conditions (a) and (b) in \secref{defcft} and we leave the 
verification to the reader.

\sssec{}
 Let $\A$ be as in \secref{defcft}. For each $X\in\Mf$ we can construct a vector 
space $\conf(X,\A)$ and globally
over $\Mf$ this construction yields an $\O_\Mf$-module $\conf(\Xf,\A)$.

\medskip

According to the discussion in \secref{confsummary}, we have a map of left 
D-modules on $\Xf$:
$$\langle\ldots\rangle:{}^l\A\twoheadrightarrow \pi^*(\conf(\Xf,\A)).$$

\begin{prop}   \label{nocentral}
There exists a (unique) left $D_\Mf$-module structure on $\conf(\Xf,\A)$ such that 
when $\pi^*(\conf(\Xf,\A))$ is viewed
as an $\gf/\gf^+_0$-module via $\gf/\gf^+_0\to\Theta_{\Xf}$, the map 
$\langle\ldots\rangle$ is compatible
with the $\gf/\gf^+_0$-module structure.
\end{prop}

\begin{proof}

Since the map $\langle\ldots\rangle$ is surjective, there is at most one 
$\gf/\gf^+_0$-module structure on
$\pi^*(\conf(\Xf,\A))$ compatible with that on $^l\A$.

\medskip

To see that it is correctly defined, we must show that if $c$ is a local section 
of $^l\A$ of the form
$$c=(h\boxtimes\on{id})(\{f_0\cdot a\boxtimes b\}),$$
where $f_0\cdot a\boxtimes b\in p_2{}_*(\A\boxtimes{}^l\A(\infty\cdot\Delta))$ and 
if 
$f\cdot 1\boxtimes \xi\in p_1{}_*(\O\hat\boxtimes\Theta(\infty\cdot\Delta))$ is a 
local section of $\gf/\gf^+_0$,
the section $c':=\tau(f\cdot 1\boxtimes \xi,c)$ is a sum of sections of the form
$$(h\boxtimes\on{id})(\{f_0'\cdot a'\boxtimes b'\}),$$
for $f'_0\cdot a'\boxtimes b'$ a local section of 
$p_2{}_*(\A\boxtimes{}^l\A(\infty\cdot\Delta))$.

\smallskip

This, however, follows from the next lemma.

\begin{lem}
Let $f_0\cdot a\boxtimes b$ and $f\cdot 1\boxtimes \xi$ and $c'$ be as above and 
let
$a''$ be the section of $p_2{}_*(\A\boxtimes\O(\infty\cdot\Delta))$ defined as
$$a''=(\on{id}\boxtimes\on{id}\boxtimes\on{Res}_{x_2=x_3})
(\nabla_{loc}(f_0\cdot a\boxtimes 1)\otimes(1\boxtimes f\boxtimes \xi)).$$
Then
$$c'=(h\boxtimes\on{id})(\{f_0\cdot a\boxtimes \tau(f\cdot 1\boxtimes \xi,b)\})+
(h\boxtimes\on{id})(\{a''\otimes(1\boxtimes b\})).$$
\end{lem}

\bigskip

Next, we must show that the action of $\gf/\gf^+_0$ on $\pi^*(\conf(\Xf,\A))$ 
factors through $\Theta_{\Xf}$.
The map $\gf/\gf^+_0\to\Theta_{\Xf}$ is surjective and its 
kernel is identified with
$$p_1{}_*(\O\boxtimes\Theta(\infty\cdot\Delta))\subset 
p_1{}_*(\O\hat\boxtimes\Theta(\infty\cdot\Delta)),$$
(in other words, the fiber of the above kernel at $x\in X\in\Mf$ 
is identified with the subspace
$H^0(X\setminus x,\Theta_X)\subset \Theta\underset{\O}\otimes\hat\K_x$).

\medskip

For $f\cdot 1\boxtimes\xi\in p_1{}_*(\O\boxtimes\Theta(\infty\cdot\Delta))$ and 
for 
a local section $a$ of $^l\A$ we have by definition:

$$\langle \tau(f\cdot 1\boxtimes\xi,a)\rangle=
\langle (\on{id}\boxtimes\on{Res}_\Delta)(f\cdot\nabla_{loc}(a)\otimes 
(1\boxtimes\xi)) \rangle-
\langle(\on{id}\boxtimes h)(\{f\cdot a\boxtimes T(\xi)\})\rangle.$$

However, the first term vanishes by the residue formula, since
$f\cdot\nabla_{loc}(a)\otimes (1\boxtimes\xi)$ is a section of 
$p_1{}_*({}^l\A\boxtimes\Omega(\infty\cdot\Delta))$
and the second term vanishes by the very definition of $\conf(\Xf,\A)$.

\bigskip

We have shown, therefore, that $\pi^*(\conf(\Xf,\A))$ is a D-module on $\Xf$ such 
that the action of vertical
vector fields on it is the ``tautological'' one. Therefore, the $\O$-module 
$\conf(\Xf,\A)$ on $\Mf$ is
naturally a $D_\Mf$-module.

\end{proof}

Here is another characterization of the connection on $\conf(\Xf,\A)$:

\smallskip

Let $\nabla$ denote the covariant derivative acting on $\pi^*_k(\conf(\Xf,A))$; in 
particular, 
there is a map\footnote{In the formula below, 
$\langle f\cdot a_1\boxtimes\ldots
\boxtimes a_k\boxtimes {}^lT\rangle^{k+1}$ will
denote a section of $j_*j^*(\O^{\boxtimes k}\boxtimes\Omega^2)$ characterized by 
the property that
for $\eta\in\Omega^{-2}$:
$$\langle f\cdot a_1\boxtimes\ldots\boxtimes 
a_k\boxtimes {}^lT\rangle\otimes 
(1\boxtimes\ldots\boxtimes 1\boxtimes\eta):=
\langle f\cdot a_1\boxtimes\ldots\boxtimes a_k\boxtimes {}^lT(\eta)\rangle\in 
j_*j^*(\O^{\boxtimes k}).$$}
$$\nabla:j_*j^*(\pi^*_k(\conf(\Xf,A)))\to 
p_{1,\dotsc,k}{}_*(j_*j^*(\pi_k^*(\conf(\Xf,\A))
\boxtimes\Omega^2)).$$

\begin{prop}   \label{anotherchar}
Let $f\cdot a_1\boxtimes\ldots\boxtimes a_k$ be a section of 
$j_*j^*({}^l\A^{\boxtimes k})$ on $\Xf^k\setminus \Delta$. 
Then $$\nabla(\langle f\cdot a_1\boxtimes\ldots\boxtimes 
a_k\rangle)=(\langle\ldots\rangle\boxtimes\on{id})
(\nabla_{loc}(f\cdot a_1\boxtimes\ldots\boxtimes a_k))+
\langle f\cdot a_1\boxtimes\ldots\boxtimes a_k\boxtimes {}^lT\rangle.$$
\end{prop}

\begin{proof}

When $k=1$, the statement follows immediately from the definition of the action 
$\tau$.

\medskip

To treat the case of arbitrary $k$ one proceeds as follows: 

\medskip

Step 1: 

\noindent Consider the diagram
$$
\CD
j_*j^*({}^l\A^{\boxtimes k})  @>>> p_{1,\dotsc,k}{}_*(j_*j^*({}^l\A^{\boxtimes 
k}\boxtimes\Omega^2))  \\
@V{\langle\ldots\rangle}VV    @V{\langle\ldots\rangle\boxtimes\on{id}}VV  \\
j_*j^*(\pi_k^*(\conf(\Xf,\A)))  @>{?}>>  
p_{1,\dotsc,k}{}_*(j_*j^*(\pi_k^*(\conf(\Xf,\A))\boxtimes\Omega^2)),
\endCD
$$
where the upper horizontal arrow sends a section $f\cdot 
a_1\boxtimes\ldots\boxtimes a_k$
of $j_*j^*({}^l\A^{\boxtimes k})$ to
$$\nabla_{loc}(f\cdot a_1\boxtimes\ldots\boxtimes a_k)+f\cdot 
a_1\boxtimes\ldots\boxtimes a_k\boxtimes {}^lT.$$

\begin{lem}
There exists a unique map 
$$\nabla':j_*j^*(\pi_k^*(\conf(\Xf,\A)))\to 
p_{1,\dotsc,k}{}_*(j_*j^*(\pi_k^*(\conf(\Xf,\A))\boxtimes\Omega^2))$$
which completes the above diagram to a commutative diagram.
\end{lem}

\medskip

Step 2: 

\noindent It remains to show that the map $\nabla'$ constructed above coincides 
with the covariant derivative $\nabla$.
Consider the diagram
$$
\CD
j_*j^*({}^l\A\boxtimes\O^{\boxtimes k-1}) @>>> 
p_{1,\dotsc,k}{}_*(j_*j^*({}^l\A\boxtimes\O^{\boxtimes k-1}\boxtimes\Omega^2)) \\
@V{\langle\ldots\rangle\boxtimes\on{id}}VV 
@V{\langle\ldots\rangle\boxtimes\on{id}\boxtimes\on{id}}VV  \\
j_*j^*(\pi_k^*(\conf(\Xf,\A)))  @>{?}>>  
p_{1,\dotsc,k}{}_*(j_*j^*(\pi_k^*(\conf(\Xf,\A))\boxtimes\Omega^2)),
\endCD
$$
where the upper horizontal arrow is given by the same formula as above.

\smallskip

On the one hand, according to Step 1 and \secref{confsummary},
this diagram becomes commutative if we set the lower horizontal arrow to be 
$\nabla'$. On the other hand, the assertion
of our proposition for $k=1$ implies that this diagram becomes commutative if we 
set the lower horizontal arrow
to be $\nabla$. Since the left vertical arrow is a surjection (by 
\propref{confmany}), $\nabla'=\nabla$.

\end{proof}

\ssec{Introducing the central charge}  
 \sssec{}  \label{linalg}
 We shall start with some linear algebra preliminaries. 

\smallskip

Consider the canonical decreasing sequence of ideals 
in the Lie algebra $\Vir^+_0$:
$$
0\subset\ldots\subset \Vir^+_n\subset \ldots\subset 
\Vir^+_1\subset \Vir^+_0.
$$
For $k\in\NN$ let $L_k$ denote the $\Aut^+_0$-module $\Vir^+_k/\Vir^+_{k+1}$. We 
have canonical 
isomorphisms $L_k\otimes L_{k'}\simeq L_{k+k'}$ and for $k\in\NN$ we define 
$L_{-k}$ as $L_k^*$.

\begin{lem} \label{allthesame}
The following two categories are equivalent:

\smallskip

\noindent{\em (a)}
The category $\on{Ext}_{\Aut^+_0}(L_{-1},L_1)$ of extensions 
$0\to L_1\to E\to L_{-1}\to 0$ of $\Aut^+_0$-modules.

\smallskip

\noindent{\em (b)}
The category of Harish-Chandra pairs $(\Vir',\Aut^+_0)$, where $\Vir'$ is a 
central extension of $\Vir$:
$$0\to\CC\to\Vir'\to\Vir\to 0.$$

\end{lem}

In fact, both categories are groupoids such that every object has no non-trivial 
automorphisms; moreover
the set of isomorphism classes of objects in each of them is canonically 
identified with $\CC$.\footnote{These identifications are chosen as follows:

If $0\to L_1\to E\to L_{-1}\to 0$ is an object of 
$\on{Ext}_{\Aut^+_0}(L_{-1},L_1)$, the action of $\Vir^+_0$
on $E$ produces a map $L_2\otimes L_{-1}\to L_1$, i.e. a scalar $c\in\CC$, since 
$L_2\otimes L_{-1}$ is apriori 
identified with $L_1$. We set the class of this $E$ in $\CC$ to be $-2c$.

If $\Vir'$ is an object of the second category, the adjoint action of 
$\CC^*\subset \Aut^+_0$ on $\Vir'$ defines
a decomposition of the latter as a vector space: 
$\Vir'=\underset{i}\oplus L_i\oplus \CC$. The Lie bracket 
on $\Vir'$ induces a map $L_2\otimes L_{-2}\to\CC$, i.e. a scalar $c\in\CC$. We 
set the class of this $\Vir'$ 
to be equal to be $-2c$.} This implies the assertion of the lemma.

\medskip

The above equivalence of categories has the following additional property:

\smallskip

To an object $0\to L_1\to E\to L_{-1}\to 0$ of $\on{Ext}_{\Aut^+_0}(L_{-1},L_1)$ 
one can attach 
a local $\O$-module $\Theta':=\on{Ass}(E)$ on $\Xf$ (cf. \secref{anotherlocal}). 
By the construction, $\Theta'$ is an extension
$$0\to\Omega\to\Theta'\to\Theta\to 0.$$ 

\smallskip

Now let $(x\in X)$ be a point of $\Xf$. Consider the corresponding extension of 
$\hat\K_x$-modules:
$$0\to \Omega{\underset{\O}\otimes}\hat\K_x \to 
\Theta'{\underset{\O}\otimes}\hat\K_x \to 
\Theta{\underset{\O}\otimes}\hat\K_x\to 0$$ and consider the vector space
$$\Theta'{\underset{\O}\otimes}\hat\K_x/\on{ker}(\on{Res}:\Omega{\underset{\O}\otimes}\hat\K_x\to\CC).$$

\smallskip

The action of $\Theta{\underset{\O}\otimes}\hat\K_x$ on 
$\Theta'{\underset{\O}\otimes}\hat\K_x$ by Lie derivations 
(cf. \secref{Lieder}) makes the above vector space into a Lie algebra.

\begin{lem}
Every choice of an identification $\alpha:\hat\O\simeq\hat\O_x$ gives rise to an 
isomorphism 
$$\Vir'\simeq 
\Theta'{\underset{\O}\otimes}\hat\K_x/\on{ker}(\on{Res}:\Omega{\underset{\O}\otimes}\hat\K_x\to\CC),$$
where $\Vir'$ is the central extension of $\Vir$ corresponding to $E$. A 
modification of $\alpha$ by means of $g\in\Aut^+_0$
corresponds to the adjoint action of $g$ on $\Vir'$.
\end{lem}

\medskip

Let us now describe explicitly the extension $0\to\Omega\to\Theta'_2\to\Theta\to 
0$ that 
corresponds to the object $E_2\in \on{Ext}_{\Aut^+_0}(L_{-1},L_1)$ whose class in 
$\CC$ equals $2$:

\smallskip

Consider the local $\O$-module
$$p_2{}_*(\O\boxtimes\O(\Delta)/\O\boxtimes\O(-2\cdot\Delta));$$
the section $1\boxtimes 1$ gives rise to an embedding of $\O$ into it. 

\begin{lem}
The local $\O$-module $\Theta'_2$ is identified with
$$(p_2{}_*(\O\boxtimes\O(\Delta)/\O\boxtimes\O(-2\cdot\Delta)))/\O.$$
\end{lem}

\smallskip

In what follows, for $c\in\CC$ we shall denote by $\Theta_c$ (resp., $\Vir'_c$)
the corresponding extension of $\Theta$ by $\Omega$ (resp., the 
corresponding central extension of $\Vir$). The construction
of Example 2 of \secref{firstexamples} produces from $\Theta'_c$ a 
Lie-* algebra $\T'_c$ on $\Xf$ that satisfies the conditions of 
\lemref{locchiral}.

\smallskip

Let $\P_c$ denote the Serre dual to $\Theta'_c$. This is again a local $\O$-module 
on $\Xf$, which
is an extension of $\O$ by $\Omega^2$. It is easy to see that $\P_c$ identifies 
with the
$\Omega^2$-torsor of $c$-projective connections. In particular, for a fixed 
$X\in\Mf$, the extensions
$$0\to\Omega\to \Theta'_c\to\Theta\to 0 \text{ and } 0\to\Omega^2\to\P_c\to\O\to 
0$$
are (non-canonically) split.

\sssec{}  \label{central}
 A CFT of central charge $c$ consists of:

\smallskip

\noindent (1) A local $\O$-module $\A$ on $\Xf$ endowed with a structure of chiral 
algebra.

\smallskip

\noindent (2) A map of local $\O$-modules 
$T:\Theta'_c\to\A$ (called the energy-momentum tensor), such
that the composition $\Omega\to\Theta'_c\to\A$ coincides with the map $\unit$.

\medskip

The pair $(\A,T)$ must satisfy the following conditions:

\smallskip

\noindent (a) The chiral bracket on $\A$: 
$$j_*j^*(\A\boxtimes\A){\overset{\{,\}}\longrightarrow}\Delta_{!}(\A)$$
is a map of local $\O$-modules on $\Xf^2$.

\smallskip

\noindent (b) Let $\xi$ be a vertical vector field on $\Xf$ and let $\xi'$ be some 
lift of $\xi$ to a section
of $\Theta'_c$. Consider the map $\A\to\A$ given by
$$l\to (\on{id}\boxtimes h)(\{l,T(\xi')\})-\on{Lie}_\xi(l).$$
We need this map to coincide with the action of $\xi$ on $\A$ coming from the 
right D-module
structure on $\A$. (Note that the above expression a priori does not depend on the 
choice of $\xi'$.) 

\bigskip

As in the case of central charge $0$, it follows that the map $\A\otimes D\to \A$ 
defining the D-module 
structure on $\A$ is a map of local $\O$-modules and that the D-module map
$\T'_c\to\A$ induced by $T:\Theta'_c\to\A$ is a Lie-* algebra homomorphism.

\sssec{}  \label{centralalgebroid}
 Let $\tilgf'_c$ denote the Lie algebroid $\Vir'_c\otimes\O_{\hat\Xf}$ on 
$\hat\Xf$. Let 
$\tilgf_{glob}$ denote the kernel $\tilgf\to\Theta_{\hat\Xf}$; obviously 
$\tilgf_{glob}$ is an ideal
in $\tilgf$ and its fiber at a point $(X,x,\alpha)\in\hat\Xf$ 
is identified with $H^0(X\setminus x,\Theta)$.

\begin{prop}
There is a canonical lifting:
$$\tilgf_{glob}\to\tilgf'_c.$$
Moreover, $\tilgf_{glob}$ becomes an ideal in $\tilgf'_c$.
\end{prop}

\begin{proof}

To construct the embedding $\tilgf_{glob}\hookrightarrow\tilgf'_c$ we must exhibit 
for each $(X,x,\alpha)\in\hat\Xf$
as above a map
$$H^0(X\setminus x,\Theta)\to 
\Theta'_c{\underset{\O}\otimes}\hat\K_x/\on{ker}(\on{Res}:\Omega\underset{\O}\otimes\hat\K_x\to\CC).$$

\smallskip

We obviously have a map
$H^0(X\setminus x,\Theta'_c)\to \Theta'_c{\underset{\O}\otimes}\hat\K_x,$ and the 
composition
$$H^0(X\setminus x,\Omega)\to \Theta'_c{\underset{\O}\otimes}\hat\K_x\to 
\Theta'_c{\underset{\O}\otimes}\hat\K_x/\on{ker}(\on{Res}:\Omega\underset{\O}\otimes\hat\K_x\to\CC)$$
vanishes by the residue formula.

\smallskip

Since the extension $0\to\Omega\to\Theta'_c\to\Theta\to 0$ is 
(non-canonically) split,
the projection $H^0(X\setminus x,\Theta'_c)\to H^0(X\setminus x,\Theta)$ is 
surjective
and we obtain a well-defined map
$$H^0(X\setminus x,\Theta)\to \Theta'_c{\underset{\O}\otimes}\hat\K_x/\on{ker}
(\on{Res}:\Omega\underset{\O}\otimes\hat\K_x\to\CC)$$
as needed.

\medskip

To prove that $\tilgf_{glob}$ is an ideal in $\tilgf'_c$ we must show that it is 
stable under the action of
the formal group $\Aut$ on $\hat\Xf$. This follows from the description of this 
action that was given in 
\secref{intralgebroid}.

\end{proof}

Obviously, $\tilgf'_c$ is a Harish-Chandra Lie algebroid with respect to the group 
$\Aut^+_0$; let
$\gf'_c$ and $\gf'_c/\gf^+_0$ denote the corresponding Lie algebroids on $\Xf$. 
Clearly, we have short exact
sequences:
$$0\to\O\to\gf'_c\to\gf\to 0 \text{ and } 0\to\O\to\gf'_c/\gf^+_0\to\gf/\gf^+_0\to 
0$$ 
with $\O$ being an ideal and the adjoint action of $\gf$ on $\O$ being the 
canonical one. As an
$\O$-module, $\gf'_c$ is identified with
$$p_1{}_*(\O\hat\boxtimes\Theta'_c(\infty\cdot\Delta))/
\on{ker}(p_1{}_*(\O\hat\boxtimes\Omega(\infty\cdot\Delta))\to\O)$$
and similarly for $\gf'_c/\gf^+_0$.

\medskip

Now let $\gf_{glob}$ denote the kernel of $\gf\to \Theta_{\Xf}$. It follows from 
the proposition above that
$\gf_{glob}$ embeds naturally into $\gf'_c$ and is an ideal thereof.
 Moreover, the quotient 
$${\Theta'_{\Xf}}_c:=\gf'_c/(\gf^+_0+\gf_{glob})$$
is a Picard Lie algebroid on $\Xf$. 

\smallskip

Note that we have a canonical embedding $\Vir^+\to\Vir'_c$ which gives rise to an 
embedding of Lie algebroids
$\Theta\to \gf'_c/\gf^+_0$. In particular, there is a homomorphism of Lie 
algebroids
$\Theta\to {\Theta'_{\Xf}}_c$.

\begin{prop} \label{algebroiddescent}
There exists a canonical Picard Lie algebroid ${\Theta'_{\Mf}}_c$ on $\Mf$ such 
that 
${\Theta'_{\Xf}}_c$ is its pull-back (in the sense of \secref{Dmod})
under the map $\pi$. As an $\O_\Mf$-module, ${\Theta'_{\Mf}}_c$ is
identified with $R^1\pi_*(\Theta'_c)$.
\end{prop}

\begin{proof}

As an $\O$-module, the Picard Lie algebroid ${\Theta'_{\Xf}}_c$ 
is identified with
$$R^1p_1{}_*(\O\boxtimes\Theta'_c(-\Delta))/\on{ker}(R^1p_1{}_*(\O\boxtimes\Omega(
-\Delta))\to
R^1p_1{}_*(\O\boxtimes\Omega))$$
and the quotient ${\Theta'_{\Xf}}_c/\Theta$ is therefore 
identified with $\pi^*(R^1\pi_*(\Theta'_c))$.

\smallskip

Let $N(\Theta)$ denote the normalizer of $\Theta$ in ${\Theta'_{\Xf}}_c$. 

\begin{lem}
$N(\Theta)/\Theta\subset {\Theta'_{\Xf}}_c/\Theta$ corresponds to 
$\pi^{\bullet}(R^1\pi_*(\Theta'_c))\subset \pi^*(R^1\pi_*(\Theta'_c))$.
\end{lem}

Therefore, ${\Theta'_{\Mf}}_c:=R^1\pi_*(\Theta'_c)$ has a canonical structure of  
Lie algebroid on $\Mf$ and 
we have a commutative square:
$$
\CD
\Theta @>>> \Theta'_{\Xf_c} @>>> 
\pi^*(\Theta'_{\Mf_c}) @>>> 0 \\
@VVV    @VVV      @VVV  \\
\Theta @>>> \Theta_\Xf @>>> \pi^*(\Theta_\Mf) @>>> 0,
\endCD
$$
which is what we had to prove.
\end{proof}

\sssec{}
 Similarly to the case $c=0$, we shall show now that our CFT chiral algebra $^l\A$ 
carries a natural 
action of the Lie algebroid $\gf'_c/\gf^+_0$. To do this, we must construct a map
$\tau:\gf'_c/\gf_0^+\underset{\CC}\otimes{}^l\A\to{}^l\A$ such that the 
composition
$$\O\underset{\CC}\otimes{}^l\A\hookrightarrow 
\gf'_c/\gf_0^+\underset{\CC}\otimes{}^l\A\to{}^l\A$$
coincides with the tautological map $\O\underset{\CC}\otimes{}^l\A\to{}^l\A$.

\medskip

We have a map $\tau_{loc}:\gf'_c\underset{\CC}\otimes{}^l\A\to{}^l\A$ that factors 
through
$$\gf'_c\underset{\CC}\otimes{}^l\A\to\gf\underset{\CC}\otimes{}^l\A\to {}^l\A,$$
the second arrow being defined by the same formula as in \secref{intralgebroid} 
(it comes from the 
structure on $^l\A$ of a local $\O$-module).

\medskip

We have a map $\tau_{em}:\gf'_c\underset{\CC}\otimes{}^l\A\to{}^l\A$ that attaches 
to a local section
$f\cdot 1\hat\boxtimes \xi'\in 
p_1{}_*(\O\hat\boxtimes\Theta'_c(\infty\cdot\Delta))\simeq \gf'_c$ and 
a local section $a\in{}^l\A$ 
the section $$(\on{id}\boxtimes h)(\{f\cdot a\boxtimes T(\xi')\})\in{}^l\A.$$

\medskip

The difference $\tau:=\tau_{loc}-\tau_{em}$ is therefore a map 
$\gf'_c\underset{\CC}\otimes{}^l\A\to{}^l\A$ 
which factors (due to condition (b) of \secref{central}) as 
$$\gf'_c/\gf_0^+\underset{\CC}\otimes{}^l\A\to{}^l\A.$$

\smallskip

As in the previous subsection we have the following assertion:

\begin{lem}
The map $\tau:\gf'_c/\gf_0^+\underset{\CC}\otimes{}^l\A\to{}^l\A$ is a left action 
of the Lie algebroid 
$\gf'_c/\gf_0^+$ on $^l\A$. 
\end{lem}

\bigskip

Consider now the $\O_\Mf$-module $\conf(\Xf,\A)$.

\begin{prop} \label{KZ}

The $\O_\Mf$-module $\conf(\Xf,\A)$ on $\Mf$ has a (unique) structure of left
${\Theta'_{\Mf}}_c$-module such that when we view $\pi^*(\conf(\Xf,\A))$ as a 
$\gf'_c/\gf^+_0$-module via 
$\gf'_c/\gf^+_0\to{\Theta'_{\Xf}}_c$,
the projection
$$\langle\ldots\rangle:{}^l\A\to \pi^*(\conf(\Xf,\A))$$
is compatible with the $\gf'_c/\gf^+_0$--action.

\end{prop}

The proof of this assertion is the same as in the case $c=0$.

\bigskip

For $k\in\NN$ let ${\Theta'_{\Xf_c^k}}$ denote the pull-back of the Picard Lie 
algebroid ${\Theta'_\Mf}_c$ to
$\Xf^k$ (in the sense of \secref{Dmod});
by definition, for $k=1$ we re-obtain on $\Xf$ the Picard Lie-algebroid 
${\Theta'_{\Xf}}_c$. Note that
the corresponding twisted cotangent bundle on $\Xf^k\setminus\Delta$ is identified 
with
$$\on{ker}(p_{1,\dotsc,k}{}_*(\O^{\boxtimes k}\boxtimes \P_c(\Delta))\to 
p_{1,\dotsc,k}{}_*(\O^{\boxtimes k+1}(\Delta))).$$

\smallskip

In particular, the ${\Theta'_{\Xf_c^k}}$-action on 
$\pi_k^*(\conf(\Xf,\A))$ yields a covariant derivative map 
$$\nabla:j_*j^*(\pi_k^*(\conf(\Xf,\A)))\to 
p_{1,\dotsc,k}{}_*(j_*j^*(\pi_k^*(\conf(\Xf,\A))\boxtimes \P_c)).$$

\smallskip

The following is a generalization to the case $c\neq 0$ of \propref{anotherchar}:

\begin{prop}
Let $f\cdot a_1\boxtimes\ldots\boxtimes a_k$ be a local section of 
$j_*j^*(^l\A^{\boxtimes k})$ on $\Xf^k\setminus\Delta$.
Then the section 
$\nabla(f\cdot a_1\boxtimes\ldots\boxtimes a_k)$ of 
$j_*j^*(\pi_k^*(\conf(\Xf,\A))\boxtimes \P_c)$
equals
$$(\langle\ldots\rangle\boxtimes\on{id})(\nabla_{loc}(f\cdot 
a_1\boxtimes\ldots\boxtimes a_k))+
\langle f\cdot a_1\boxtimes\ldots\boxtimes a_k\boxtimes {}^lT\rangle,$$
where the insertion of $^lT$ in the correlation function has the same meaning as 
in \propref{anotherchar}.
\end{prop}

\newpage

\centerline{\bf Chapter III. Examples} 

 \ssec{Heisenberg and Kac-Moody algebras}
\sssec{}  \label{KM}
 Let $\gff$ and $Q$ be as in Example 1 of \secref{firstexamples} and let 
$\B(\gff,Q)$ be the corresponding Lie-* algebra; 
consider its chiral universal enveloping algebra $\U(\B(\gff,Q))$. 

\smallskip

There are two maps of $\Omega$ into $\U(\B(\gff,Q))$: one is the map 
$\unit:\Omega\to \U(\B(\gff,Q))$ and the other
is the composition
$${\widetilde\unit}:\Omega\to\B(\gff,Q)\to \U(\B(\gff,Q)).$$

We define the chiral algebra $\A(\gff,Q)$ as a quotient of $\U(\B(\gff,Q))/I$, 
where $I$ is the ideal\footnote{An ideal in a chiral algebra $\A$ is
by definition a chiral submodule of $\A$, when the latter is viewed as a module 
over itself. For a D-submodule
of $\A$ there exists a minimal ideal that contains it, 
called the ``ideal generated by the given D-sub-module''. If
$I\subset \A$ is an ideal, the quotient D-module $\A/I$ has a natural structure of 
chiral algebra.}
 generated by the
image of 
$\unit-{\widetilde\unit}\colon \Omega
\longrightarrow\U(\B(\gff,Q))$.

\medskip

We are going to show now that under certain conditions, $\A(\gff,Q)$ defines a CFT 
data.

\smallskip

It is clear that $\B(\gff,Q)$ has a natural structure of local $\O$-module on 
$\Xf$, that satisfies the conditions 
of \lemref{locchiral}. 
Hence, $\A(\gff,Q)$ is a chiral algebra on $\Xf$ 
which has a structure of local $\O$-module.
Therefore, in order for $\A(\gff,Q)$ to define a CFT structure, we must construct 
a map of local $\O$-modules
$$T_{\gff,Q}:\Theta'_c\to \A(\gff,Q)$$
for some $c\in\CC$ that satisfies condition (b) of \secref{central}. 

\smallskip

All our constructions will be local, therefore, we will be working on a fixed 
curve $X$.

\medskip

Let $B\in \gff\otimes \gff$ be an $\on{ad}_{\gff}$-invariant symmetric tensor that 
satisfies the following conditions:

\begin{itemize}

\item If $B=\underset{i,j}\Sigma \,b_{i,j}\cdot u_i\otimes u_j$, then the 
endomorphism of $\gff$ given by
$$v\to \underset{i,j}\Sigma \,b_{i,j}\cdot[u_i,[u_j,v]]$$
is a scalar operator $\kappa(B)\cdot \on{Id}_{\gff}$.

\item The map $B\circ Q:\gff\to \gff^*\to \gff$ is a scalar operator that 
satisfies
$$2(B\circ Q)+\kappa(B)=\on{Id}_{\gff}.$$

\end{itemize}

\medskip

First, we define a map $\widetilde T_{\gff,Q}:\Theta'_2\to \A(\gff,Q)$ as follows:

\smallskip

\noindent Recall that in \secref{central}, $\Theta'_2$ was constructed as a 
quotient of the $\O$-module\break 
$p_2{}_*(\O\boxtimes\O(\Delta))$. Let $\widetilde{\widetilde T}_{\gff,Q}$ be the 
map
$$p_2{}_*(\O\boxtimes\O(\Delta))\to \A(\gff,Q)$$ that sends a section $f(x,y)\in 
\O\boxtimes\O(\Delta)$ to
$$(h\boxtimes\on{id})(\underset{i,j}\Sigma \,b_{i,j}\cdot\{u_i\boxtimes u_j\cdot 
f(x,y)\})\in \A(\gff,Q).$$

\begin{lem}
For any symmetric $B\in \gff\otimes \gff$, the map
$\widetilde{\widetilde T}_{\gff,Q}: p_2{}_*(\O\boxtimes\O(\Delta))\to \A(\gff,Q)$ 
factors as
$$p_2{}_*(\O\boxtimes\O(\Delta))\twoheadrightarrow \Theta'_2\overset{\widetilde 
T_{\gff,Q}}\longrightarrow \A(\gff,Q).$$
Moreover, the composition
$$\Omega\hookrightarrow \Theta'_2\overset{\widetilde 
T_{\gff,Q}}\longrightarrow\A(\gff,Q)$$
coincides with $\on{Tr}(B\circ Q)\cdot\unit$.
\end{lem}

\smallskip

We now set $c=2\on{Tr}(B\circ Q)$ and define the map $T_{\gff,Q}:\Theta'_c\to 
\A(\gff,Q)$ as a composition
$$\Theta'_c\overset{?}\to \Theta'_2{\overset{\widetilde 
T_{\gff,Q}}\longrightarrow}\A(\gff,Q),$$
where the first arrow is the unique isomorphism that fits into the commutative 
diagram
$$
\CD
0 @>>> \Omega @>>> \Theta'_c @>>> \Theta @>>> 0 \\
@VVV     @V{\frac{1}{\on{Tr}(B\circ Q)}}VV    @V{?}VV   @V{\on{id}}VV     @VVV  \\
0 @>>> \Omega @>>> \Theta'_2 @>>> \Theta @>>> 0.
\endCD
$$

\smallskip

Clearly, the map $T_{\gff,Q}$ is compatible with the local $\O$-module structure.

\begin{prop}  \label{sugawara}
For $B$ and $c$ as above, the map $T_{\gff,Q}:\Theta'_c\to \A(\gff,Q)$ satisfies 
condition (b) of the CFT data.
\end{prop}

Before giving a proof of this proposition, let us consider two examples:

\medskip

\noindent{\bf Example 1} 

\smallskip

Let $\gff=\hf$ be abelian; in this case $B$ exists if and only if $Q$ is 
non-degenerate and $B=(2\cdot Q)^{-1}$. 
We have $c=\on{dim}(\hf)$.

\medskip

\noindent{\bf Example 2}

\smallskip

Now let $\gff$ be semisimple. The condition for the existence of $B$ is that 
$\kappa(Q^{-1})+2\neq 0$, i.e.
that $Q$ is not equal to $-\frac12\,K$, where $K$ is the Killing form. 

\smallskip

Therefore, in this case the construction we described above is just an invariant 
way to write down
Sugawara's construction. Moreover, if $Q=q\cdot K$
we see that the central charge in this case equals $q\cdot\on{dim}(\gff)/(q+1/2)$.

\sssec{}
 Let us now prove \propref{sugawara}:

\begin{proof}

\medskip

We have to show that for a section $f(x,y)$ of $\O\boxtimes\O(\Delta)$ and for any 
section
$l\in\A(\gff,Q)$, we have:
$$(h\boxtimes\on{id})(\{T_{\gff,Q}(\xi_f),l\})=-\on{Lie}_{\xi_f}(l)-l\cdot 
\xi_f,$$
where $\xi_f$ is a vector field on $X$ corresponding to $f(x,y)$ via the 
identification
$$\O\boxtimes\O(\Delta)/\O\boxtimes\O\simeq\Theta.$$

\smallskip

First of all, since $\B(\gff,Q)$ generates $\A(\gff,Q)$, the
Jacobi identity implies that it is enough to prove the
above equality for $l\in \B(\gff,Q)\subset \A(\gff,Q)$. 
Moreover, we can assume that $l$ is of the form $l=v\otimes 1\in \gff\otimes D_X$.
Therefore, we have to prove the equality
$$(h\boxtimes h\boxtimes\on{id})
(\underset{i,j}\Sigma \,b_{i,j}\{\{u_i\boxtimes u_j\cdot 
f(x,y)\},v\})=-v\otimes\xi_f\in\B(\gff,Q)
\subset\A(\gff,Q).$$

\smallskip

It is easy to see that the assertion holds for $f(x,y)\in \O\boxtimes\O$, 
therefore, we can assume that
$f(x,y)=-f(y,x)$. Using the fact that the tensor $\underset{i,j}\Sigma 
\,b_{i,j}\cdot u_i\otimes u_j$ is symmetric
and the Jacobi identity, the above expression can be rewritten as:
$$\underset{i,j}\Sigma\, 2b_{i,j}\cdot(h\boxtimes h\boxtimes\on{id})
(\{u_i\boxtimes\{f(x,y)\cdot u_j\boxtimes v\}\}).$$

By the definition of the Lie-* bracket on $\B(\gff,Q)$, the last expression equals
\begin{align*}
&\underset{i,j}\Sigma\, 2b_{i,j}\cdot(h\boxtimes\on{id})(\{f(x,y)\cdot 
u_i\boxtimes [u_j,v]\})+\\
&+\underset{i,j}\Sigma\, 2b_{i,j}\cdot Q(u_j,v)\cdot(h\boxtimes h\boxtimes\on{id})
(\{f(x,y)\cdot u_i\boxtimes \one'_{y=z}\}),
\end{align*}
where $\one'_{y,z}$ denotes the canonical antisymmetric section of 
$\Delta_{!}(\Omega)\subset \Delta_{!}(\A(\gff,Q))$.

\medskip

To prove the proposition we must show that 

\smallskip

\noindent a) 
$\underset{i,j}\Sigma\, 2b_{i,j}\cdot(h\boxtimes\on{id})(\{f(x,y)\cdot 
u_i\boxtimes [u_j,v]\})=
-\underset{i,j}\Sigma \,b_{i,j}\cdot [u_i,[u_j,v]]\otimes\xi_f$ and 

\smallskip

\noindent b)
that for $u\in \gff$, $(h\boxtimes h\boxtimes\on{id})(\{f(x,y)\cdot u\boxtimes 
\one'_{y,z}\})=-u\otimes\xi_f$.

\medskip

Using the fact that the tensor $\underset{i,j}\Sigma\, b_{i,j}\cdot u_i\otimes 
u_j$ is $\on{ad}_{\gff}$-invariant,
the expression $$\underset{i,j}\Sigma\, 
2b_{i,j}\cdot(h\boxtimes\on{id})(\{f(x,y)\cdot u_i\boxtimes [u_j,v]\})$$
can be rewritten as
$$\underset{i,j}\Sigma\, b_{i,j}\cdot(h\boxtimes\on{id})(\{f(x,y)\cdot 
u_i\boxtimes [u_j,v]\}-
\{f(x,y)\cdot [u_j,v]\boxtimes u_i\})$$
and assertion a) follows from the next lemma:

\begin{lem} \label{funnycommutator}
Let $\A$ be a chiral algebra and let $\underset{i,j}\Sigma \,a^i_1\boxtimes a^j_2$ 
be a section of $\A\boxtimes\A$. 
Assume that $\underset{i,j}\Sigma\,\{a^i_1\boxtimes 
a^j_2\}\in\A\subset\Delta_{!}(\A)$. Now let $f(x,y)$ 
be an antisymmetric section of $\O\boxtimes\O(\Delta)$. We have:
$$\underset{i,j}\Sigma\,(\{a^i_1\boxtimes a^j_2\cdot f(x,y)\}-\{a^j_2\boxtimes 
a^i_1\cdot f(x,y)\})=
\underset{i,j}\Sigma\,\{a^i_1\boxtimes a^j_2\}\cdot 
(\xi_f,-\xi_f)\in\Delta_{!}(\A),$$
where $\xi_f$ has the same meaning as above.
\end{lem}

\medskip

The expression
$(h\boxtimes h\boxtimes\on{id})
(\{f(x,y)\cdot u\boxtimes \one'_{y=z}\})$ equals
$$(h\boxtimes\on{id})(\{u\boxtimes d_y(f(x,y))\})$$
and assertion b) follows from the following general fact:

\begin{lem} \label{funnyresidue}
Let $\M$ be a right D-module on a curve $X$, let $m\in \M$ be its local section 
and let $f(x,y)$ be a section of 
$\O\boxtimes\O(\Delta)$. Then under the map 
$$(h\boxtimes\on{id})\circ\on{can}_\M:j_*j^*(\M\boxtimes\Omega)\to 
\Delta_{!}(\M)\to\M$$
the section $m\boxtimes d_y(f(x,y))$ goes over to $-m\cdot\xi_f$.
\end{lem}

\end{proof}

\ssec{The linear dilaton}  
\sssec{}  \label{dilaton}
 Now let $\gff=\CC$. We are going to consider 
a twist of the CFT corresponding to $\A(\CC,Q)$, called the linear dilaton.

\smallskip

We start with the following observation:

\begin{lem}
Consider the category $\on{Ext}_{\Aut^+_0}(L_0,L_1)$
of extensions 
$$0\to L_1\to E\to L_0\to 0$$ of $\Aut^+_0$-modules.
We have a canonical isomorphism\footnote{We fix this isomorphism
in such a way that if $0\to L_1\to E\to L_0\to 0$ is an object of 
$\on{Ext}_{\Aut^+_0}(L_0,L_1)$ corresponding
to $\lambda\in\CC$, the map $L_1\simeq L_1\otimes L_0\to L_1$ induced by the 
action of $L_1\subset\Vir^+_0$ on $E$ 
is multiplication by $\lambda$.} $\pi_0(\on{Ext}_{\Aut^+_0}(L_0,L_1))\simeq\CC$.
\end{lem}

Thus, for $\lambda\in\CC$ one obtains a canonical extension $F_\lambda$:
$$0\to\Omega\to F_\lambda\to\O\to 0$$
in the category of local $\O$-modules on $\Xf$. For instance, if $\lambda=-2n$, 
where $n$ is an integer, we have
a canonical isomorphism
$$F_\lambda\simeq \on{Diff}_1(\Omega^n,\Omega^{n+1}).$$

\medskip

For a pair of complex numbers $(\lambda,Q)$ we now define a Lie-* algebra 
structure on the D-module
$$\B(\lambda,Q):=F_\lambda \otimes D/\on{ker}(\Omega\otimes D\to\Omega)$$
by the rule that the Lie-* bracket is the composition
$$\B(\lambda,Q)\boxtimes\B(\lambda,Q)\to D\boxtimes D \to \Delta_{!}(\Omega),$$
where the last arrow sends $1\boxtimes 1\in D\boxtimes D$ to $Q\cdot \one'\in 
\Delta_{!}(\Omega)$.

\smallskip

We define the chiral algebra $\A(\lambda,Q)$ as the quotient of 
$\U(\B(\lambda,Q))$ by the standard relation
$\unit=\widetilde\unit$. 

\medskip

Here is another useful point of view on the chiral algebra $\A(\lambda,Q)$:

\begin{prop}  \label{twistdilaton}
Consider the chiral algebra $\A(\CC,Q)$ ($Q$ is viewed as a quadratic form on 
$\CC$).

\smallskip

\noindent {\em (a)} There is a natural map of sheaves $\Omega\to 
\on{Aut}(\A(\CC,Q))$.

\smallskip

\noindent {\em (b)} The chiral algebra $\A(\lambda,Q)$ is
identified naturally with 
a twist of $\A(\CC,Q)$ with respect to the $\Omega$-torsor $F_\lambda$.
\end{prop}

\begin{proof}

For a section $\omega\in\Omega$ we define the corresponding automorphism of 
$\B(\CC,Q)$ by the rule that it
acts as the identity on $\Omega\subset \B(\CC,Q)$ and that it sends the section 
$1\in D\in\B(\CC,Q)$ 
to $\omega\in\Omega\subset\B(\CC,Q)$; clearly, it extends to an automorphism of 
$\A(\CC,Q)$. This establishes point (a)
of the proposition, whereas point (b) follows immediately from the definitions.

\end{proof}

\sssec{}  \label{stressdilaton}
 We shall now show that the chiral algebra $\A(\lambda,Q)$ defines a CFT with 
$c=1+3\lambda^2/Q$.

\smallskip

Obviously, $\A(\lambda,Q)$ carries a structure of local $\O$-module on $\Xf$ which 
is compatible
with the chiral bracket. Thus, to define the corresponding CFT data we must 
exhibit a map of local $\O$-modules
$T_{\lambda,Q}:\Theta_c\to\A(\lambda,Q)$ with 
$c$ as above which will satisfy condition (b) of \secref{central}.

\begin{prop}  \label{almoststress}
There exists a unique map $\widetilde T_{\lambda,Q}:\Theta\to 
\A(\lambda,Q)/\Omega$ such that the following
holds:

\smallskip

\noindent For any local section $l\in \A(\lambda,Q)$, any vector field $\xi$ on 
$X$ and any section $\xi'$
of $\A(\lambda,Q)$ that projects to $\widetilde T_{\lambda,Q}(\xi)$ we have:
$$(h\boxtimes\on{id})(\{\xi',l\})=-\on{Lie}_\xi(l)-l\cdot\xi.$$

\end{prop}

\begin{proof}

Note first of all that the left hand side of the above expression is a priori 
independent of the choice of $\xi'$.
Moreover, the uniqueness statement is clear, since $\Omega$ is the only 
$\O$-submodule in $\A(\lambda,Q)$
that annihilates the Lie-* bracket (this follows, for example, from the 
corresponding fact in the untwisted case
(i.e. when $\lambda=0$) in view of \propref{twistdilaton}(b)).

\smallskip

The existence assertion is local and let us choose a coordinate $\z$ on some open 
sub-set $U\subset X$; 
let $\partial_\z$
denote the corresponding vector field on $U$. It is clear that the choice of $\z$ 
as above trivializes the
$\Omega$-torsor $F_\lambda$ over $U$; let $\phi_\z$ denote the corresponding
isomorphism
$$\phi_\z:\A(\CC,Q)\to\A(\lambda,Q)$$
defined over $U$.

\smallskip

We define the map 
$\widetilde T_{\lambda,Q}$ by the condition that 
$$
\widetilde T_{\lambda,Q}(\partial_\z)=
\phi_\z(T_{\CC,Q}(\partial'_\z)-\lambda/2Q\cdot\partial_\z),
$$
where $\partial'_\z$ in the first term is any lifting of $\partial_\z$ to a 
section of $\Theta'_1$ and
$\partial_\z$ in the second term is viewed as a section of $D_X\subset\B(\CC,Q)$.

\medskip

Let us check now that the map $\widetilde T_{\lambda,Q}$ defined in the above way 
satisfies the condition
of the proposition.

\smallskip

Let $f$ be a function on $U$. It is enough to show that
$$(h\boxtimes\on{id})
(\{\phi_\z(T_{\CC,Q}(\partial'_\z)\cdot f-\lambda/2Q\cdot\partial_\z f)\boxtimes
\phi_\z(1)\})=-\phi_\z(1)\cdot f\partial_\z-\on{Lie}_{f\partial_\z}(\phi_\z(1)).$$

\smallskip

Since $\phi_\z$ is a homomorphism of chiral algebras and since $T_{\CC,Q}$ is the 
energy momentum tensor for 
$\A(\CC,Q)$, we have:
$$(h\boxtimes\on{id})(\{\phi_\z(T_{\CC,Q}(\partial'_\z)\cdot 
f)\boxtimes\phi_\z(1)\})=
-\phi_\z(1)\cdot f\partial_\z$$
and the above assertion reduces to the equality
$$\lambda\cdot d(\partial_\z(f))=2\on{Lie}_{f\partial_\z}(\phi_\z(1)),$$
which follows from the definition of $F_\lambda$.

\end{proof}

The map $\widetilde T_{\lambda,Q}$ constructed above is compatible with the local 
$\O$-module structure
due to the uniqueness part of the assertion of \propref{almoststress}.
We now define $\Theta'$ to be the extension $0\to \Omega\to\Theta'\to\Theta\to 0$ 
induced from the extension
$0\to \Omega\to \A(\lambda,Q)\to \A(\lambda,Q)/\Omega\to 0$ by means of the map 
$\widetilde T_{\lambda,Q}$. 

\smallskip

By the construction, $\Theta'$ carries in a natural way a structure of local 
$\O$-module on $\Xf$
and we have a map $T_{\lambda,Q}:\Theta'\to \A(\lambda,Q)$ that satisfies 
condition (b) of \secref{central}.

\begin{prop}  \label{dilatonstress}
The local $\O$-module $\Theta'$ constructed above is
identified with $\Theta_c$ for $c=1+3\lambda^2/Q$.
\end{prop}

\begin{proof}

The proof of the proposition will be based on the following lemma:

\smallskip

\begin{lem} \label{compofc}
Let $0\to\Omega\to\Theta'\to\Theta\to 0$ be an arbitrary extension in the category 
of local $\O$-modules on $\Xf$.
Assume that for any pair $(x\in X)$ the following holds: if 
$\z$ is a coordinate around $x$ and if $\widetilde{\partial_\z}$ is a 
local section of $\Theta'$ that projects to $\partial_\z$, the value of 
$\on{Lie}_{\z^3\partial_\z}(\partial'_\z)$ at $x$
equals $-c/2\cdot d\z$. Then $\Theta'$ is canonically isomorphic to $\Theta'_c$.
\end{lem}

\smallskip

Let $U$, $\z$ and $\phi_\z$ be as in \propref{almoststress}. 
For a function $f$ on $U$ consider the map from $\A(\CC,Q)$ to itself given by
$$l\to\phi_\z^{-1}\circ\on{Lie}_{f\partial_\z}\circ\phi_\z(l)-\on{Lie}_{f\partial_
\z}(l).$$
We claim that it coincides with the map
$$l\to \lambda/2Q(h\boxtimes\on{id})(\{\partial_\z\cdot f\boxtimes l\}).$$
Indeed, both maps are derivations of the chiral algebra structure on $\A(\CC,Q)$ 
(in particular, they commute 
with the action of $D_X$) and coincide on the generators.

\smallskip

Let us now apply \lemref{compofc} to the extension $\Theta'$ as in the formulation 
of the proposition. By the
construction, we can take $\widetilde{\partial_\z}$ to be
$$\phi_\z(T_{\CC,Q}(\partial'_\z)-\lambda/2Q\cdot\partial_\z)\in \A(\lambda,Q),$$
where $\partial'_\z$ is the corresponding section of $\Theta'_1$.

\smallskip

We have:
\begin{align*}
&\phi_\z^{-1}(\on{Lie}_{\z^3\partial_\z}(\widetilde{\partial_\z}))=\lambda/2Q(h\boxtimes\on{id})
(\{\partial_\z\cdot \z^3\boxtimes T_{\CC,Q}(\partial'_\z)\})-  \\
&-\lambda/2Q(h\boxtimes\on{id})
(\{\partial_\z\cdot \z^3\boxtimes \lambda/2Q\cdot\partial_\z\})+
\on{Lie}_{\z^3\partial_\z}
(T_{\CC,Q}(\partial'_\z))-\on{Lie}_{\z^3\partial_\z}(\lambda/2Q\cdot\partial_\z).
\end{align*}

\smallskip

It is easy to see that the values at $x$ of the first and of the 
fourth terms of the above formula are $0$. The value at $x$ of the second term is 
$-3\lambda^2/2Q\cdot d\z$,
by definition. Since the chiral algebra $\A(\CC,Q)$ defines a CFT of central 
charge $1$,
the value at $x$ of the third term is $-\frac12\,d\z$. 

\end{proof}

\ssec{The bc-system} 
\sssec{}  \label{bc}
 Let $\M$ be a locally free finitely generated D-module on a curve $X$ (e.g. 
$\M=\L\otimes D_X$, where
$\L$ is a torsion-free quasi-coherent sheaf). Let $\M^*$ denote the Verdier-dual 
D-module, i.e.
$$\M^*:=\on{Hom}_{D_X}(\M,D_X\otimes \Omega).$$

\smallskip

Consider the D-module $\B_{bc}(\M):=\Omega\oplus \M\oplus\M^*$. We define a 
$\ZZ$-grading (and hence a $\ZZ_2$-grading)
on $\B_{bc}(\M)$ by declaring that $\Omega$ is homogeneous
of degree $0$ and $\M$ and $\M^*$ are of degrees $-1$ and 
$1$ respectively. We introduce a (super) Lie-* algebra structure on $\B_{bc}(\M)$ 
as follows:

\smallskip

\noindent $\Omega$ is central and the map $\M\boxtimes\M^*\to 
\Delta_{!}(\B_{bc}(\M))$ is the composition
$$\M\boxtimes\M^*\to D_X\otimes \Omega\simeq\Delta_{!}(\Omega)\hookrightarrow 
\Delta_{!}(\B_{bc}(\M)).$$

\smallskip

The (super) chiral algebra $\A_{bc}(\M)$ is by definition the quotient of 
$\U(\B_{bc}(\M))$ modulo the standard relation $\unit-\widetilde\unit=0$. 

\smallskip

The chiral algebra $\A_{bc}(\M)$ carries a $\ZZ$-grading 
$\A_{bc}(\M)=\underset{i}\oplus \A_{bc}(\M)^i$
that comes from the $\ZZ$-grading on $\B_{bc}(\M)$ and a filtration 
$\A_{bc}(\M)=\underset{i\in\ZZ^+}\cup
\A_{bc}(\M)_i$ that comes from the canonical filtration on $\U(\B_{bc}(\M))$ 
(cf. \secref{unenv}).

\medskip

In this subsection we will compute the space of conformal blocks of 
$\A_{bc}(\M)$ (assuming that $X$ is complete).

\smallskip

Note that $DR^{-1}(X,\M)=0$, since $\M$ is locally free as a D-module. Let us
denote by $\on{det}DR(\M)$ the
$1$-dimensional vector space 
$$\Lambda^{top}(DR^0(X,\M))\otimes(\Lambda^{top}(DR^1(X,\M)))^{-1}.$$ 

\smallskip

\begin{prop}  \label{confbc}
The space $\conf(X,\A_{bc}(\M))$ is $1$-dimensional. It is
concentrated in the homogeneity degree 
$\on{dim}(DR^0(X,\M))-\on{dim}(DR^1(X,\M))$ 
and is canonically isomorphic
to $(\on{det}DR(\M))^{-1}$.
\end{prop}

\begin{proof}

Let $x$ be a point of $X$. We shall compute the space $\conf(X,\A_{bc}(\M))$ 
using \propref{actualcompuatations}, which implies that
$$\conf(X,\A_{bc}(\M))=\on{coker}(DR^0(X\setminus x,\M\oplus \M^*)\otimes 
A_{bc}(\M)_x)\to
A_{bc}(\M)_x).$$

\medskip

The space $DR^0(\Spec(\hat\K_x),\M)$ is a topological vector space of Tate type 
(cf. \cite{BFM}) and we have a canonical isomorphism
$$DR^0(\Spec(\hat\K_x,\M^*))\simeq (DR^0(\Spec(\hat\K_x),\M)))^*.$$

\smallskip

Therefore, $V:=DR^0(\Spec(\hat\K_x),\M\oplus\M^*)$ acquires
a non-degenerate quadratic form such that
$$V_1:=DR^0(\Spec(\hat\O_x),\M\oplus \M^*) \text{ and } V_2:=DR^0(X\setminus 
x,\M\oplus \M^*)$$
are two Lagrangian subspaces of $V$, which are compact and co-compact,
respectively.

\smallskip 

We have the following general statement:

\smallskip

\begin{lem}  \label{spin}
Let $V$ be a Tate topological vector space endowed with a symmetric 
nondegenerate
pairing $Q:V\otimes V\to\CC$, let $V_1$ 
be a Lagrangian compact subspace of $V$ and let $S$ be the corresponding 
representation of the Clifford 
algebra. Let $V_2\subset V$ be a Lagrangian co-compact subspace of $V$. 
We have a canonical isomorphism
$$\on{coker}(V_2\otimes S\to S)\simeq\Lambda^{top}(V/V_1+V_2).$$
\end{lem}

This implies the existence of an isomorphism as in the statement of the 
proposition
using Verdier duality: $$DR^1(X,\M^*)\simeq (DR^0(X,\M))^*.$$

\medskip

The fact that the isomorphism 
$$\conf(X,\A_{bc}(\M))\simeq (\on{det}DR(\M))^{-1}$$ does not depend on the choice 
of $x\in X$ is easy to show
by generalizing the above argument to the case of several points 
$x_1,\dotsc,x_n\in X$ and by establishing
the appropriate compatibility when one of them is deleted.

\end{proof}

\smallskip

\lemref{spin} above implies also the following statement:

\smallskip

\noindent Let $x_1,\dotsc,x_{k_1},y_1,\dotsc,y_{k_2}$ be distinct points of $X$ 
and consider the composition:
$$\underset{i=1}{\overset{k_1}\otimes} M_{x_i}\underset{j=1}{\overset{k_2}\otimes} 
M_{y_j}\to
\underset{i=1}{\overset{k_1}\otimes} A_{bc}(\M)_{x_i} 
\underset{j=1}{\overset{k_2}\otimes} A_{bc}(\M^*)_{y_j}\to
\conf(X,\A_{bc}(\M))\simeq (\on{det}DR(\M))^{-1}.$$
(Here $M_x$ denotes the fiber of $\M$ at $x\in X$ and similarly for $\M^*$.)

\begin{cor} \label{bccorrelator}
The above map 
$M_{x_1}\otimes \ldots \otimes M_{x_{k_1}}\otimes M^*_{x_1}\otimes \ldots \otimes 
M^*_{x_{k_2}}\to$\break 
$(\on{det}DR(\M))^{-1}$ is zero unless $k_1=\on{dim}(DR^1(X,\M))$ and 
$k_2=\on{dim}(DR^0(X,\M))=\on{dim}(DR^0(X,\M^*))$
and in the latter case it coincides with the composition
\begin{align*}
&\underset{i=1}{\overset{k_1}\otimes} M_{x_i}\underset{j=1}{\overset{k_2}\otimes} 
M_{y_j}\to
DR^1(X,\M)^{\otimes k_1}\otimes DR^1(X,\M^*)^{\otimes k_2}\to \\
&\Lambda^{k_1}(DR^1(X,\M))\otimes \Lambda^{k_2}(DR^1(X,\M^*))\simeq 
(\on{det}DR(\M))^{-1}
\end{align*}
\end{cor}

\sssec{}
 Let $\A_{bc}(n)$ denote the chiral algebra $\A_{bc}(\M)$ for $\M=\Omega^n\otimes 
D$. In this subsection
we shall show that it defines a CFT with central charge equal to 
$c_n:=-2(6n^2-6n+1)$. 
In the physics literature this CFT is often referred to as the $bc$-system.

\medskip

According to \lemref{locchiral} $\A_{bc}(n)$ carries a structure of local 
$\O$-module on $\Xf$ that satisfies
condition (a) of \secref{central}. Therefore, to define the corresponding CFT we 
must construct
the energy-momentum tensor $T_{bc,n}:\Theta'_{c_n}\to\A_{bc}(n)$ that satisfies 
condition (b) of \secref{central}.

\smallskip

\begin{prop} \label{almoststressbc}
There exists a unique map $\widetilde T_{bc,n}:\Theta\to\A_{bc}(n)/\Omega$ such 
that for any
section $l\in\A_{bc}(n)$, any vector field $\xi$ and a section $\xi'$ of 
$\A_{bc}(n)$ that projects to
$\widetilde T_{bc,n}(\xi)$ we have:
$$(h\boxtimes\on{id})(\{\xi',l\})=-l\cdot\xi-\on{Lie}_\xi(l).$$
\end{prop}

\begin{proof}

It is easy to see that $\Omega$ is the unique $\O$-sub-module of $\A_{bc}(n)$ 
which annihilates the Lie-*
bracket. This implies the uniqueness part in the assertion of 
\propref{almoststressbc}. 

\smallskip

Therefore, the existence part is local; let $\z$ be a coordinate over some 
$U\subset X$. 

\smallskip

For $f(x,y)\in\O\boxtimes\O(\Delta)$, let $\xi_f$ denote the corresponding vector 
field and 
set $\widetilde T_{bc,n}(\xi_f)$ to be equal to the image of
$$\xi'_f:=-(h\boxtimes \on{id})(n\{d\z^{\otimes n}\cdot\partial_\z\boxtimes 
d\z^{\otimes 1-n}\cdot f(x,y)\}
-(n-1)\{d\z^{\otimes n}\boxtimes d\z^{\otimes 1-n}\cdot\partial_\z\cdot 
f(x,y)\})$$
under the projection $\A_{bc}(n)\to \A_{bc}(n)/\Omega$.

\smallskip

It is easy to see that the above expression is well-defined 
(i.e. does not depend on the choice of $f(x,y)$ that gives rise to $\xi_f$)
and that, moreover, it equals the projection to $\A_{bc}(n)/\Omega$ of
$$-(\on{id}\boxtimes h)(n\{d\z^{\otimes n}\cdot\partial_\z\boxtimes d\z^{\otimes 
1-n}\cdot f(x,y)\}
-(n-1)\{d\z^{\otimes n}\boxtimes d\z^{\otimes 1-n}\cdot\partial_\z\cdot 
f(x,y)\}).$$

\smallskip

Let us check that the map $\widetilde T_{bc,n}$ constructed above satisfies the 
condition of the proposition.
It is enough to check that for $\xi_f$ and $\xi_f'$ as above
$$(h\boxtimes\on{id})(\{\xi_f',l\})=-l\cdot\xi_f-\on{Lie}_{\xi_f}(l)$$
for $l$ being a nonvanishing section of either $\Omega^n$ or $\Omega^{1-n}$. 
As the situation is symmetric in $n$ and $1-n$ we shall treat the case of 
$l=d\z^{\otimes n}\in\Omega^n\subset \A_{bc}(n)$.

\medskip

By applying the Jacobi identity we obtain:
\begin{align*}
&(h\boxtimes\on{id})(\{\xi_f',l\})=
-n(h\boxtimes h\boxtimes\on{id})(\{d\z^{\otimes n}\cdot\partial_\z \boxtimes 
\{d\z^{\otimes 1-n}
\boxtimes d\z^{\otimes n}\cdot f(x,y)\}\})+ \\
&(n-1)(h\boxtimes h\boxtimes\on{id})(\{d\z^{\otimes n}\boxtimes 
\{d\z^{\otimes 1-n}\cdot\partial_\z\boxtimes d\z^{\otimes n}\cdot f(y,x)\}\}).
\end{align*}

\smallskip

An easy computation shows that the first term in the above expression equals 
$-ndz^{\otimes n}\cdot\xi_f-\on{Lie}_{\xi_f}(dz^{\otimes n})$. The second term 
equals $(n-1)dz^{\otimes n}\cdot\xi_f$, 
according to \lemref{funnyresidue}. This implies the assertion of the proposition.

\end{proof}

\smallskip

The uniqueness assertion in \propref{almoststressbc} guarantees that the map 
$\widetilde T_{bc,n}:\Theta\to\A_{bc}(n)/\Omega$ is compatible with the 
$\O$-module structure.
We define the local $\O$-module $\Theta'(n)$ to be the extension
$$0\to \Omega\to\Theta'(n)\to\Theta\to 0$$
induced from the extension
$0\to\Omega\to \A_{bc}\to\A_{bc}(n)/\Omega\to 0$ by means of the map $\widetilde 
T_{bc,n}$. By definition, we have 
the energy-momentum tensor map
$$T_{bc,n}:\Theta'(n)\to \A_{bc}(n).$$

\smallskip

It can be shown by 
a direct computation (as in \propref{stressdilaton}), that $\Theta'(n)\simeq 
\Theta'_{c_n}$. In the next subsection
we shall show a way to establish such an isomorphism without making any 
computations: we shall deduce it from the
so-called Bose-Fermi correspondence.

\begin{cor}
The line bundle on $\Mf$ whose fiber at $X\in\Mf$ is $\det(R\Gamma(X,\Omega^n))$ 
carries a canonical structure 
of a $\Theta'_{\Mf_{-c_n}}$-module.
\end{cor}

\begin{proof}

This assertion is simply a combination of \propref{confbc} and \propref{KZ}.

\end{proof}

\begin{rem} \label{canclass}

Take $n=-1$ (we have $c_{-1}=-26$) and consider the union of the 
connected components of $\Mf$ that correspond to 
curves of genus $g\geq 2$. In this case the line bundle $X\to 
\det(R\Gamma(X,\Omega^{-1}))$ is identified with 
$\Omega^{top}_\Mf$. Therefore, the ring of twisted differential operators on $\Mf$ 
corresponding to
the Picard Lie algebroid $\Theta'_{\Mf_{26}}$ 
is identified in this case with $D_\Mf^{op}$ in such a way that
$(\conf(\Xf,\A_{bc}(-1)))^*$ becomes the canonical right D-module over $\Mf$.

\end{rem}

\sssec{}
 Consider the local $\O$-module 
$p_2{}_*(\Omega^n\boxtimes\Omega^{1-n}(\Delta)/\Omega^n
\boxtimes\Omega^{1-n}(-\Delta)$.
We have a short exact sequence:
$$0\to\Omega\to 
p_2{}_*(\Omega^n\boxtimes\Omega^{1-n}(\Delta)/
\Omega^n\boxtimes\Omega^{1-n}(-\Delta)\to\O\to 0$$
and it is easy to show that we have in fact a canonical isomorphism 
$$p_2{}_*(\Omega^n\boxtimes\Omega^{1-n}
(\Delta)/\Omega^n\boxtimes\Omega^{1-n}(-\Delta)\simeq F_{2n-1},$$
(cf. \secref{dilaton}).

\smallskip

We define a map $\widetilde J:F_{2n-1}\to \A_{bc}(n)$ by setting
$$\widetilde J(\omega^n\boxtimes\omega^{1-n}\cdot f(x,y))=
(h\boxtimes\on{id})(\{\omega^n\boxtimes\omega^{1-n}\cdot f(x,y)\}).$$ 

\smallskip

Clearly, $\widetilde J$ is a map of local $\O$-modules. It is commonly referred to 
as the ``fermion number current''
due to the following property:

\begin{prop} \label{fermionnumber}
Let ${\mathfrak f}_1$ be a section of $F_{2n-1}$ that maps to $1$ 
under the projection $F_{2n-1}\to\O$
and set ${\mathfrak j}_1:=\widetilde J({\mathfrak f}_1)$. We have
$$
(h\boxtimes\on{id})(\{{\mathfrak j}_1,l\})=i\cdot l 
\text{ for } l\in\A_{bc}(n)^i.
$$
\end{prop}

\begin{proof}

It enough to check the assertion on the generators, i.e. for 
$l\in\Omega^n$ and for $l\in\Omega^{1-n}$, in which case it
follows immediately from the Jacobi identity.

\end{proof}

The following result is a part of the twisted Bose-Fermi correspondence:

\begin{thm} \label{bosefermi}
{\em (a)} The map $\widetilde J$ gives rise to a map of chiral algebras
$$J:\A(2n-1,-1)\to \A_{bc}(n).$$
The map $J$ is compatible with the local $\O$-module structure.

\noindent {\em (b)} We have a canonical isomorphism
$\Theta'(n)\simeq\Theta'_{c_n}$ and the composition
$$\Theta'_{c_n}\overset{T_{2n-1,-1}}\longrightarrow 
\A(2n-1,-1)\overset{J}\longrightarrow \A_{bc}(n)$$
coincides under this identification with the map $T_{bc,n}$.\footnote{Of course, 
$c_n:=-2(6n^2-6n+1)=1-3(2n-1)^2$.}
\end{thm}

\begin{proof}

To prove point (a) of the theorem it is enough to show that if ${\mathfrak f}_g\in 
F_{2n-2}$ is a section that maps 
to $g\in\O$ under the projection $F_{2n-1}\to \O$, then
$$(h\boxtimes\on{id})(\{J({\mathfrak f}_g),{\mathfrak j}_1\})=-dg.$$
Set ${\mathfrak j}_g=J({\mathfrak f}_g)$. As in \propref{fermionnumber} we have:
$$(h\boxtimes\on{id})(\{{\mathfrak j}_g,l\})=g\cdot l \text{ for } 
l\in\Omega^{1-n} \text{ and }
(h\boxtimes\on{id})(\{{\mathfrak j}_g,l\})=-g\cdot l \text{ for } 
l\in\Omega^{n}.$$
The assertion of point (a) now follows from the Jacobi identity.

\medskip

In order to prove point (b) observe that it is enough to show that the map 
$$J\circ\widetilde T_{2n-1,-1}:\Theta\to \A(2n-1,-1)/\Omega\to \A_{bc}(n)/\Omega$$
coincides with the map $\widetilde T_{bc,n}$. This will be done in two steps:

\medskip

\noindent {\bf Step 1}: We first show that $J\circ\widetilde T_{2n-1,-1}$ is 
proportional to $\widetilde T_{bc,n}$.

\smallskip

Indeed, $\widetilde T_{bc,n}$ is obviously a map $\Theta\to \A_{bc}(n)_2\cap 
\A_{bc}(n)^0$. We claim that
the same is true for $J\circ\widetilde T_{2n-1,-1}$. A priori, the latter is a map
$\Theta\to \A_{bc}(n)_4\cap \A_{bc}(n)^0$.

\smallskip

However, it is easy to see that its symbol, i.e. the composition
$$\Theta\to \A_{bc}(n)_4\cap \A_{bc}(n)^0/\A_{bc}(n)_3\cap \A_{bc}(n)^0$$
vanishes. Therefore, the image of $J\circ\widetilde T_{2n-1,-1}$ lies in 
$\Theta\to \A_{bc}(n)_3\cap \A_{bc}(n)^0$.
However, by parity considerations, $\A_{bc}(n)_3\cap \A_{bc}(n)^0=\A_{bc}(n)_2\cap 
\A_{bc}(n)^0$.

\smallskip

Observe now that as an $\O$-module, $\A_{bc}(n)_2\cap \A_{bc}(n)^0$ 
is identified with
$$(\Omega^n\otimes D)\overset{!}\otimes(\Omega^{1-n}\otimes D),$$ according to 
\lemref{filtration}(c). 

\smallskip

The assertion of Step 1 now follows from the next lemma:

\begin{lem}
The space $\on{Hom}(\Theta, (\Omega^n\otimes 
D)\overset{!}\otimes(\Omega^{1-n}\otimes D))$ in the category
of local $\O$-modules is $1$-dimensional.
\end{lem}

\noindent {\bf Step 2}:

\smallskip

We know by now that $J\circ\widetilde T_{2n-1,-1}=\epsilon\cdot \widetilde 
T_{bc,n}$ for some $\epsilon\in\CC$
and it remains to show that $\epsilon=1$.

\smallskip

For $\xi\in\Theta$ let $\xi'$ be any section of $\A_{bc}(n)$ that projects to 
$T_{bc,n}(\xi)$.
Since $T_{bc,n}$ is an energy-momentum tensor for $\A_{bc}(n)$, on the one hand we 
have
$$(h\boxtimes\on{id})(\{\epsilon\cdot\xi'\boxtimes J(l)\})=\epsilon(-J(l)\cdot 
\xi-J(\on{Lie}_\xi(l)))\,\, 
\forall l\in \A_{bc}(n) .$$

\smallskip

On the other hand, since $J$ is a homomorphism of chiral algebras and 
$T_{2n-1,-1}$ is an 
energy-momentum tensor for $\A(2n-1,-1)$, we have
$$(h\boxtimes\on{id})(\{\epsilon\cdot\xi'\boxtimes J(l)\})=-J(l)\cdot 
\xi-J(\on{Lie}_\xi(l))\,\, \forall l\in \A(2n-1,-1).$$

The two equations are compatible only when $\epsilon=1$.

\end{proof}

\newpage

\centerline{\bf Chapter IV. BRST and String Amplitudes}
\ssec{The BRST complex}
\sssec{}  \label{extension}
 Let $\B$ be a Lie-* algebra over a curve $X$, which is locally free and finitely 
generated as a D-module.
We shall now attach to $\B$ a canonical central extension $$0\to \Omega\to \B''\to 
\B\to 0$$
of Lie-* algebras.

\smallskip

The construction will be based on the following general assertion:

\begin{lem} \label{hom}
Let $\M_1$, $\M_2$ be two D-modules on $X$ with $\M_1$ being locally free and 
finitely generated. Then:

\smallskip

\noindent{\em (a)} For a third D-module $\M$ on $X$, there is a canonical 
isomorphism:
$$\on{Hom}_{D(X^2)}(\M\boxtimes \M_1,\Delta_{!}(\M_2))\simeq 
\on{Hom}_D(\M,\M_1^*\overset{!}\otimes\M_2).$$

\smallskip

\noindent{\em (b)}
The canonical map (from point (a)) 
$\M_1^*\overset{!}\otimes\M_2\boxtimes\M_1\to\Delta_{!}(\M_2)$ induces
an isomorphism
$h(\M_1^*\overset{!}\otimes\M_2)\to \on{Hom}_D(\M_1,\M_2)$.
\end{lem}

\begin{cor}
For a Lie-* algebra $\B$ we have canonical maps 
$$\on{co-ad}:\B\boxtimes\B^*\to\Delta_{!}(\B^*) \text{ and } \on{co-br}:\B^*\to 
\B^*\overset{!}\otimes\B^*.$$
\end{cor}

We shall call the two maps of the above corollary ``the co-adjoint action'' and 
``the co-bracket'', respectively.

\medskip

Let us for a moment view $\B$ as a plain D-module on $X$ (i.e. forget the Lie-* 
algebra structure)
and consider the chiral algebra $\A_{bc}(\B)$ of \secref{bc}. Let $i_\B$ (resp., 
$i^*_\B$) denote the canonical embedding
of $\B$ (resp., of $\B^*$) into $\A_{bc}(\B)$.

\smallskip

It follows from \lemref{filtration}(c) that the intersection 
$\A_{bc}(\B)_2\cap\A_{bc}(\B)^0$ is an extension of D-modules
$$0\to\Omega\to \A_{bc}(\B)_2\cap\A_{bc}(\B)^0\to \B\overset{!}\otimes\B^*\to 0.$$

\smallskip

According to \lemref{hom} above, the Lie-* bracket on $\B$ gives rise to a map
$$S_\B:\B\to \B^*\overset{!}\otimes\B\hookrightarrow \A_{bc}(\B)/\Omega.$$

\begin{prop}  \label{endom}
For two local sections $b_1$ and $b_2$ of $\B$ let $S_\B(b_1)'$ and $S_\B(b_2)'$ 
be local sections of $\A_{bc}(\B)$
that project to $S_\B(b_1)$ and $S_\B(b_2)$, respectively. Then the section 
$\{S_\B(b_1)'\boxtimes S_\B(b_2)'\}$ of $\Delta_{!}(\A_{bc}(\B))$ projects to the 
section $S_\B(\{b_1,b_2\})$ of
$\Delta_{!}(\A_{bc}(\B)/\Omega)$. 
\end{prop}

\begin{proof}

It is enough to prove that
$$(h\boxtimes h)(\{S_\B(b_1)'\boxtimes S_\B(b_2)'\})=h(S_\B(\{b_1,b_2\}))\in 
h(\A_{bc}(\B)/\Omega)$$ for
any $b_1$ and $b_2$ as above.

\smallskip

For a section $l$ of $h(\B^*\overset{!}\otimes\B)$ and for some 
lift of $l$ to a section $l'$ of\break
$h(\A_{bc}(\B)_2\cap\A_{bc}(\B)^0)$ consider the map $\B\to \A_{bc}(\B)$ 
given by $$b\to (h\boxtimes\on{id})(\{l'\boxtimes i_\B(b)\}).$$

\smallskip

It follows from the definition of the chiral bracket on $\A_{bc}(\B)$ that this 
map takes values in $\B\subset \A_{bc}(\B)$
and coincides with the canonical map $h(\B^*\overset{!}\otimes\B)\otimes\B\to\B$ 
of point (b) of \lemref{hom}.

\smallskip

The assertion of the proposition now follows from the Jacobi identity combined 
with \lemref{hom}(b), since 
both $(h\boxtimes h)(\{S_\B(b_1)'\boxtimes S_\B(b_2)'\})$ and 
$h(S_\B(\{b_1,b_2\}))$ belong to 
$h(\B^*\overset{!}\otimes\B)\subset h(\A_{bc}(\B)/\Omega)$.

\end{proof}

We now define $\B''$ to be the induced extension:
$$
\CD
0 @>>> \Omega @>>> \B'' @>>> \B @>>> 0 \\
@VVV @V{\on{id}}VV @V{S''_\B}VV @V{S_\B}VV @VVV \\
0 @>>> \Omega @>>> \A_{bc}(\B)_2\cap\A_{bc}(\B)^0 @>>> \B^*\overset{!}\otimes\B 
@>>> 0.
\endCD
$$

\propref{endom} above implies that $\B''$ acquires a canonical Lie-* algebra 
structure which is a central
extension of that of $\B$.

\sssec{}  \label{differential}
 Let $\B'$ denote the extension 
$0\to\Omega\to \B'\to\B\to 0$ obtained from $\B''$ 
by acting by $-1$ on
$\Omega$. Let $\U(\B)'$ denote the quotient of $\U(\B')$ by 
the standard relation $\unit-\widetilde\unit=0$
and consider now the tensor product chiral algebra 
$\A_{BRST}(\B):=\U(\B)'\overset{!}\otimes\A_{bc}(\B)$.
\footnote{It follows
from the definitions that, for two chiral algebras $\A_1$ and $\A_2$, 
the tensor product $\A_1\overset{!}\otimes A_2$ 
acquires a natural chiral algebra structure.} 
The chiral algebra $\A_{BRST}(\B)$ inherits from $\A_{bc}(\B)$
a $\ZZ$-- and a $\ZZ_2$--grading.

\smallskip

From the construction of $\B'$ it follows that we have a canonical map of Lie-* 
algebras $S_{\B,BRST}:\B\to \A_{BRST}(\B)$.

\medskip

Our next goal will be to construct a canonical element $\delta_\B\in 
h(\A_{BRST}(\B))$ with 
the following properties:

\begin{itemize}

\item

$\delta_\B$ has degree $1$ with respect to the $\ZZ$-grading on $\A_{BRST}(\B)$
and $$(h\boxtimes h)(\{\delta_\B,\delta_\B\})=0.$$

\item

$(h\boxtimes\on{id})(\{\delta_\B,S_{\B,BRST}(b)\})=0$ for any $b\in\B$. 

\item

$(h\boxtimes\on{id})(\{\delta_\B,i_\B(b)\})=-S_{\B,BRST}(b)$ for any $b\in\B$.

\end{itemize}

In particular, for any chiral $\B'$-module $\L$, the map 
$(h\boxtimes\on{id})(\{\delta_\B,\cdot\})$
would define a differential on the tensor product 
$\L\overset{!}\otimes\A_{bc}(\B)$, which is homogeneous of degree $1$ 
and whose square equals $0$.

\bigskip

We have two maps $G_1,G_2:j_*j^*(\B\boxtimes\B^*)\to\Delta_{!}(\A_{BRST}(\B))$:

\smallskip

\noindent
The map $G_1$ is defined as minus the composition
$$
j_*j^*(\B\boxtimes\B^*){\overset{S_{\B,BRST}\boxtimes 
i_\B^*}\longrightarrow}
j_*j^*(\A_{BRST}(\B)\boxtimes 
\A_{BRST}(\B)){\overset{\{,\}}\longrightarrow}
\Delta_{!}(\A_{BRST}(\B)).
$$

\smallskip

\noindent
To define the map $G_2$ note that the chiral bracket on $\A_{bc}(\B)$ gives rise 
to a map 
$\B^*\overset{!}\otimes\B^*\to \A_{bc}(\B)$ and let $i_{\B,2}^*$ denote the 
composition
$\B^*\overset{\on{co-br}}\longrightarrow \B^*\overset{!}\otimes\B^*\to 
\A_{bc}(\B)$. We set the map $G_2$ to be the composition:
$$j_*j^*(\B\boxtimes\B^*){\overset{i_\B\boxtimes 
i_{\B,2}^*}\longrightarrow}j_*j^*(\A_{bc}(\B)\boxtimes \A_{bc}(\B))
\overset{\{,\}}\longrightarrow\Delta_{!}(\A_{bc}(\B))\to 
\Delta_{!}(\A_{BRST}(\B)).$$

\begin{lem}
The restriction of $G_1$ (resp., of $G_2$) to $\B\boxtimes\B^*\subset 
j_*j^*(\B\boxtimes\B^*)$ coincides with the map
$$\B\boxtimes\B^*\overset{\on{co-ad}}\longrightarrow\Delta_{!}(\B^*)
\overset{i_\B^*}\longrightarrow \Delta_{!}(\A_{BRST}(\B)).$$
\end{lem}

In particular, the difference $G_1-G_2$ is a well-defined map
$$\B\overset{!}\otimes\B^*\to \A_{BRST}(\B).$$

\smallskip

We set $\delta_\B$ to be the image under $G_1-G_2$ of the canonical element 
$\on{id}_\B\in \on{End}_D(\B)\simeq
h(\B\overset{!}\otimes\B^*)$

\begin{prop}
The element $\delta_{BRST}$ satisfies conditions (1)--(3) above.
\end{prop}

\begin{proof}

Property (3) is a straightforward calculation. 

\smallskip

To prove property (2), it is enough to show that
$$
(h\boxtimes h)(\{\delta_{BRST},
S_{\B,BRST}(b)\})=0\in h(\A_{BRST}(\B)) 
$$
for any $b\in\B$.

\smallskip

Consider the action of $h(\B)$ on the D-module $\B\boxtimes\B^*$ given by the 
formula:
$$b'\cdot (b\boxtimes b^*)\to (h\boxtimes\on{id})(\{b',b\})\boxtimes b^*+ 
b\boxtimes (h\boxtimes \on{id})(\on{co-ad}(b'\boxtimes b^*)).$$
This action extends to an action of $h(\B)$ on $j_*j^*(\B\boxtimes\B^*)$.

\smallskip

Consider also the action of $\B$ on $\A_{BRST}(\B)$ given by 
$$
b\cdot l\to (h\boxtimes\on{id})(\{S_{\B,BRST}(b),l\}).
$$

\begin{lem}
The maps $G_1$ and $G_2$ commute with the above $h(\B)$--actions.
\end{lem}

Property (2) now follows from the fact that the section $\on{id}_\B\in 
h(\B\overset{!}\otimes\B^*)$
is $h(\B)$--invariant.

\smallskip

To establish property (1) note that the element $\delta_{BRST}^2$ 
satisfies
$\{\delta_{BRST}^2,i_\B(b)\}=0$ for all 
$b\in\B$. The assertion now follows from the next lemma:

\begin{lem}
Let $\M$ be as in \secref{bc} and let $l\in h(\A_{bc}(\M))$ be an element of 
positive degree. Assume that
$(h\boxtimes\on{id})(\{l,i_\M(m)\})=0\,\,\,\forall m\in h(\M)$. Then $l=0$.
\end{lem}

\end{proof}

\ssec{String amplitudes}
\sssec{}
 Let us now apply the discussion of the previous subsection to the case 
$\B=\T=\Theta\otimes D_X$.

\smallskip

In this case the chiral algebra $\A_{bc}(\B)$ is identified with 
the $bc$-system chiral algebra $\A_{bc}(n)$ for
$n=-1$.

\begin{lem}
The canonical central extension $\B''$ of \secref{extension} is
identified for $\B=\T$ with $\T'_{-26}$. Moreover,
under this identification, the map $S''_\B$ goes over to the energy-momentum 
tensor $T_{bc,-1}$.
\end{lem}

This assertion follows immediately from the definitions (using the equality 
$c_{-1}:=
-2(6\cdot(-1)^2-6\cdot(-1)+1=-26$).

\medskip

Now let $(\A,T)$ be a CFT data of central charge $26$. Consider the chiral algebra 
$\A_{comp}:=\A\overset{!}\otimes\A_{bc}(-1)$. The energy momentum tensor for $\A$ 
gives rise to
a homomorphism of chiral algebras:
$$\varphi:\A_{BRST}(\T)\to \A_{comp},$$
in particular, by composing it with the map $S_{\T,BRST}$ we obtain a map 
$T_{comp}:\Theta\to\T\to \A_{comp}$.

\smallskip

\begin{prop}
The pair $\A_{comp},T_{comp}$ defines a CFT data of central charge $0$.
\end{prop}

\begin{proof}

Since both $\A$ and $\A_{bc}(-1)$ carry a structure of local $\O$-module over 
$\Xf$, so does $\A_{comp}$.

\smallskip

The fact that the pair $(\A_{comp},T_{comp})$ satisfies conditions (a) and (b) of 
\secref{defcft} follows from
the corresponding facts for $\A$ and for $\A_{bc}(-1)$.

\end{proof}

\sssec{}
 Let $\delta_{comp}\in h(\A_{comp})$ denote the image of $\delta_{BRST}\in 
h(\A_{BRST}(\Theta))$ under the homomorphism
$\varphi$. 

\smallskip

Note that the $\O_\Mf$-module $X\to\conf(X,\A_{comp})$ carries a canonical 
structure of left $D_\Mf$-module, 
according to \propref{nocentral}. The next proposition describes the connection 
between the differential 
$\delta_{comp}$ and the De Rham differential on $\conf(\Xf,\A_{comp})$:

\begin{prop} \label{brstchar}
Let $f\cdot a_1\boxtimes\ldots\boxtimes a_k$ be a $\ZZ_2$-homogeneous section of 
$j_*j^*({}^l\A_{comp}^{\boxtimes k})$ on 
$\Xf^k\setminus \Delta$. Then\footnote{Here $^li_\T$ denotes the map 
$\Omega^{-2}\to{}^l\A_{comp}$ obtained from the canonical map
$i_\T:\Theta\to \A_{comp}$
and its insertion into the correlation functions has the same meaning as in 
\propref{anotherchar}.}
\begin{align*}
&\nabla(\langle f\cdot 
a_1,\dotsc,a_k\rangle)=(\langle\ldots\rangle\boxtimes\on{id})
(\nabla_{loc}(f\cdot a_1,\dotsc,a_k))+ \\
&+\underset{i=1}{\overset{k}\Sigma}\,
(-1)^{\on{deg}(a_i)+\ldots+\on{deg}(a_k)}\langle f\cdot 
a_1,\dotsc,\delta_{BRST}(a_i),\dotsc,a_k,{}^li_\T\rangle.
\end{align*}
\end{prop}

\begin{proof}

Using \propref{anotherchar}, the assertion of the proposition reduces to the fact 
that
$$\underset{i=1}{\overset{k}\Sigma}\,
(-1)^{\on{deg}(a_i)+\ldots+\on{deg}(a_k)}\langle f\cdot 
a_1,\dotsc,\delta_{BRST}(a_i),\dotsc,a_k,{}^li_\T\rangle=
\langle f\cdot a_1,\dotsc,a_k,{}^lT\rangle.$$

\smallskip

This, however, follows from the equality
\begin{align*}
&\underset{i=1}{\overset{k}\Sigma}\,
(-1)^{\on{deg}(a_1)+\ldots+\on{deg}(a_{i-1})}\langle f\cdot 
a_1,\dotsc,\delta_{BRST}(a_i),\dotsc,a_k,{}^li_\T\rangle+ \\
&+(-1)^{\on{deg}(a_1)+\ldots+\on{deg}(a_k)}\langle f\cdot 
a_1,\dotsc,a_k,\delta_{BRST}({}^li_\T)\rangle=0
\end{align*}
and property (3) of the element $\delta_{BRST}$ (cf. \secref{differential}).

\end{proof}

\medskip

Since $c=26$, the space of conformal blocks $\conf(\Xf,\A)$ carries a canonical 
right D-module structure over
$\Mf$ (cf. \remref{canclass}). Moreover, over the connected components of $\Mf$ 
corresponding to curves of genus $g\geq 2$,
we have an isomorphism $$\conf(X,\A)\simeq 
\conf(X,\A_{comp})\otimes\Omega^{top}_\Mf,$$
in view of the following lemma:

\begin{lem}
For two chiral algebras $\A_1$ and $\A_2$ and a complete curve $X$ we have a 
canonical isomorphism
$$\conf(X,\A_1\overset{!}\otimes\A_2)\simeq \conf(X,\A_1)\otimes\conf(X,\A_2).$$
\end{lem}

\smallskip

Therefore, \propref{brstchar} expresses the 
right D-module structure on $\conf(\Xf,\A)$
in terms of the fields of the composite 
theory (i.e. of $\A_{comp})$. 

\newpage

\centerline{\bf Chapter V. Further Constructions}

\ssec{Chiral algebras via the Ran space} \label{ran}
 The definition of chiral algebras that we used above was a rather 
straightforward 
axiomatization of the OPE operation on
quantum fields. We shall now discuss a different approach to chiral algebras,  
which is much less obvious from the point of view of QFT. 

\smallskip

In what follows we shall use the following conventions:

\smallskip

For a curve $X$ and a surjection of finite sets $\phi:\I\twoheadrightarrow \J$, 
$\Delta_\phi$ will denote
the corresponding diagonal embedding $X^\J\to X^\I$.

\smallskip

If $\I=\I_1\cup\ldots\cup \I_k$ is a decomposition of $I$ into a disjoint union of 
finite subsets
we shall denote by $j_{\I_1,\dotsc,\I_k}$ the embedding of the open subset 
$X^{\I_i\neq \I_j}$ of $X^\I$ that consists of 
points $$x_1^1,\dotsc,x_{|\I_1|}^1,\dotsc,x_1^k,\dotsc,x_{|\I_k|}^k\in 
X^{\I_1}\times\ldots\times X^{\I_k}\simeq X^\I$$ 
with $x_{i_1}^{j_1}\neq x_{i_2}^{j_2}$ whenever $j_1\neq j_2$.

\smallskip

When $k=|\I|$, we shall replace the notation $j_{\I_1,\dotsc,\I_k}$ by $j_\I$ or 
simply by $j$.

\sssec{} \label{An}
 For a chiral algebra $\A$ consider the corresponding left D-module $^l\A$ and for 
a finite set $\I$ consider
the left D-module $^l\A^\I$ over $X^\I$.

\smallskip

Now let $\phi:\I\twoheadrightarrow \J$ be a surjection of finite sets with 
$|\J|=|\I|-1$. 
The chiral bracket on $\A$ yields a well-defined map
$$\{,\}_\phi:j_\I{}_*j_\I^*({}^l\A^\I)\to 
\Delta_{\phi}{}_{!}(j_\J{}_*j_\J^*({}^l\A^\J)).$$

\smallskip

The fact that the chiral bracket on $\A$ satisfies the Jacobi identity implies 
that we can form a complex
$C^{\bullet}$ of left D-modules on $X^n$, with $C^{-k}$ being
$$\underset{\phi:\{1,\dotsc,n\}\surj \I,|\I|=k}\oplus 
\Delta_{\phi}{}_{!}(j_\I{}_*j_\I^*({}^l\A^\I)).$$

\smallskip

Let $^l\A^{(n)}$ denote the $-n$-th cohomology of this complex, i.e. $^l\A^{(n)}$ 
is a sub-D-module of $^l\A^n$
obtained by intersecting the kernels of the maps $\{,\}_\phi$ for all possible 
pairs $(\phi:\{1,\dotsc,n\}\to \I, |\I|=n-1)$.
For example, $^l\A^{(2)}$ is the left D-module corresponding to $\A^{(2)}$ of 
\secref{confblocks}. 

\smallskip

Since $^l\A^{(n)}$ is equivariant with respect to the action of the symmetric 
group on $X^n$, we can form a D-module 
$^l\A^{(\I)}$ for any finite non-empty set $\I$.

\smallskip

\begin{prop}
$H^i(C^{\bullet})=0 \text{ for } i\neq -n$.
\end{prop}

\begin{proof}

To simplify the notation, we shall give a proof in the case $n=3$. In the general 
case, the proof is
completely analogous.

\smallskip

Our complex (for the corresponding right D-modules) looks as follows:
\begin{align*}
j_{\{1,2,3\}}{}_*j_{\{1,2,3\}}^*(\A^3)
&\overset{d_{-3}}\longrightarrow 
\Delta_{x_1=x_2}{}_{!}(j_*j^*(\A^2))\oplus
\Delta_{x_2=x_3}{}_{!}(j_*j^*(\A^2))\\
\oplus 
\Delta_{x_3=x_1}{}_{!}(j_*j^*(\A^2)) 
&\overset{d_{-2}}\longrightarrow\Delta_{x_1=x_2=x_3}{}_{!}(\A)\to 0
\end{align*}
and we have to prove that it is exact at the $-2$-nd term (the fact that it is 
exact at the $-1$-st term is obvious,
since $\A$ has a unit).

\smallskip

We claim, first of all, that any section of $C^{-2}$ is equivalent modulo the 
image of $d_{-3}$
to a section that belongs to $\Delta_{x_2=x_3}{}_{!}(j_*j^*(\A^2))$. This again 
follows from 
the fact that $\A$ has a unit.

\smallskip

Now let $l$ be a section 
$$l\in\Gamma(X\times X,j_*j^*(\A^2))\subset\Gamma(X\times X\times 
X,\Delta_{x_2=x_3}{}_{!}(j_*j^*(\A^2)))$$
with $d_{-2}(l)=0$; it remains to show that $l$ lies in the image of $d_{-3}$.
 
\smallskip

Choose (locally) a $1$-form
$\omega$ on $X$ and a function $f(x_1,x_2)\in \O\boxtimes\O(\Delta)$ such that 
$f(x_1,x_2)\cdot\omega=1\on{mod}
\O\boxtimes\O(-\Delta)$ (we have used here the identification $\Omega\simeq 
\O\boxtimes\O(-\Delta)/\O\boxtimes\O(-2\cdot\Delta)$).

\smallskip

It now follows from the Jacobi identity that $l=d_{-3}(l\boxtimes\omega\cdot 
f(x_2,x_3))$.

\end{proof}

\begin{cor} \label{generalfactorization}
For a surjection $\phi:\I\surj \J$, we have a canonical isomorphism
$$\Delta_\phi^{!}({}^l\A^{(\I)})\simeq {}^l\A^{(\J)}[|\I|-|\J|]$$ and 
for a decomposition $\I=\I_1\cup\ldots\cup \I_k$ we have a canonical isomorphism 
$$j_{\I_1,\dotsc,\I_k}^*({}^l\A^{(\I)})\simeq 
j_{\I_1,\dotsc,\I_k}^*({}^l\A^{(\I_1)}\boxtimes\ldots\boxtimes {}^l\A^{(\I_k)}).$$ 
Moreover, the above isomorphisms are compatible
in the obvious sense. 
\end{cor}

We shall call the isomorphisms of \corref{generalfactorization} ``factorization 
isomorphisms''.

\smallskip

\begin{rem}

Consider the topological space $Ran(X)$ that consists of all finite non-empty 
subsets of $X$. The above corollary implies that
the collection of $^l\A^{(\I)}$'s can be viewed as ``D-module'' over $Ran(X)$: its 
fiber over a finite subset of $X$ 
that corresponds to a map $\I\to X$ is, by definition, the stalk of $^l\A^{(\I)}$ 
at the corresponding point of $X^\I$.

\end{rem}

\sssec{} \label{reconstruction}
 Next we want to show that the chiral algebra structure can be completely 
recovered from the system of
the $^l\A^{(n)}$'s viewed as plain $\O_{X^n}$-modules.

\smallskip

Let $^l\A$ be an $\O$-module on $X$. Assume now that for every finite set $\I$ we 
are given an $\O_{X^\I}$-module 
$^l\A'{}^{(\I)}$ (such that $^l\A={}^l\A'{}^{(\{1\})}$) together
with a compatible system of factorization isomorphisms:
$$\Delta_\phi^*({}^l\A'{}^{(\I)})\simeq {}^l\A'{}^{(\J)}\text{ for a surjection 
}\phi:\I\surj \J\text{ and }$$
$$j_{\I_1,\dotsc,\I_k}^*({}^l\A^{(\I)})\simeq 
j_{\I_1,\dotsc,\I_k}^*({}^l\A'{}^{(\I_1)}\boxtimes\ldots\boxtimes 
{}^l\A'{}^{(\I_k)})\text { for a decomposition }\I=\I_1\cup\ldots\cup \I_k.$$

\smallskip

Assume also that there exists a map $\unit:\O\to {}^l\A$ with the following 
properties:

\smallskip

\begin{itemize}

\item
For $\I=\I_0\cup i$ the map 
$$\on{id}\boxtimes \unit:j_{\I_0,\{i\}}^*({}^l\A'{}^{(\I_0)}\boxtimes \O)\to
j_{\I_0,\{i\}}^*({}^l\A'{}^{(\I_0)}\boxtimes {}^l\A)\simeq
j_{\I_0,\{i\}}^*({}^l\A'{}^{(\I)})$$ extends to a map
$\on{id}\boxtimes\unit:{}^l\A'{}^{(\I_0)}\boxtimes\O\to {}^l\A'{}^{(\I)}$.

\item
Take $\I=\I_0\cup i_1\cup i_2$ and $\J=\I_0\cup i_1$ and let $\phi:\I\to \J$ be a 
map that contracts $\{i_1,i_2\}\to i_1$.
We need that the composition
$${}^l\A'{}^{(\J)}\simeq \Delta_\phi^*({}^l\A'{}^{(\J)}\boxtimes 
\O)\overset{\on{id}\boxtimes\unit}
\longrightarrow \Delta_\phi^*({}^l\A'{}^{(\I)})\simeq {}^l\A'{}^{(\J)}$$
is the identity map.

\end{itemize}

\begin{thm} \label{reconstr}
Let $({}^l\A,{}^l\A'{}^{(\I)},\unit)$ be data as above. 
Then there exists a (canonical) chiral algebra 
structure on $\A:={}^l\A\otimes\Omega^{-1}$ such that the corresponding D-modules 
$^l\A^{(\I)}$ are identified as
$\O_{X^\I}$-modules with the ${}^l\A'{}^{(\I)}$'s.
\end{thm}

\smallskip

We shall not give here a complete proof of this theorem. Instead, we will explain 
the key points:

\begin{proof} (sketch)

Let $\A'{}^{(n)}$ denote the $\O_{X^n}$-module corresponding to 
$\I=\{1,\ldots,n\}$. 

\smallskip

We shall first show how to endow $^l\A$ with a left D-module structure:

\smallskip

The embedding of $\O_{X^2}$-modules 
$\unit\boxtimes\on{id}:\O\boxtimes{}^l\A\hookrightarrow{}^l\A'{}^{(2)}$ induces an 
isomorphism
between their restrictions to the diagonal. This fact implies that
the restrictions of the above two $\O_{X^2}$-modules to the second infinitesimal 
neighborhood of 
$\Delta(X)\subset X\times X$ are isomorphic as well, i.e. that
$$\O\boxtimes {}^l\A/\O\boxtimes{}^l\A(-2\cdot\Delta)\simeq 
{}^l\A'{}^{(2)}/{}^l\A'{}^{(2)}(-2\cdot\Delta).$$

\smallskip

Symmetrically, we have an isomorphism
$$^l\A\boxtimes\O/{}^l\A\boxtimes\O(-2\cdot\Delta)\simeq 
{}^l\A'{}^{(2)}/{}^l\A'{}^{(2)}(-2\cdot\Delta)$$
and by transitivity, we obtain an isomorphism of $\O_{X^2}$-modules:
$$\varphi:\O\boxtimes {}^l\A/\O\boxtimes{}^l\A(-2\cdot\Delta)\to
^l\A\boxtimes\O/{}^l\A\boxtimes\O(-2\cdot\Delta).$$

\smallskip

This enables us to define a connection on $^l\A$: to a section $l\in{}^l\A$ we 
associate a section
$\nabla(l)\in {}^l\A\otimes\Omega$ as follows:

\smallskip

Note that $l\boxtimes 1-\varphi(1\boxtimes l)$ is a section of 
$^l\A\boxtimes\O/{}^l\A\boxtimes\O(-2\cdot\Delta)$
which vanishes modulo $^l\A\boxtimes\O(-\Delta)$. We set $\nabla(l)$ to be the 
image
of $l\boxtimes 1-\varphi(1\boxtimes l)$ under the identification
$$^l\A\boxtimes\O(-\Delta)/{}^l\A\boxtimes\O(-2\cdot\Delta)\simeq 
{}^l\A\otimes\Omega.$$

\smallskip

We leave it to the reader to check the correctness of this definition (this can be 
done using $^l\A'{}^{(3)}$).
Moreover, the above construction can be generalized to produce left D-module 
structures on all the $^l\A'{}^{(\I)}$'s in
a way compatible with the factorization isomorphisms.

\medskip

Let us now show how $\A$ acquires a chiral algebra structure:

\smallskip

The short exact sequence of D-modules
$$0\to {}^l\A'{}^{(2)}\to j_*j^*({}^l\A'{}^{(2)})\to 
\Delta_{!}\Delta^{!}({}^l\A'{}^{(2)})[1]\to 0$$
gives rise under our identifications to a map $j_*j^*({}^l\A\boxtimes{}^l\A)\to 
\Delta_{!}({}^l\A)$.

\smallskip

Using the D-module $^l\A'{}^{(3)}$, it can be checked that this operation 
satisfies the Jacobi identity when
we pass to the corresponding right D-modules.

\medskip

Finally, we have to identify $^l\A'{}^{(\I)}$ with $^l\A^{(\I)}$.

\smallskip

It follows from the construction that there is a natural embedding 
$^l\A'{}^{(\I)}\hookrightarrow{} ^l\A^{(\I)}$.
Moreover, for every  $\phi:\I\surj \J$ this map induces an isomorphism
$$j_\J^*\Delta^{!}_\phi(^l\A'{}^{(\I)})\simeq j_\J^*(^l\A'{}^{(\J)})\simeq 
j_\J^*(^l\A^{(\J)})
\simeq j_\J^*\Delta^{!}_\phi(^l\A^{(\I)}).$$
This implies that the embedding $^l\A'{}^{(\I)}\hookrightarrow{} ^l\A^{(\I)}$ is 
an isomorphism, too.

\end{proof}

\begin{rem}

One may have observed that the construction of the chiral universal enveloping 
algebra discussed in \secref{unenv}
was also in the spirit of \thmref{reconstr}. In fact, it can be completely 
reformulated in terms of the system 
of $\O_{X^n}$-modules $\U(\B)^{(n)}$, which we leave as an exercise to 
the reader. 

\end{rem} 

\sssec{}
 Let $\A$ be a chiral algebra over a complete curve $X$. It follows from 
\secref{confsummary} that the
map $$\langle\ldots\rangle:j_*j^*({}^l\A\boxtimes\ldots\boxtimes{}^l\A)\to
j_*j^*(\O\boxtimes\ldots\boxtimes\O)\otimes\conf(X,\A)$$
restricts to a D-module map
$$^l\A^{(n)}\to \O_{X^n}\otimes\conf(X,\A).$$ Hence, we obtain a map of
vector spaces $DR^n(X^n,{}^l\A^{(n)})\to \conf(X,\A)$.

\begin{prop} \label{confran}
For $n\geq 2$, the map $DR^n (X^n,{}^l\A^{(n)})\to \conf(X,\A)$ is an 
isomorphism.
\end{prop} 

\begin{proof}

We can compute $DR^n(X^n,{}^l\A^{(n)})$ using the resolution $C^{\bullet}$ of 
\secref{An}. We have:
$$DR^n(X^n,{}^l\A^{(n)})\simeq\on{coker}(DR^1(X^n,C^{-2})\to DR^1(X^n,C^{-1})).$$

\smallskip

However, $DR^1(X^n,C^{-1})$ is just $DR^1(X,\A)$ and $DR^1(X^n,C^{-2})$ is a 
direct sum of several copies of
$DR^1(X\times X,j_*j^*(\A\boxtimes \A))$. The assertion of the proposition follows 
from the fact that
$$\conf(X,\A)\simeq DR^1(X,\A^{(2)})\simeq \on{coker}(DR^1(X^2,j_*j^*(\A\boxtimes 
\A))\to DR^1(X,\A)).$$

\end{proof}

\ssec{Geometry of the affine Grassmannian} 

 Let $G$ be an algebraic group and let $\gff$ be the corresponding 
Lie algebra. Consider the
Lie-* algebra $\B(\gff,0)$ of Ex.1 in \secref{firstexamples}. In the next two 
subsections 
we shall describe a geometric construction of the corresponding chiral 
universal enveloping
algebra $\A(\gff,0)$ via the so-called affine Grassmannian of the group $G$.
This description will be used later on for the construction of the free bosonic 
theory.

\sssec{}  \label{intrgrass}
 Consider the group $G(\hat\O)$ and the ind-group $G(\hat\K)$ that classify 
maps $$\Spec(\hat\O)\to G \text{ and }
\Spec(\hat\K)\to G,$$ respectively. Consider now the affine Grassmannian 
$\Gr_G:=G(\hat\K)/G(\hat\O)$
corresponding to the group $G$. The ind-scheme $\Gr_G$ may be highly 
non-reduced; for instance, when $G=H$ is a torus, 
the corresponding reduced scheme is identified with the discrete set of 
co-characters of $H$, while $\Gr_H$ is 
infinite-dimensional.

\smallskip

Note that the group $\Aut^+_0$ of automorphisms of $\hat\O$ acts naturally on 
$\Gr_G$.

\medskip

For a curve $X$ consider the group scheme $G_\O(X)$ over $X$ whose fiber at 
$x\in X$ is $G(\hat\O_x)$ and a group ind-scheme
$G_\K(X)$, whose fiber at $x\in X$ is $G(\hat\K_x)$. The quotient 
$\Gr_G(X):=G_\K(X)/G_\O(X)$ is a global version of
the affine Grassmannian $\Gr_G$ considered above. Let $r:\Gr_G(X)\to X$ 
denote the 
natural 
projection and let $\unit_{\Gr}$ denote the unit section $X\to\Gr_G(X)$.

\smallskip

\begin{lem} \label{functorgr}
The ind-scheme $\Gr_G(X)$ represents the functor whose value on a scheme $S$ is 
the data of a map $f:S\to X$, a 
$G$-torsor $\P$ on $S\times X$ and a trivialization
$$\alpha:\P|_{S\times X\setminus \Gamma_f}\simeq \P_0|_{S\times X\setminus 
\Gamma_f},$$
where $\Gamma_f$ is the graph of the map $f$ and $\P_0$ is the trivial 
$G$-torsor.
\end{lem}

This is a version of the Beauville-Laszlo theorem (\cite{BL}).

\smallskip

\begin{prop}  \label{connongr}
The ind-scheme $\Gr_G(X)$ possesses a natural connection along $X$. The section 
$\unit_{\Gr}:X\to \Gr_G(X)$
is preserved by this connection.
\end{prop}

\begin{proof}

Let $I$ be a local Artinian scheme and let $f$ be a map $I\times S\to X$, where 
$S$ is an arbitrary scheme. 
Let $\Spec(\CC)\simeq I_0\subset I$ be the reduced (point) scheme and let let 
$\tilf_0$ be a map
$I_0\times S\to\Gr_G(X)$ such that the map $r\circ \tilf_0:I_0\times S\to X$ 
coincides with the composition 
$I_0\times S\to I\times S\to X$.

\smallskip

To define a connection on $\Gr_G(X)$ along $X$, we must associate to a triple 
$(I,f,\tilf_0)$ as above
a map $\tilf:I\times S\to\Gr_G(X)$ such that $r\circ\tilf=f$ that extends 
$\tilf_0$.

\smallskip

This can be done as follows: 

\smallskip

\noindent The corresponding map $I\times S\to X$ has already
been given to us (this is our $f$) and we set the $G$-torsor
$\P^I$ over $I\times S\times X$ that would correspond to $\tilf$ to be pulled pack 
under $I\times S\times X\to S\times X$ from the 
$G$-torsor
$\P$ that corresponds to $\tilf_0$. However, $I\times S\times X\setminus 
\Gamma_f\simeq (S\times X\setminus \Gamma_{f_0})\times I$
and the data of trivialization for $\P^I$ comes from the corresponding data for 
$\P$.

\end{proof}

\smallskip

When $X$ is complete let $\Bun_G(X)$ denote the moduli stack of $G$-bundles over 
$X$. 
We have a natural projection $k:\Gr_G(X)\to\Bun_G(X)$.

\begin{lem}
The projection $k:\Gr_G(X)\to\Bun_G(X)$ is preserved by the connection on 
$\Gr_G(X)$ along $X$.
\end{lem}

\medskip

It follows from the definitions that the dependence of $\Gr_G(X)$ on the curve $X$ 
is local
(in the sense of \secref{local}):

\smallskip

We have the ind-stack $\Gr_G(\Xf)$ fibered over $\Xf$, that is, for a family $X^S$ 
of curves over a base $S$
we can form a scheme $\Gr_G(X^S)$ over $X^S$. Now let
$X^{S,I},x^{S,I},\phi^{S,I}$ be as in \secref{local}. We have a canonical 
isomorphism:
$$\Gr_G^n(X^{S^I})\setminus r^{-1}(x^{S\times I})\simeq
(\Gr_G(X^S)\setminus r^{-1}(x^S))\times I.$$

\smallskip

The ind-scheme $\Gr_G(\Xf)$ carries a connection along the fibers of the 
projection $\Xf\to\Mf$ and it is easy to
see that the above isomorphism is compatible with this connection.

\sssec{}
 For a finite set $\I$ we define the Beilinson-Drinfeld 
Grassmannian $\Gr^\I_G(X)$ to be the ind-scheme 
representing the following functor:

\smallskip

\noindent For a scheme $S$, $\on{Hom}(S,\Gr^\I_G(X))$ is the data consisting 
of a map $f^\I:S\to X^\I$, a $G$-torsor
$\P$ on $S\times X$ and a trivialization
$$\alpha:
\P|_{S\times X\setminus \{\Gamma_{f^\I_{i_1}},\dotsc,\Gamma_{f^\I_{i_n}}\}}
\simeq \P_0|_{S\times X\setminus 
\{\Gamma_{f^\I_{i_1}},\dotsc,\Gamma_{f^\I_{i_n}}\}},$$
where for $i_k\in \I$, $f^\I_{i_k}$ is the composition of $f^\I$ with the 
projection on the $i_k$-th factor 
$X^\I\to X$.

\smallskip

Let $r^\I$ denote the natural projection $\Gr^\I_G(X)\to X^\I$ and let 
$\unit^\I_{\Gr}:X^\I\to \Gr^\I_G(X)$
denote the corresponding unit section.

\smallskip 

When $\I=\{1,\dotsc,n\}$ we shall denote $\Gr^\I_G(X)$ by $\Gr^n_G(X)$ and the 
corresponding maps $r^\I$ and  
$\unit^\I_{\Gr}$ by $r^n$ and $\unit^n_{\Gr}$, respectively.

\medskip

A remarkable feature of the ind-scheme $\Gr^\I_H(X)$ is the following 
factorization property:

\begin{prop} \label{factorgr}
For a surjection of finite sets $\I\surj \J$, we have a natural isomorphism:
$$X^\J\underset{X^\I}\times \Gr_G(X)^\I\simeq \Gr_G(X)^\J$$ and for a 
decomposition $\I=\I_1\cup\ldots\cup \I_k$ we have a 
natural isomorphism
$$X^{\I_i\neq \I_j}\underset{X^\I}\times \Gr_G(X)^\I\simeq X^{\I_i\neq 
\I_j}\underset
{X^{\I_1}\times\ldots\times X^{\I_k}}\times \Gr_G(X)^{\I_1}\times\ldots\times 
\Gr_G(X)^{\I_k}.$$
\end{prop}

\begin{proof}

The first assertion follows immediately from the definitions. To prove the second 
one let us assume for simplicity
that $\I=\{1,2\}$, $\I_1=\{1\}$ and $\I_2=\{2\}$.

\smallskip

Given an object $(x_1,x_2,\P,\alpha)$ of $\Gr^2_H(X)$ we construct an object of 
$$(x_1,\P_1,\alpha_1)\times (x_2,\P_2,\alpha_2)\in\Gr_G(X)\times\Gr_G(X)$$ as 
follows:

\smallskip 

\noindent The $G$-torsor $\P_1$ (resp., $\P_2$) is set to be isomorphic to $\P_0$ 
over $X\setminus x_1$ 
(resp., over $X\setminus x_2$) and to be isomorphic to $\P$ over $X\setminus x_2$
(resp., over $X\setminus x_1$). The datum of $\alpha$ for $\P$ provides the gluing 
data for
$\P_1$ and $\P_2$. The data of $\alpha_1$ and $\alpha_2$ follow from the 
construction.

\smallskip

Vice versa, given an object 
$(x_1,\P_1,\alpha_1,x_2,\P_2,\alpha_2)\in\Gr_G(X)\times\Gr_G(X)$
we define $(x_1,x_2,\P,\alpha)\in\Gr^2_G(X))$ by setting $\P$ to equal $\P_1$ over 
$X\setminus x_2$ and
to equal $\P_2$ over $X\setminus x_1$. The isomorphisms $\alpha_1$ and $\alpha_2$ 
provide the gluing data
for $\P$ together with the trivialization $\alpha$ over $X\setminus\{x_1,x_2\}$

\end{proof}

\begin{rem}

It follows from \propref{factorgr} above that the fiber of $\Gr^2_G(X)$ over a 
point $(x_1,x_2)\in X\times X$
is identified with $\Gr_G\times\Gr_G$ whenever $x_1\neq x_2$ and with $\Gr_G$ when 
$x_1=x_2$. This is a purely
infinite-dimensional phenomenon (in the finite-dimensional situation the dimension 
of fibers cannot
drop under a specialization). Moreover, it is not difficult to prove that the map 
$r^\I$ is formally smooth.

\end{rem}

\smallskip

As in the case $|\I|=1$, it is easy to show that $\Gr^\I_G(X)$ possesses a natural 
connection along $X^\I$ and
that the isomorphisms constructed in \propref{factorgr} are compatible with these 
connections. 
When $X$ is complete we have a projection $k^\I:\Gr^\I_G(X)\to\Bun_G(X)$ which (as 
in the case of $|\I|=1$)
is preserved by the connection on $\Gr^\I_G(X)$ along $X^\I$. (When 
$\I=\{1,\dotsc,n\}$ we shall replace the notation $k^\I$
simply by $k^n$.)

\smallskip

Again, as in the case of $\I=\{1\}$, there exists an ind-stack $\Gr^\I_G(\Xf)$ 
fibered over $\Xf^I$ which has
a locality property in the sense that was specified at the end of 
\secref{intrgrass}. It follows from the definitions
that the local structure on the $\Gr^\I_G(X)$'s is compatible with the connections 
along the $X^\I$'s and with the
factorization isomorphisms of \propref{factorgr}.

\sssec{} \label{twgrassmannian}
 We shall now introduce a twisted version of the Beilinson-Drinfeld Grassmannian, 
which will be used for the construction
of the bosonic chiral algebra in the next section.

\smallskip

For an integer $n$ we introduce the iterated Beilinson-Drinfeld Grassmannian 
$\tGr^n_G(X)$ as follows:

\smallskip

$\tGr^n_G(X)$ is the ind-scheme representing the functor whose value on a scheme 
$S$ is the data
consisting of a map $f^n:S\to X^n$ and a collection of $n$
$G$-torsors $\P_i$ over $S\times X$ together with isomorphisms 
$$\alpha_i:\P_i|_{S\times X\setminus \Gamma_{f^n_i}}\simeq\P_{i-1}|_{S\times 
X\setminus \Gamma_{f^n_{i-1}}},$$
(here $f^n_i$ is the composition of the map $f^n$ with the projection on the 
$i$-th factor $X^n\to X$.)

\smallskip

Let $\widetilde r^n$ denote the projection $\tGr^n_G(X)\to X^n$. It is easy to see 
that $\tGr^n_G(X)$ is a fibration
over ${\tGr}^{n-1}_G(X)\times X$ with typical fiber $\Gr_G$. 

\smallskip

More precisely, consider the direct product ${\tGr}^{n-1}_G(X)\times X$ and let 
${\widetilde G^n_\O(X)}$ denote the pull-back
of the group scheme $G_\O(X)$ with respect to the projection 
$$
{\tGr}^{n-1}_G(X)\times X\to X.
$$
On the one hand, there is a canonical ${\widetilde G^n_\O(X)}$-torsor over 
${\tGr}^{n-1}_G(X)\times X$, whose fiber at a 
point $(x_1,\dotsc,x_{n-1},\P_1,\dotsc,\P_{n-1},x_n)$ is the restriction of 
$\P_{n-1}$ to the
formal disc around $x_n$. On the other hand, ${\widetilde G^n_\O(X)}$ acts on the 
${\tGr}^{n-1}_G(X)\times X$-scheme
${\tGr}^{n-1}_G(X)\times\Gr_G(X)$. It is easy to see that the 
${\tGr}^{n-1}_G(X)\times X$-scheme $\tGr^n_G(X)$ is a twist
of ${\tGr}^{n-1}_G(X)\times\Gr_G(X)$ with respect to the above-mentioned 
${\widetilde G^n_\O(X)}$-torsor.

\medskip

We have a canonical projection $w_n:\tGr^n_G(X)\to\Gr^n_G(X)$ 
that sends the data of $(x_1,\dotsc,x_n,\P_1,\dotsc,\P_n,
\alpha_1,\dotsc,\alpha_n)$ to $(x_1,\dotsc,x_n,\P_n,\alpha)$, 
where $\alpha$ is obtained by composing
$\alpha_1,\dotsc,\alpha_n$.
By repeating the proof of \propref{factorgr}, we see that the map $w_n$ is an 
isomorphism over $X^n\setminus\Delta$.

\smallskip

As in the case of $\Gr^n_G(X)$, there exists a scheme 
$\tGr^n_G(\Xf)$ over $\Xf^n$ which
is local in the sense of \secref{intrgrass}.

\ssec{Chiral algebra attached to the affine Grassmannian}  
\sssec{}  \label{affgrchalg}
 Consider the right D-module $\unit_{\Gr}{}_{!}(\Omega)$ on
$\Gr_G(X)$.\footnote{When $Z$ is a strict ind-scheme (i.e. $Z$ is a union of
closed finite dimensional subschemes) and $Z'\subset Z$ is a closed
finite-dimensional subscheme of $Z$ it makes perfect sense to talk about
right D-modules on $Z$ supported on $Z'$ (cf. \secref{Dmod}).} As an
$\O_{\Gr_G(X)}$-module, $\unit_{\Gr}{}_{!}(\Omega)$ is a union of
$\O$-modules on $\Gr_G(X)$ that are supported scheme-theoretically on an
increasing family of finite-dimensional infinitesimal neighborhoods of
$\unit_{\Gr}(X)$ in $\Gr_G(X)$.

\smallskip

Consider now the $\O$-module direct image 
$$\A_{\Gr_G}:=r_*(\unit_{\Gr}{}_{!}(\Omega)).$$ 
The connection on $\Gr_G(X)$ along $X$ defines a right D-module structure on 
$\A_{\Gr_G}$ and the section $\unit_{\Gr}$ defines 
an embedding $\Omega\hookrightarrow \A_{\Gr_G}$. Let $^l\A_{\Gr_G}$ be the 
corresponding left D-module.

\smallskip

More generally, for a finite set $\I$ let $\A'{}^{(\I)}_{\Gr_G}$ be the right 
D-module on $X^\I$ defined as
$r^\I_*(\unit^\I_{\Gr}{}_{!}(\Omega_{X^\I}))$ and let $^l\A'{}^{(\I)}_{\Gr_G}$ 
denote the corresponding left D-module.

\medskip

Let us now view the $^l\A'{}^{(\I)}_{\Gr_G}$'s as plain $\O$-modules on $X^\I$. 
The next assertion follows from 
\propref{factorgr}:

\begin{lem}
The system $\I\to{}^l\A'{}^{(\I)}_{\Gr_G}$ satisfies the conditions of 
\thmref{reconstr}.
\end{lem}

\smallskip

By applying  \thmref{reconstr} we obtain, therefore, a chiral algebra structure on 
$\A_{\Gr_G}$ with the 
$^l\A'{}^{(\I)}_{\Gr_G}$'s being identified with the corresponding 
$^l\A_{\Gr_G}^{(\I)}$'s
as $\O$-modules. 

However, since the factorization isomorphisms of \propref{factorgr} are compatible 
with the connections
along $X^\I$, it is easy to see that the isomorphisms $\A_{\Gr_G}^{(\I)}\simeq 
\A'{}^{(\I)}_{\Gr_G}$ are compatible
with the D-module structure.

\sssec{}
 Consider the fiber $A_{\Gr_G}{}_x$ of the chiral algebra $\A_{\Gr_G}$ at a point 
$x\in X$. The action of the group 
$G(\hat\K_x)$ on $\Gr_G=G(\hat\K_x)/G(\hat\O_x)$ defines a structure of 
$(\gff\otimes\hat\K_x,G(\hat\O_x))$-module
on $A_{\Gr_G}{}_x$. We shall denote by $$\Act_x:(\gff\otimes\hat\K_x)\otimes 
A_{\Gr_G}{}_x\to A_{\Gr_G}{}_x$$
the corresponding action map.

\smallskip

It is easy to see, moreover, that $A_{\Gr_G}{}_x$ is in fact canonically 
isomorphic to the vacuum
module $\on{Ind}^{\gff\otimes\hat\K_x}_{\gff\otimes\hat\O_x}(\CC)$.

\medskip

Our goal now is to prove the following theorem:

\begin{thm} \label{kmisom}
There exists a natural isomorphism of chiral algebras 
$\A_{\Gr_G}\simeq\A(\gff,0)$. Moreover, for every $x\in X$
the map $\Act_x$ goes over to the action of
$\gff\otimes\hat\K_x\simeq DR^0(\Spec(\hat\K_x),\B(\gff,0))$ on $A(\gff,0)_x$.
\end{thm}

\begin{rem}

A similar description of the chiral algebra $\A(\gff,Q)$ can be given when $Q\neq 
0$. This will be done explicitly
in \secref{descrbose} when $G=H$ is a torus.

\end{rem}

\medskip

We shall deduce the assertion of the theorem from the following result:

\begin{prop} \label{basickm}
There is a natural map 
$$
\Act:j_*j^*(\B(\gff,0)\boxtimes 
\A_{\Gr_G})\to\Delta_{!}(\A_{\Gr_G})
$$
such that

\smallskip

\noindent{\em(a)} $\B(\gff,0)$ acts on $\A_{\Gr_G}$ by derivations of the chiral 
algebra structure.

\smallskip

\noindent{\em(b)} The induced map 
$$(h\boxtimes\on{id})(\Act):DR^0(\Spec(\hat\K_x),\B(\gff,0))\otimes 
A_{\Gr_G}{}_x\to A_{\Gr_G}{}_x$$
coincides with the map $\Act_x$.
\end{prop}

\begin{proof} (Of \thmref{kmisom})

According to property (b) of the map $\Act$, the composition
$$j_*j^*(\B(\gff,0)\boxtimes \Omega)\overset{\on{id}\times\unit}\longrightarrow 
j_*j^*(\B(\gff,0)\boxtimes \A_{\Gr_G})\to
\Delta_{!}(\A_{\Gr_G})$$
factors as
$$j_*j^*(\B(\gff,0)\boxtimes \Omega)\twoheadrightarrow \B(\gff,0)\to \A_{\Gr_G},$$
i.e. we obtain a map of D-modules 
$$emb:\B(\gff,0)\to \A_{\Gr_G}.$$ On the level of fibers, this map corresponds to 
the
embedding
$\hat\K_x/\hat\O_x\otimes\gff\hookrightarrow 
\on{Ind}^{\gff\otimes\hat\K_x}_{\gff\otimes\hat\O_x}(\CC).$

\smallskip

By again applying point (b) of \propref{basickm}, we obtain that 
$emb:\B(\gff,0)\to \A_{\Gr_G}$
is a map of chiral $\B(\gff,0)$-modules. 

\smallskip

Therefore, to prove the theorem it remains to show that the two maps
$$j_*j^*(\B(\gff,0)\boxtimes \A_{\Gr_G})\to\Delta_{!}(\A_{\Gr_G}),$$ namely $\Act$ 
and
$(\{,\})\circ (emb\times \on{id})$, coincide. This is, however, an almost 
immediate consequence of point (a) 
of \propref{basickm}.

\medskip

Let $b\boxtimes a\cdot f(x,y)$ be a section of $j_*j^*(\B(\gff,0)\boxtimes 
\A_{\Gr_G})$. It is enough to show that 
$$(\on{id}\times h)(\Act(b\boxtimes a\cdot f(x,y)))=
(\on{id}\times h)(\{emb(b)\boxtimes a\cdot f(x,y)\})$$
for any $a,b$ and $f(x,y)$ as above.

\smallskip

Let $\omega$ be a non-vanishing $1$-form on $X$ and let $h(x,y)$ be a function on 
$X\times X\setminus \Delta$ with a simple
pole along the diagonal and with $\on{Res}_\Delta(\omega\otimes h(x,y))=1$. 
Consider the section
$$b\boxtimes \omega\boxtimes a\cdot f(x,z)\cdot h(x,y)\in 
j_*j^*(\B(\gff,0)\boxtimes \A_{\Gr_G}\boxtimes \A_{\Gr_G}).$$

\smallskip

The fact that $\B(\gff,0)$ acts on $\A_{\Gr_G}$ by derivations of the chiral 
algebra structure implies that
\begin{align*}
&\{\Act(b\boxtimes\omega\cdot h(x,y))\boxtimes a\cdot f(x,z)\}= \\
&\Act(b\boxtimes\{\omega\boxtimes a\cdot h(x,y)\cdot f(x,z)\})-
\sigma_{1,2}\{\omega,\Act(b\boxtimes a\cdot h(y,x)\cdot f(y,z))\}.
\end{align*}

\smallskip

Note, first of all, that the expression $\Act(b\boxtimes\omega\cdot h(x,y))$ is 
equal by definition
to $emb(b)$ and that the first term on the RHS vanishes.

\smallskip

By applying $(\on{id}\boxtimes h\boxtimes h)$ to both sides of the above formula 
we obtain the needed result.

\end{proof}

\sssec{}
 The proof of \propref{basickm} will be based on the following construction:

\smallskip

Assume that $X$ is affine and consider the formal group 
$\on{Maps}(X,G)$ that classifies maps $X\to G$. In addition, given a finite set 
$\I$
we shall consider a formal group-scheme $\tMaps{^\I}(X,G)$ 
over $X^\I$ that corresponds to the functor whose value 
on a pair $(S,f^\I:S\to X^\I)$ is the group of regular maps 
$$(S\times X\setminus \{\Gamma(f_{i_1}^\I),\dotsc,\Gamma(f_{i_n}^\I)\})\to G.$$

\smallskip

It follows from the definitions that $\tMaps{^\I}(X,G)$ carries a natural 
connection along $X^\I$. 
The corresponding sheaf of Lie algebras 
$\on{Lie}(\tMaps{^\I}(X,G))$ is a left 
D-module on $X^\I$ which is identified with $(p_\I{}\circ 
j_{\I,\{i\}})_*(\gff\otimes \O_{X^{\I\cup i}})$,\footnote{In this formula, 
$(p_{\I}\circ j_{\I,\{i\}})_*$ denotes the direct image 
in the category of $\O$-modules.}
where $p_\I$ denotes the projection $X^{\I\cup i}\to X^\I$.

\smallskip

For a surjection of finite sets $\phi:\I\surj \J$ we have a natural map 
$$X^{\J}\underset{X^{\I}}\times\tMaps{^\I}(X,G)\to
\tMaps{^\J}(X,G)$$ and for a decomposition $\I=\I_1\cup\ldots\cup\I_k$, we have a 
map
$$\tMaps{^\I}(X,G)|(X^{\I_i\neq\I_j})\to \tMaps{^{\I_1}}(X,G)\times\ldots\times 
\tMaps{^{\I_k}}(X,G)|(X^{\I_i\neq\I_j}).$$

\medskip

Observe now that there is a natural action 
$$\Act^\I:\tMaps{^\I}(X,G)\underset{X^\I}\times \Gr^\I_G(X)\to \Gr^\I_G(X):$$

Take 
$$(x^\I\in X^\I, g:X\setminus x^\I\to G)\in\tMaps{^\I}(X,G) \text{ and }
(x^\I,\P,\alpha:\P\simeq\P_0|_{X\setminus x^\I}\in \Gr^\I_G(X).$$
To define the corresponding new point $\Act^\I((x^\I,g),(x^\I,\P,\alpha))$ of 
$\Gr^\I_G(X)$ we need to specify a new
$G$-torsor $\P'$ with a trivialization $\alpha'$ over $X\setminus x^\I$.

\smallskip

The $G$-torsor $\P'$ is set to be isomorphic to $\P$ over $X\setminus x^\I$ 
with
$\alpha':\P'\to\P_0|_{X\setminus x^\I}$ being obtained from $\alpha$. A local 
section $s:X\setminus x^\I\to\P_0$
is now set to give rise to a regular section of $\P'$ (via $\alpha'$) if 
$\alpha^{-1}\circ g\circ s:X\setminus x^\I\to\P$ is regular.

\medskip

It is easy to see, moreover, that the action 
$$\on{Act}^\I:\tMaps{^\I}(X,G)\underset{X^\I}\times \Gr^\I_G(X)\to \Gr^\I_G(X)$$
described above is compatible with the connection along $X^\I$
and with the factorization isomorphisms. Hence, the sheaves of Lie algebras 
$\on{Lie}(\tMaps{^\I}(X,G))$ act on the 
D-modules $\A^{(\I)}_{\Gr_G}$ 
in a way compatible with the factorization
isomorphisms. In particular, the chiral bracket
$j_*j^*(\A_{\Gr_G}\boxtimes \A_{\Gr_G})\to \Delta_{!}(\A_{\Gr_G})$
is compatible with the action of 
$\on{Lie}(\tMaps{^{\{1,2\}}}(X,G))$ in the sense that the diagram:
$$
\CD
\on{Lie}(\tMaps{^{\{1,2\}}}(X,G))\otimes j_*j^*(\A_{\Gr_G}\boxtimes \A_{\Gr_G}) 
@>{\on{Act}^{\{1,2\}}}>> 
j_*j^*(\A_{\Gr_G}\boxtimes \A_{\Gr_G}) \\
@V{\on{id}\otimes \{,\}}VV  @VVV  \\
\Delta_{!}(\Delta^{!}(\on{Lie}(\tMaps{^{\{1,2\}}}(X,G)))[1]\otimes \A_{\Gr_G})  \\
@VVV    @VVV  \\
\Delta_{!}(\on{Lie}(\tMaps{^{\{1\}}}(X,G))\otimes \A_{\Gr_G}) 
@>{\on{Act}^{\{1\}}}>> \Delta_{!}(\A_{\Gr_G})
\endCD
$$
commutes.

\smallskip

In the proof of \propref{basickm} we shall use the following observation:

\smallskip

Let $\on{Lie}(\tMaps{^{\{1\}}}(X,G))_x$ denote the fiber of the sheaf 
$\on{Lie}(\tMaps{^{\{1\}}}(X,G))$ at $x\in X$. 

\begin{lem} \label{twoactions}
The action $\on{Act}^{\{1\}}:\on{Lie}(\tMaps{^{\{1\}}}(X,G))_x\otimes 
A_{\Gr_G}{}_x\to A_{\Gr_G}{}_x$
coincides with the one that comes from the
embedding $$\on{Lie}(\tMaps{^{\{1\}}}(X,G))_x\simeq \gff\otimes\O(X\setminus 
x)\hookrightarrow \gff\otimes\hat\K_x$$
and the canonical action $\on{Act}_x:(\gff\otimes\hat\K_x)\otimes A_{\Gr_G}{}_x\to 
A_{\Gr_G}{}_x$.
\end{lem}

\begin{proof} (Of \propref{basickm})

Consider the completion $p_1{}_*(\gff\otimes\O\hat\boxtimes\O)$ as a left D-module 
on $X$.

\smallskip

To specify a map $j_*j^*(\B(\gff,0)\boxtimes\A_{\Gr_g})\to \Delta_{!}(\A_{\Gr_g})$
is the same as to specify a continuous map 
$p_1{}_*(\gff\otimes\O\hat\boxtimes\O)\otimes \A_{\Gr_g}\to \A_{\Gr_g}$
compatible with the right D-module structure.

\smallskip

Without restricting the generality, we can assume that $X$ is affine. In this case 
$\on{Lie}(\tMaps{^{\{1\}}}(X,G))\simeq p_1{}_*(\gff\otimes\O\boxtimes\O)$
is dense in $p_1{}_*(\gff\otimes\O\hat\boxtimes\O)$ and it is enough to construct 
a continuous map
$$p_1{}_*(\gff\otimes\O\boxtimes\O)\otimes \A_{\Gr_g}\to \A_{\Gr_g}.$$

\smallskip

However, such a map has been already constructed: this is the map 
$\on{Act}^{\{1\}}$. The fact that property (a) of
\propref{basickm} holds follows from the fact that the map $\on{Act}^{\{1,2\}}$ 
commutes with the chiral bracket 
(cf. the commutative diagram above). Property (b) of \propref{basickm} follows 
from \lemref{twoactions}.
\end{proof}

\newpage

\centerline{\bf Chapter VI. The Free Bosonic Theory}
\ssec{The canonical line bundle}
 Let $H$ be a torus, let $\Lambda$ denote the lattice of $1$-parameter subgroups 
in $H$ and let
$\hf$ be the Lie algebra of $H$. Consider the affine Grassmannians
$\Gr_H$, $\Gr_H(X)$ and $\Gr^\I_H(X)$ corresponding to the group $H$.

\smallskip

The reduced schemes $_{red}\Gr_H$ and $_{red}\Gr_H(X)$ are finite dimensional and 
are identified with $\Lambda$ and
$\Lambda\times X$, respectively. For $\lambda\in\Lambda$ we shall denote by 
$\Gr_H^\lambda$ (resp., $\Gr_H^\lambda(X)$)
the corresponding connected component of $\Gr_H$ (resp., of $\Gr_H(X)$); we shall 
denote by $\unit_{\Gr}^\lambda$
the canonical section $X\simeq {}_{red}\Gr_H^\lambda \hookrightarrow 
\Gr_H^\lambda$.

\smallskip

As for $\Gr^\I_H(X)$, its set of connected (resp., irreducible) 
components is identified with $\Lambda$ (resp., with
$\Lambda^\I$). The reduced scheme corresponding to each irreducible component 
projects isomorphically onto
$X^\I$. We leave it to the reader to work out the intersection pattern of various 
irreducible components of
$_{red}\Gr^\I_H(X)$. For $\overline\lambda\in\Lambda^\I$ we shall denote by 
$\unit_{\Gr}^{\overline\lambda}$
the corresponding section $X^\I\to {}_{red}\Gr^\I_H(X)\to \Gr^\I_H(X)$.

\sssec{} \label{groupextension}
 Recall the groups $H(\hat\O)$ and $H(\hat\K)$ of \secref{intrgrass}. Now let 
now $Q$ be an integral-valued even symmetric form on $\Lambda$. We shall denote by 
the same 
character $Q$ the corresponding quadratic form on $\hf$. In addition, we shall 
choose once and for all a $2$-cocycle
$\epsilon$ of the group $\Lambda$ with values in $\ZZ_2$ whose class in 
$H^2(\Lambda,\ZZ_2)$ corresponds to $Q\on{mod}2$. 

\smallskip

It is known (cf. \cite{De}, \cite{BP}) that to the data of $(\Lambda,Q,\epsilon)$ 
one can associate in a canonical way
a central extension $H(\hat\K)'$ of the ind-algebraic group $H(\hat\K)$:
$$1\to \CC^*\to H(\hat\K)'\to H(\hat\K)\to 1.$$

\smallskip

The group $H(\hat\K)'$ has the following properties:

\begin{itemize}

\item
We have a canonical lifting of the embedding $H(\hat\O)\to H(\hat\K)$ to an 
embedding $H(\hat\O)\to H(\hat\K)'$.
In particular, we obtain an $H(\hat\K)'$-equivariant line bundle $R_Q$ on 
$\Gr_H=H(\hat\K)/H(\hat\O)$.

\item
The Lie algebra of $H(\hat\K)'$ can be identified with 
the Heisenberg algebra
$\on{Heis}(\hf,Q):=\hf\otimes \hat\K\oplus \CC$.

\item
For a complete curve curve $X$ and $n$ points $x_1,\dotsc,x_n\in X$ consider the 
corresponding central extension
$$1\to\CC^*\to (H(\hat\K_{x_1})\times\ldots\times H(\hat\K_{x_n}))'
\to H(\hat\K_{x_1})\times\ldots\times H(\hat\K_{x_n})\to 1.$$
We have a canonical lifting of the embedding 
$$\on{Maps}(X\setminus \{x_1,\dotsc,x_n\},H)\subset 
H(\hat\K_{x_1})\times\ldots\times H(\hat\K_{x_n})$$
to a homomorphism of formal group schemes $\on{Maps}(X\setminus 
\{x_1,\dotsc,x_n\},H)\to H(\hat\K)_n'$.

\end{itemize}

\medskip

To state the additional property of $H(\hat\K)'$ we need to introduce some 
notation. 

\smallskip

For a character $\mu\in\hf^*$, let $V^{\mu}$ denote the Weyl module over 
$\on{Heis}(\hf,Q)$, i.e. 
$$V^{\mu}\simeq\on{Ind}^{\on{Heis}(\hf,Q)}_{\hf\otimes\hat\O\oplus\CC}(\CC^{\mu}),
$$
where $\CC^\mu$ denotes the $1$-dimensional representation of 
$\hf\otimes\hat\O\oplus\CC$ corresponding
to $\mu:\hf\to\CC$. When $\mu=Q(\lambda,\cdot)$ for 
$\lambda\in\Lambda\subset\Lambda\underset{\ZZ}\otimes\CC\simeq\hf$, 
we endow $V^{\mu}$ with an action of the group $\Aut^+_0$ by 
requiring
that it be compatible with the $\Aut^+_0$-action on $\on{Heis}(\hf,Q)$ and by 
setting
$\CC^{\mu}\simeq L_{-Q(\lambda,\lambda)/2}$ as an $\Aut^+_0$-module (cf. 
\secref{linalg}). 

\smallskip

Now note that the vector space $\Gamma(\Gr^\lambda_H,R_Q)$ carries an 
$\Aut^+_0$-action (this is due to the fact that
the construction of $H(\hat\K)'$ and hence of $R_Q$ is canonical).

\smallskip

The fourth property of $H(\hat\K)'$ reads as follows:

\begin{itemize}

\item 

We have an isomorphism of $\on{Heis}(\hf,Q)$-modules: 
$\Gamma(\Gr^\lambda_H,(R_Q)^{-1})^*\simeq V^{\mu}$,
where $\mu=Q(\lambda,\cdot)$. This isomorphism is compatible with the 
$\Aut^+_0$-action.

\end{itemize}

\sssec{} \label{linebundle}
 Globally over a (not necessarily complete) curve $X$, the form $Q$ gives rise to 
a central extension 
$$1\to\CC^*\to H'_\K(X)\to H_\K(X)\to 1$$ of the ind-group scheme $H_\K(X)$. 
Let $R_Q(X)$ denote the line bundle $$R_Q(X):=H'_\K(X)/H_\O(X)$$ over 
$\Gr_H(X)=H_\K(X)/H_\O(X)$. Property (2) of the extension $H(\hat\K)'$ implies 
that
$R_Q(X)$ is equivariant with respect to the action of the group scheme $H_\O(X)$ 
on $\Gr_H(X)$.

\smallskip

Property (4) of the group $H(\hat\K)'$ implies the following assertion:

\begin{lem} \label{hw}
We have a canonical isomorphism 
$$\unit_{\Gr}^\lambda{}^*(R_Q(X))\simeq \Omega^{-Q(\lambda,\lambda)/2}.$$
\end{lem}

\medskip

A generalization of the above construction yields a line bundle $R^\I_Q(X)$ over 
the pre-image in $\Gr^\I_H(X)$ of 
$X^\I\setminus\Delta$. Our present goal is to extend this line bundle to the whole 
of $\Gr^\I_H(X)$, i.e.
to the locus where the base points $x^\I_{i_1},\dotsc,x^\I_{i_n}\in X$ collide 
with one another.

\smallskip

Assume now that $X$ is complete. Property (3) of the extension $H(\hat\K)'$ 
implies the following assertion:

\begin{lem}
For a complete curve $X$ there exists a canonical line bundle $\R_Q(X)$ over 
$\Bun_H(X)$ such that for every $\I$, 
the restriction of the line bundle $k^\I{}^*(\R_Q(X))$ to $\Gr^\I_H(X)\setminus 
(r^\I)^{-1}(\Delta)$
is identified with $R_Q^\I$.
\end{lem}

\smallskip

We now set $R^\I_Q(X)$ to be the line bundle over $\Gr^\I_H(X)$ equal to 
$k^\I{}^*(\R_Q(X))$.
(When $\I=\{1,\dotsc,n\}$, we shall replace the notation $R^\I_Q(X)$ simply by 
$R^n_Q(X)$.)

\smallskip

Since the projection $k^\I:\Gr^\I_H(X)\to \Bun_H(X)$ respects the connection on 
$\Gr^\I_H(X)$ along $X^\I$,
we obtain a connection along $X^\I$ on the line bundle $R^\I(Q)$.

\medskip

\begin{prop} \label{factorline}
Under the isomorphisms of \propref{factorgr}, the system of line bundles $\I\to 
R^\I_Q(X)$ satisfies the following
factorization property: 

\smallskip

\noindent {\em (a)} For a surjection of finite sets $\phi:\I\surj \J$ the first 
isomorphism of \propref{factorgr}
underlies an isomorphism of line bundles: 
$R^\I_Q(X)|_{\Gr_H^\J(X)}\simeq R^\J_Q(X)$.

\smallskip

\noindent {\em (b)} For a decomposition $\I=\I_1\cup\ldots\cup \I_k$ the second 
isomorphism of \propref{factorgr}
underlies an isomorphism of line bundles:
\begin{align*}
&R^\I_Q(X)|_{X^{\I_i\neq \I_j}\underset{X^\I}\times\Gr^\I_H(X)}\simeq \\
&\simeq R_Q^{\I_1}(X)\times\ldots\times R_Q^{\I_1}(X)|_
{X^{\I_i\neq \I_j}\underset{X^{\I_1}\times\ldots\times 
X^{\I_k}}\times\Gr^{\I_1}_H(X)\times\ldots\times\Gr^{\I_k}_H(X)}.
\end{align*}

\smallskip

Moreover, the above isomorphisms are compatible with the connections on the line 
bundles $R^\I_Q(X)$'s along the 
$X^\I$'s.

\end{prop}

\begin{proof}

The isomorphism of point (a) follows from the definitions. To prove point (b) note 
that the two line bundles are
obviously isomorphic over the pre-image of $X^\I\setminus\Delta$ in $\Gr^\I_H(X)$ 
and we only have to show that
this isomorphism extends to the whole of $\Gr^\I_H(X)$. 
This follows 
easily (by degree considerations) from the following assertion:

\begin{lem}
For ${\overline\lambda}=\lambda_1,\dotsc,\lambda_n$ the line bundle 
$\unit_{\Gr}^{{\overline\lambda}}{}^*(R^n_Q(X))$
over $X^n$ is identified canonically with
$$\Omega^{-Q(\lambda_1,\lambda_1)/2}\boxtimes\ldots\boxtimes 
\Omega^{-Q(\lambda_n,\lambda_n)/2}
(\underset{i<j}\Sigma Q(\lambda_i,\lambda_j)\cdot\Delta_{x_i=x_j}).$$
\end{lem}

\end{proof}

\sssec{}
 We now want to show that the dependence of the line bundles $R^\I_Q(X)$ on $X$ is 
local in the sense of
\secref{intrgrass}. This fact is not immediately obvious, since in the definition 
of $R^\I_Q(X)$ we gave above
we used the fact that $X$ is complete. A way to overcome this difficulty will be 
via the ind-schemes $\tGr^n_H(X)$
introduced in \secref{twgrassmannian}.

\medskip

Using the iterative description of the ind-scheme $\tGr^n_H(X)$ given in 
\secref{twgrassmannian} 
we can produce a line bundle ${\widetilde R_Q^n}(X)$ over $\tGr^n_H(X)$. This 
construction works for curves that are
not necessarily complete; therefore we can form 
line bundle ${\widetilde R^n_Q}(\Xf)$ 
over $\tGr^n_H(\Xf)$ which is
local in the sense of \secref{intrgrass}. In more detail, let 
$X^{S,I},x^{S,I},\phi^{S,I}$ be as in \secref{intrgrass} and
let $\tGr^n_H(X^S)\setminus x^S$ denote the open subscheme of 
$\tGr^n_H(X^S)$ equal to the pre-image
in $\tGr^n_H(X^S)$ of 
$$(X^S\setminus x^S)\underset{S}\times\ldots
\underset{S}\times (X^S\setminus x^S)\subset 
X\underset{S}\times\ldots\underset{S}\times X$$
under the projection $r^n\circ w_n$.

\smallskip

Then the isomorphism of ind-schemes 
$$\tGr^n_H(X^{S\times I})\setminus x^{S\times I}\simeq 
(\tGr^n_H(X^S)\setminus x^S)\times I$$
underlies an isomorphism of line bundles:
$${\widetilde R_Q^n}(X^{S\times I})
|_{\tGr^n_H(X^{S\times I})\setminus x^{S\times I}}\simeq
({\widetilde R_Q^n}(X^S)|_{\tGr^n_H(X^S)\setminus x^S})\times I.$$

\smallskip

Note that over the pre-image of $X^n\setminus\Delta$
in $\tGr^n_H(X)$ (which is identified with the pre-image of $X^n\setminus\Delta$ 
in $\Gr^n_H(X)$) the line bundle 
${\widetilde R_Q^n}(X)$ is identified naturally with $R^n_Q(X)$.

\begin{prop}
The line bundle $R^n_Q(X)$ is uniquely characterized by the following property: 

There exists a canonical isomorphism $$w_n^*(R_Q^n)\simeq {\widetilde R_Q^n}$$
which reduces to the identity isomorphism over $X^n\setminus\Delta$.
\end{prop}

The proof of this assertion is not difficult to deduce from the construction of 
the line bundle ${\widetilde R_Q^n}(X)$
and we will omit it.

\medskip

\begin{cor} \label{loclinebundle}
The line bundle $R^\I_Q(\Xf)$ over $\Gr^\I_H(\Xf)$ is local (in the above sense). 
Moreover, the factorization isomorphisms 
of \lemref{factorline} are compatible with the local structure.
\end{cor}

\ssec{Construction of the bosonic chiral algebra}
\sssec{}  \label{definebose}
 For $(H,Q)$ as above consider the chiral algebra $\A(\hf,Q)$ of 
Example 1 of \secref{KM}. 
Our goal in this subsection is to construct a bigger chiral algebra $\A(H,Q)$ that
will contain $\A(\hf,Q)$ as a subalgebra.

\smallskip

Recall that if $Z$ is a strict ind-scheme and $R$ is a line bundle over $Z$ it 
makes sense to talk about $R$-twisted
right D-modules on $Z$ which are set-theoretically supported on finite-dimensional 
subschemes of $Z$.

\smallskip

Thus, consider the $R_Q(X)$-twisted right D-module 
$$\unit_{\Gr}^\lambda{}_{!}(\Omega)\otimes R_Q(X)$$
on $\Gr_H(X)$. Let $\A(H,Q)^\lambda$ denote its $\O$-module
direct image onto $X$, i.e.
$$\A(H,Q)^\lambda:=r_*(\unit_{\Gr}^\lambda{}_{!}(\Omega)\otimes R_Q(X)).$$

\smallskip

The connection on the pair $(\Gr_H(X),R_Q(X))$ along $X$ defines on 
$\A(H,Q)^\lambda$
a structure of right D-module on $X$. Finally, let $\A(H,Q)$ be the direct sum 
$\underset{\lambda\in\Lambda}\oplus \A(H,Q)^\lambda$. The sections 
$\unit_{\Gr}^\lambda:X\to \Gr_H(X)^\lambda$
give rise to $\O$-module maps 
$$\unit^\lambda:\Omega^{-Q(\lambda,\lambda)/2+1}\hookrightarrow \A(H,Q);$$ 
for $\lambda=0$ the map $$\unit:\Omega\hookrightarrow \A(H,Q)^0\hookrightarrow 
\A(H,Q)$$
is in fact a map of D-modules.

\smallskip

More generally, for a finite set $\I$ and for an element 
$\overline\lambda\in\Lambda^\I$ consider 
on $X^\I$ the right D-module
$$\A(H,Q)'{}^{\overline\lambda}:=r^\I_*(\unit_{\Gr}^{\overline\lambda}{}_{!}
(\Omega_{X^\I})\otimes R^\I_Q(X)).$$
Let $\A(H,Q)'{}^{(\I)}$ denote the direct sum 
$\underset{\overline\lambda\in\Lambda^\I}\oplus \A(H,Q)'{}^{\overline\lambda}$ and 
let
$^l\A(H,Q)'{}^{(\I)}$ be the corresponding left D-module on $X^\I$.

\medskip

Consider the sheaf $r^\I_*((R^\I_Q(X))^{-1})$ over $X^\I$. Since $\Gr_H^\I(X)$ is 
an inductive limit of schemes finite over $X^\I$,
the above direct image is a projective limit of coherent sheaves on $X^\I$; 
moreover, it carries a natural left D-module structure
due to the connection on $R^\I_Q(X)$ along $X^\I$. It is easy to see that we have 
in fact a canonical isomorphism of left D-modules:
$$^l\A(H,Q)'{}^{(\I)}\simeq 
\on{Hom}_{\O_{X^\I}}(r^\I_*((R^\I_Q(X))^{-1}),\O_{X^\I}).$$
 
\medskip

The following assertion follows from \propref{factorline}:

\begin{lem}
The system of $\O$-modules $\I\to {}^l\A(H,Q)'{}^{(\I)}$ satisfies the conditions 
of \thmref{reconstr}.
\end{lem}

\smallskip

Thus, $\A(H,Q)$ has a natural structure of chiral algebra such that for every $\I$
$$^l\A(H,Q)^{(\I)}\simeq {}^l\A(H,Q)'{}^{(\I)},$$ as $\O$-modules. The last 
assertion of \propref{factorline} implies that
the above isomorphisms preserve the D-module structure as well.

\medskip

For $x\in X$, let $\Gr_H(X)_x$ denote the fiber of $\Gr_H(X)$ over $x$ and let 
$R_Q(X)_x$ denote the restriction
of $R_Q(X)$ to this fiber. (Once we choose an isomorphism $\hat\O\simeq\hat\O_x$, 
the pair $(\Gr_H(X)_x,R_Q(X)_x)$ can
be identified with $(\Gr_H,R_Q)$.). By the construction, the fiber 
$A(H,Q)^\lambda_x$ of $\A(H,Q)^\lambda$ at $x$ 
is identified with 
$$
\Gamma(\Gr^\lambda_H(X)_x,(R_Q(X)_x)^{-1})^*.
$$

\smallskip

Therefore, the Harish-Chandra pair $(\hf\otimes\hat\K_x\oplus\CC,H(\hat\O_x))$ 
acts on
$A(H,Q)_x$; we shall denote by
$$\Act:(\hf\otimes\hat\K_x\oplus\CC)\otimes A(H,Q)_x\to A(H,Q)_x$$
the corresponding action map. 

\smallskip

Property (4) from \secref{groupextension}
implies that as a $(\hf\otimes\hat\K_x\oplus\CC,H(\hat\O_x))$-module, 
$A(H,Q)_x^\lambda$ is canonically
isomorphic to the Weyl module $V^{\mu}$ with $\mu=Q(\lambda,\cdot)$.

\smallskip

\begin{thm} \label{descrbose}
We have:

\smallskip

\noindent{\em (a)}
The chiral bracket on $\A(H,Q)$ preserves the $\Lambda$-grading, i.e.
$$j_*j^*(\A(H,Q)^\lambda\boxtimes 
\A(H,Q)^\mu)\to\Delta_{!}(\A(H,Q)^{\lambda+\mu}).$$ 
In particular, $\A(H,Q)^0$ is a chiral subalgebra in $\A(H,Q)$ and each 
$\A(H,Q)^\lambda$ is a chiral module over $\A(H,Q)^0$.

\smallskip

\noindent{\em (b)}
There exists a canonical isomorphism of chiral algebras 
$\A(H,Q)^0\simeq\A(\hf,Q)$.

\smallskip

\noindent{\em (c)}
For $x\in X$, the map $(\hf\otimes\hat\K_x\oplus\CC)\otimes A(H,Q)_x\to A(H,Q)_x$ 
induced by the chiral action of $\B(\hf,Q)$ on $\A(H,Q)$ coincides with the map 
$\Act$ (cf. above).

\end{thm}

The proof of this theorem goes along the same lines as the proof of 
\thmref{kmisom} and we will omit it.

\sssec{}
 We shall now show that the chiral algebra $\A(H,Q)$ defines a CFT of central 
charge $c=\on{dim}(H)$.

\smallskip

It follows from \corref{loclinebundle} that the $\O_{X^\I}$-modules 
$^l\A(H,Q)'{}^{(\I)}$ carry a structure of 
local $\O$-modules over $\Xf^\I$. Moreover, the factorization isomorphisms are 
morphisms of local $\O$-modules. 
Hence, according to \propref{reconstr}, the chiral algebra $\A(H,Q)$ satisfies 
condition (a) of \secref{central} and 
it remains to endow it with an energy-momentum tensor such that condition (b) of 
\secref{central} holds.

\smallskip

Consider the map $T_{H,Q}:\Theta'_{\on{dim}(\hf)}\to \A(H,Q)$ obtained as a 
composition:
$$\Theta_{\on{dim}(\hf)}\overset{T_{\hf,Q}}\longrightarrow \A(\hf,Q)\simeq 
\A(H,Q)^0\subset \A(H,Q).$$

\begin{prop}
The map $T_{H,Q}$ satisfies condition (b) of \secref{central}.
\end{prop}

\begin{proof}

For a vector field $\xi$ on a curve $X$ consider the two maps 
$\A(H,Q)^\lambda\to\A(H,Q)^\lambda$ given by
$l\to l\cdot\xi+\on{Lie}_{\xi}(l)$ and $(\on{id}\boxtimes 
h)(\{l,T_{H,Q}(\xi')\})$, where $\xi'$ is some lifting
of $\xi$ to a section of $\Theta'_{\on{dim}(\hf)}$.
Both these maps are endomorphisms of the D-module structure on $\A(H,Q)^\lambda$; 
therefore, they induce endomorphisms
of each fiber $A(H,Q)_x^\lambda$ for $x\in X$. 

\smallskip

Now choose an identification $\hat\O_x\simeq\hat\O$. Since 
$T_{\hf,Q}$ is an energy-momentum tensor for $\A(\hf,Q)$, the above construction 
yields two actions of the Lie algebra
$\Vir^+$ on $V^{\mu}$ (here $\mu=Q(\lambda,\cdot)$) that are compatible with the 
$\Vir^+$-action on 
$\on{Heis}(\hf,Q)$. 

\smallskip

Moreover, we claim that the restrictions of these two actions to $\Vir^+_0\subset 
\Vir^+$ coincide. 
Indeed, it is enough to check this fact on the generating space $\CC^\mu\subset 
V^{\mu}$ and in both cases 
$\Vir^+_0$ acts on $\CC^\mu$ as on $L_{-Q(\lambda,\lambda)/2)}$ (this is an easy 
computation using 
\lemref{hw} and \propref{normalordering}).

\smallskip

The assertion of the proposition now follows from the following general fact:

\begin{lem}
There is at most one extension of the $\Vir^+_0$-action on $V^{\mu}$ to an action 
of $\Vir^+$ which is compatible
with the action of the latter Lie algebra on $\on{Heis}(\hf,Q)$.
\end{lem}

\end{proof}

\sssec{}
 Assume now that $X$ is complete. Our present goal is to compute the space of 
conformal blocks of the chiral algebra 
$\A(H,Q)$. 

\smallskip 

\begin{thm} \label{confbose}
For a complete curve $X$ we have a canonical isomorphism
$$\conf(X,\A(H,Q))\simeq H^0(\Bun_H(X),(\R_Q(X))^{-1})^*.$$
\end{thm}

Before proving this theorem let us make a few observations:

\smallskip

\noindent{\it Observation 1:}

The stack $\Bun_H(X)$ splits into connected components 
$$
\Bun_H(X)=\underset{\lambda\in\Lambda}\cup
\Bun_H(X)^\lambda 
$$
and each $\Bun_H(X)^\lambda$ is isomorphic to the quotient of the Jacobian of $X$ 
by the trivial
action of the group $H$; this corresponds to the fact that $\forall 
\P\in\Bun_H(X)$, $\on{\Aut}(\P)\simeq H$.
 
\smallskip

\noindent{\it Observation 2:}

It is easy to infer from the properties of the extension $H(\hat\K)'$ in 
\secref{groupextension} that for every
$\P\in\Bun_H(X)^\lambda$ the group $H$ of automorphisms of $\P$ acts on the fiber 
of $\R_Q(X)$ at $\P$ via
the character $Q(\lambda,\cdot)$. This implies, in particular, that the 
contribution to $H^0(\Bun_H(X),(\R_Q(X))^{-1})$
comes only from the $0$-th component of $\Bun_H(X)$. This shows that the RHS of 
the isomorphism of \thmref{confbose}
is finite-dimensional.

\smallskip

It is equally easy to see that if $x_1,\dotsc,x_n$ are distinct 
points of $X$, then the canonical map
$A(H,Q)_{x_1}^{\lambda_1}\otimes\ldots\otimes 
A(H,Q)_{x_n}^{\lambda_n}\to \conf(X,\A(H,Q))$ is non-zero only
if $\lambda_1+\ldots +\lambda_n=0$. (This is due to the fact that the subspace of 
``constant currents''
$$\hf\subset DR^0(X,\A(\hf,Q))\subset DR^0(X,\A(H,Q))$$ acts on each 
$A(H,Q)_{x_i}^{\lambda_i}\simeq V^{\mu_i}$
by the character $\mu_i$.)

\smallskip

\noindent{\it Observation 3:}
The assertion of the theorem can be reformulated as an isomorphism between
$H^0(\Bun_H(X),(\R_Q(X))^{-1})$ and $\conf(X,\A(H,Q))^*$.

\smallskip

According to \propref{confran}, $$\conf(X,\A(H,Q))^*\simeq 
\on{Hom}_{D(X^n)}({}^l\A(H,Q)^{(n)},\O_{X^n})$$
for any $n\geq 2$. Since
$$^l\A(H,Q)^{(n)}\simeq \on{Hom}_{\O_{X^n}}(r^n_*((R^n_Q(X))^{-1}),\O_{X^n}),$$ we 
obtain that the space
$\conf(X,\A(H,Q))^*$ is identified with the space 
$$\Gamma(\Gr^n_H(X),(R_Q^n(X))^{-1})^{D_{X^n}}$$
of $D_{X^n}$-invariant sections of $(R^n_Q(X))^{-1}$ for $n$ as above.

\medskip

Let us write out the above isomorphism 
$$\conf(X,\A(H,Q))^*\simeq \Gamma(\Gr^n_H(X),(R_Q^n(X))^{-1})^{D_{X^n}}$$
more explicitly:

\smallskip

Let $\chi$ be an element of $\conf(X,\A(H,Q))^*$. Let us on the one hand consider 
the composition
\begin{align*}
&\Omega^{-Q(\lambda_1,\lambda_1)/2}\boxtimes\ldots\boxtimes 
\Omega^{-Q(\lambda_n,\lambda_n)/2}\overset
{\unit^{\lambda_1}\boxtimes\ldots\boxtimes \unit^{\lambda_n}}\longrightarrow 
{}^l\A(H,Q)^{\boxtimes n}\to \\
&j_*j^*(\O_{X^n})\otimes \conf(X,\A(H,Q)) \overset{\chi}\longrightarrow 
j_*j^*(\O_{X^n}).
\end{align*}
The above map can be regarded as a section of 
$\Omega^{Q(\lambda_1,\lambda_1)/2}\boxtimes\ldots\boxtimes 
\Omega^{Q(\lambda_n,\lambda_n)/2}$ over $X^n\setminus\Delta$.

\smallskip

On the other hand, an element $\chi'\in 
\Gamma(\Gr^n_H(X),(R_Q^n(X))^{-1})^{D_{X^n}}$
can be regarded as a map from $\Gr_H^n(X)$ to the total space of 
the line bundle 
$(\R_Q(X))^{-1}$ over $\Bun_H(X)$.
According to \lemref{hw}, the pull-back of the line bundle 
$(\R_Q(X))^{-1}$ under 
the composition
$$k^n\circ \unit_{\Gr}^{{\overline\lambda}}\circ 
j:X^n\setminus\Delta\hookrightarrow X^n\to \Bun_H(X)$$
is the line bundle 
$\Omega^{Q(\lambda_1,\lambda_1)/2}\boxtimes\ldots\boxtimes 
\Omega^{Q(\lambda_n,\lambda_n)/2}$ over $X^n\setminus\Delta$. Therefore,
the composition
$$\chi'\circ\unit_{\Gr}^{{\overline\lambda}}\circ 
j:X^n\setminus\Delta\hookrightarrow X^n\to \on{Tot}(\R_Q(X)^{-1})$$
is again a section of 
$\Omega^{Q(\lambda_1,\lambda_1)/2}\boxtimes\ldots\boxtimes 
\Omega^{Q(\lambda_n,\lambda_n)/2}$ over $X^n\setminus\Delta$.

\smallskip

It follows from the definitions that the above two sections of 
$\Omega^{Q(\lambda_1,\lambda_1)/2}\boxtimes\ldots\boxtimes 
\Omega^{Q(\lambda_n,\lambda_n)/2}(\infty\cdot\Delta)$ coincide if the
elements $\chi$ and $\chi'$ correspond to one another under the isomorphism
$$\conf(X,\A(H,Q))^*\simeq \Gamma(\Gr^n_H(X),(R_Q^n(X))^{-1})^{D_{X^n}}.$$

\begin{proof} (Of \thmref{confbose})

By the very definition of the connection on the line bundle $R^n_Q(X)$, we have a 
map 
$$H^0(\Bun_H(X),(\R_Q(X))^{-1})\to \Gamma(\Gr^n_H(X),(R_Q^n(X))^{-1})$$ whose 
image belongs to
$\Gamma(\Gr^\I_H(X),(R_Q^n(X))^{-1})^{D_{X^n}}$. Moreover, the assertion of the 
above Observation 3 combined with 
\propref{confmany} imply that the composition
$$H^0(\Bun_H(X),(\R_Q(X))^{-1})\to \Gamma(\Gr^n_H(X),(R_Q^n(X))^{-1})^{D_{X^n}}\to 
\conf(X,\A(H,Q))^*$$
is the same for any $n$.

\smallskip

Moreover, the map $k^n$ is dominant for any $n$. Therefore, the map 
$$H^0(\Bun_H(X),(\R_Q(X))^{-1})\to \conf(X,\A(H,Q))^*$$
is injective and it remains to show that it is surjective. This will be done in 
two steps.

\medskip

\noindent{\bf Step 1}

\smallskip

Let 
$$\chi'\in \Gamma(\Gr^n_H(X),(R_Q^n(X))^{-1})^{D_{X^n}}\subset 
\Gamma(\Gr^n_H(X),k^n{}^*(\R_Q^n(X))^{-1})$$
be as above. We shall first show that it is locally constant along the fibers of 
the projection $k^n$.  

\smallskip
 
Let $(x_1,\dotsc,x_n,\P,\alpha:\P\to\P_0|X\setminus \{x_1,\dotsc,x_n\})$ be a 
$\CC$-point of $\Gr^n_H(X)$ and for a local
Artinian scheme $I$ let 
$(x^I_1,\dotsc,x^I_n,\P,\alpha^I:\P^I\to\P^I_0|X\setminus\{x^I_1,\dotsc,x^I_n\})$ 
be its extension to 
an $I$-valued point of $\Gr^n_H(X)$, such that $\P^I=\P\times I$.

\smallskip

We must show that the composition
$$I\to \Gr^n_H(X)\overset{\chi'}\longrightarrow \on{Tot}(\R_Q(X)^{-1})$$ is a 
constant map
$$I\to \Spec(\CC)\to \Gr^n_H(X)\overset{\chi'}\longrightarrow 
\on{Tot}(\R_Q(X)^{-1}).$$ 
Obviously, it is enough to check this property 
when $I=\Spec(\CC[t]/t^2)$.

\smallskip

First of all, using the connection along $X^n$ on $\Gr_H^n(X)$ we can reduce the 
assertion to the case when the $x_i$
are constant maps $I\to X$ and without loss of generality we can assume that the 
points $x_1,\dotsc,x_n$ are distinct.

\smallskip

In this case the difference between the trivializations 
$\alpha\times I$ and $\alpha^I$
of $\P^I\simeq \P\times I$ over $(X\setminus \{x_1,\dotsc,x_n\})\times I$ is given 
by an infinitesimal ``gauge transformation'', 
i.e. by an action on the fiber of $\Gr^n_H(X)$ over
$(x_1,\dotsc,x_n)$ by an element of
$$H^0(X\setminus\{x_1,\dotsc,x_n\}, \hf\otimes\O)\subset 
\hf\otimes\hat\K_{x_1}\oplus\ldots\oplus\hf\otimes\hat\K_{x_n}.$$

\smallskip

Therefore, in terms of the corresponding functional $\chi:\conf(X,\A(H,Q))\to\CC$, 
we must prove that the functional on
$A(H,Q)_{x_1}\otimes\ldots\otimes A(H,Q)_{x_n}$ obtained as a composition 
$$A(H,Q)_{x_1}\otimes\ldots\otimes A(H,Q)_{x_n}\to 
\conf(X,\A(H,Q))\overset{\chi}\to\CC$$
is $H^0(X\setminus\{x_1,\dotsc,x_n\},\hf\otimes\O)$-invariant. However, this 
follows from the definition of the space of conformal
blocks in view of point (c) of \thmref{descrbose}.

\medskip

\noindent{\bf Step 2} 

\smallskip

Let us view $\chi'$ as a map from $\Gr_H^n(X)$ to the total space of 
the line bundle $(\R_Q(X))^{-1}$ (cf. Observation 3).
To prove the theorem it remains to show that for two $\CC$-points
$(x_1,\dotsc,x_n,\P,\alpha)$ 
and $(x'_1,\dotsc,x'_n,\P',\alpha')$ 
of $\Gr^n_H(X)$ with $\P=\P'$, the map $\chi'$ has the same value on 
these points.

\smallskip

Let now $\nu_1,\dotsc,\nu_k$ be a basis of $\Lambda$, let $m$ be an integer 
satisfying $m>2g(X)-2$ and let
$n=2m\cdot k$. Take an element $\overline\lambda(m)\in\Lambda^n$ equal to
$$\overline\lambda(m)=\underbrace{\lambda_1,\dotsc,\lambda_1}_{\text {$m$ 
times}},\dotsc,
\underbrace{-\lambda_1,\dotsc,-\lambda_1}_{\text {$m$ times}},\dotsc,
\underbrace{\lambda_k,\dotsc,\lambda_k}_{\text {$m$ times}},\dotsc,
\underbrace{-\lambda_k,\dotsc,-\lambda_k}_{\text {$m$ times}}.$$

\smallskip

First of all, using Observations 2 and 3 and by increasing $n$ if necessary,
we can arrange that our two points $(x_1,\dotsc,x_n,\P,\alpha)$ 
and $(x'_1,\dotsc,x'_n,\P',\alpha')$  belong to the same
irreducible component $\Gr^{\overline\lambda(m)}_H(X)$, where $m$ and 
$\overline\lambda(m)$ are as above.
Recall that $\unit^{{\overline\lambda}}_{\Gr}(X):
X^n\to\Gr^n_H(X)$ defines an isomorphism between $X^n$ and the 
reduced scheme of the corresponding irreducible component of $\Gr^n_H(X)$.

\medskip

The composition 
$$X^n\hookrightarrow 
{}_{red}\Gr^n_H(X)\hookrightarrow\Gr^n_H(X)\overset{k^n}\to\Bun_H(X)$$
factors in this case via
$$X^n\to \underbrace{X^{(m)}\times\ldots\times X^{(m)}}_{\text {$2k$ times}}\to 
\underbrace{\on{Pic}^m(X)\times\ldots\times \on{Pic}^m(X)}_{\text {$2k$ 
times}}\simeq\Bun_H(X),$$
where $X^{(m)}$ denotes the $m$-th symmetric power of $X$ and the second 
arrow is the product of the Abel-Jacobi maps
$X^{(m)}\to \on{Pic}^m(X)$.

\smallskip

We claim that it is enough to show that the map 
$\chi':X^n\to\on{Tot}(\R_Q(X)^{-1})$ factors via
$$X^n\to X^{(m)}\times\ldots\times X^{(m)}\to \on{Tot}(\R_Q(X)^{-1}):$$

\smallskip

Indeed, the map 
$$X^{(m)}\times\ldots\times X^{(m)}\to \Bun_H(X)$$ is smooth, proper and has 
connected fibers, therefore any section of the pull-back
of $(\R_Q(X))^{-1}$ to $X^{(m)}\times\ldots\times X^{(m)}$ comes from a section of 
this line bundle over $\Bun_H(X)$.

\smallskip

To prove the required factorization property of $\chi'$, we can restrict our 
attention to $X^n\setminus\Delta\subset X^n$.
Using Observation 3, all we need to show is that the composition
\begin{align*}
&(\Omega^{-Q(\lambda_1,\lambda_1)/2})^{\boxtimes m}\boxtimes\ldots\boxtimes 
(\Omega^{-Q(\lambda_n,\lambda_n)/2})^{\boxtimes m}\to 
{}^l\A(H,Q)^{\boxtimes n}\to \\
&\to j_*j^*(\O_{X^n})\otimes \conf(X,\A(H,Q)) \overset{\chi}\longrightarrow 
j_*j^*(\O_{X^n})
\end{align*}
is invariant under the group $\underbrace{S^m\times\ldots\times S^m}_{\text{ $2k$ 
times}}$ of permutations. However, this follows 
from the fact that any correlation function is symmetric in the insertions (cf. 
\secref{confsummary}).

\end{proof}

\end{document}